\documentclass[final,onefignum,onetabnum,nohypdvips]{siamart190516}

\usepackage{lipsum}
\usepackage{amsfonts}
\usepackage{graphicx}
\usepackage{epstopdf}
\usepackage{algorithmic}
\ifpdf
  \DeclareGraphicsExtensions{.eps,.pdf,.png,.jpg}
\else
  \DeclareGraphicsExtensions{.eps}
\fi

\newsiamremark{remark}{Remark}
\newsiamremark{hypothesis}{Hypothesis}
\crefname{hypothesis}{Hypothesis}{Hypotheses}
\newsiamthm{claim}{Claim}
\newsiamremark{assumption}{Assumption}
\crefname{assumption}{Assumption}{Assumptions}

\headers{Transient Times for Distributed Stochastic Gradient Methods}{Kun Huang and Shi Pu}

\title{Improving the Transient Times for Distributed Stochastic Gradient Methods\thanks{This work was partially supported  by the Shenzhen
		Research Institute of Big Data (SRIBD) under Fund No. J00120190011 and by the National Natural Science Foundation of China (NSFC) under grant No. 62003287.}}

\author{Kun Huang\thanks{School of Data Science, Shenzhen Research Institute of Big Data, The Chinese University of Hong Kong, Shenzhen
  (\email{kunhuang@link.cuhk.edu.cn}, \email{pushi@cuhk.edu.cn}).}
\and Shi Pu\footnotemark[2]}

\usepackage{amsopn}
\DeclareMathOperator{\diag}{diag}
\DeclareMathOperator{\range}{range}

\newcommand{\E}{\mathbb{E}}
\newcommand{\R}{\mathbb{R}}

\newcommand{\x}{\mathbf{x}}

\newcommand{\y}{\mathbf{y}}
\newcommand{\s}{\mathbf{s}}
\newcommand{\g}{\mathbf{g}}
\newcommand{\0}{\mathbf{0}}
\newcommand{\1}{\mathbf{1}}
\newcommand{\tx}{\tilde{\mathbf{x}}}
\newcommand{\ty}{\tilde{\mathbf{y}}}
\newcommand{\T}{\top}
\newcommand{\mo}{\mathcal{O}}

\newcommand{\mF}{\mathcal{F}}
\newcommand{\mc}[1]{\mathcal{#1}}

\usepackage{array}
\usepackage{algorithm,algorithmic}
\usepackage{graphicx,epstopdf} 
\usepackage[caption=false]{subfig} 
\newcommand{\norm}[1]{\left\Vert #1 \right\Vert}
\newcommand{\mbb}[1]{\mathcal{O}\left( #1 \right)}

\ifpdf
\hypersetup{
	pdftitle={An Example Article},
	pdfauthor={D. Doe, P. T. Frank, and J. E. Smith}
}
\fi

\begin{document}
	
	\maketitle
	
	\begin{abstract}
		We consider the distributed optimization problem
		where $n$ agents each possessing a local cost function, collaboratively minimize the average of the $n$ cost functions  over a connected network. Assuming stochastic gradient information is available, we study a distributed stochastic gradient algorithm, called exact diffusion with adaptive stepsizes (EDAS) adapted from the Exact Diffusion method \cite{yuan2018exact} and NIDS \cite{li2019decentralized} and perform a non-asymptotic convergence analysis. We not only show that EDAS asymptotically achieves the same network independent convergence rate as centralized stochastic gradient descent (SGD) for minimizing strongly convex and smooth objective functions, but also characterize the transient time needed for the algorithm to approach the asymptotic convergence rate, which behaves as $K_T=\mo\left(\frac{n}{1-\lambda_2}\right)$,  where $1-\lambda_2$ stands for the spectral gap of the mixing matrix. To the best of our knowledge, EDAS achieves the shortest transient time when the average of the $n$ cost functions is strongly convex and each cost function is smooth.
		Numerical simulations further corroborate and strengthen the obtained theoretical results.
	\end{abstract}
	
	
	
	\section{Introduction}
	In this paper, we consider the distributed optimization problem, where a set of $n$ agents each possessing a local cost function $f_i$, seek to collaboratively solve the following optimization problem: 
	\begin{equation}\label{eq:eloss}
		\underset{x\in\mathbb{R}^p}{\min}f(x)\left(=\frac{1}{n} \sum_{i=1}^n f_i(x)\right),
	\end{equation}
	where  each $f_i:\R^p\rightarrow \R$ is assumed to have Lipschitz continuous gradients, and $f$ is strongly convex.
	The agents are connected over a network, where they may only communicate and exchange information with their direct neighbors. 
	
	To solve problem (\ref{eq:eloss}), we assume that at each iteration $k$ of the algorithm we study, each agent $i\in \mathcal{N}=\{1,2,\ldots,n\}$ is able to obtain a stochastic gradient sample of the form $g_i(x_{i,k},\xi_{i,k})$,  given input $x_{i,k}\in\R^p$, that satisfies the following condition:
	\begin{assumption}
		\label{ngrad}
		For all  $k\ge 0$, 
		each random vector $\xi_{i,k}\in\R^q$ is independent across $i\in\mathcal{N}$.
		Denote by $\mF_k$ the $\sigma$-algebra generated by $\{x_{i,0},x_{i,1},\ldots,x_{i,k}\mid i\in\mathcal{N}\}$. Then,
		\begin{equation}
			\label{condition: gradient samples}
			\begin{split}
				& \mathbb{E}[g_i(x_{i,k},\xi_{i,k})\mid \mF_{k}] =  \nabla f_i(x_{i,k}),\\
				& \mathbb{E}[\|g_i(x_{i,k},\xi_{i,k})-\nabla f_i(x_{i,k})\|^2\mid \mF_{k}] \le  \sigma_i^2 ,\quad\hbox{\ for some $\sigma_i>0$}.
			\end{split}
		\end{equation}
	\end{assumption}
	
	Stochastic gradients appear in many machine learning problems. For instance, suppose $f_i(x):=\E_{\xi_i}[F_i(x,\xi_i)]$ represents the expected loss function for agent $i$, where $\xi_i$ are independent data samples gathered continuously over time. Then for any $x$ and $\xi_i$, $g_i(x,\xi_i):=\nabla F_i(x,\xi_i)$ is an unbiased estimator of $\nabla f_i(x)$. For another example, let $f_i(x):=(1/|S_i|) \sum_{\zeta_j\in S_i}F(x,\zeta_j)$ denote an empirical risk function, where $S_i$ is the local dataset of agent $i$. Under this setting, the gradient estimation of $f_i(x)$ can incur noise from various sources such as sampling and discretization \cite{li2021compressed}.
	
	Problem (\ref{eq:eloss}) has been studied extensively in the literature using various distributed algorithms \cite{tsitsiklis1986distributed,nedic2009distributed,nedic2010constrained,lobel2011distributed,jakovetic2014fast,kia2015distributed,shi2015extra,di2016next,qu2017harnessing,nedic2017achieving,xu2017convergence,pu2020push}, among which the distributed gradient descent (DGD) method proposed in \cite{nedic2009distributed} has drawn the greatest attention.
	Recently, distributed implementation of stochastic gradient algorithms has received considerable interest.
	Several works 
	have shown that distributed methods may compare with their centralized counterparts under certain conditions. For example, the work in \cite{chen2012limiting,chen2015learning,chen2015learning2} first showed that, with a sufficiently small constant stepsize, a distributed stochastic gradient method achieves comparable performance to a centralized method in terms of the steady-state mean-square-error.
	
	More recent works have demonstrated that various distributed stochastic gradient methods enjoy the so-called ``asymptotic network independence'' property, that is, the generated iterates asymptotically converge at the same network independent rate compared to centralized stochastic gradient descent (SGD) for minimizing smooth objective functions \cite{spiridonoff2020robust,pu2020asymptotic,pu2019sharp,lian2017can,koloskova2019decentralized,tang2018d,xin2020improved}. Nevertheless, it always takes a certain number of updates before this network independent convergence rate can be reached, which is always network dependent and referred to as \emph{transient time} of the algorithm. For example, when considering smooth and strongly convex objective functions, the shortest transient time that has been shown so far is achieved by the distributed stochastic gradient descent (DSGD) method \cite{pu2019sharp}, which behaves as $\mo\left(\frac{n}{(1-\lambda_2)^2}\right)$, where $1-\lambda_2$ stands for the spectral gap of the mixing matrix. For general undirected networks such as ring graphs, the transient time of DSGD is as large as $\mo\left(n^5\right)$, which can be significant for large-scale networks. Therefore, it is critical to develop new algorithms with shorter transient times.
	
	Recently, a class of closely related algorithms, termed EXTRA \cite{shi2015extra}, $D^2$ \cite{tang2018d}, Exact Diffusion \cite{yuan2018exact,yuan2020influence} and NIDS \cite{li2019decentralized} have been proposed to solve problem \eqref{eq:eloss} under both exact and stochastic gradient settings. These methods achieve linear convergence for minimizing smooth and strongly convex objective functions 
	under exact gradient information \cite{li2019decentralized}. When only stochastic gradients are available, the work in \cite{yuan2020influence} compared the performance of various distributed optimization methods under constant stepsizes and demonstrated the effectiveness of Exact Diffusion compared to other distributed stochastic gradient algorithms.
	
	In this work, we consider a distributed stochastic gradient algorithm adapted from NIDS/Exact Diffusion/$D^2$, termed EDAS (Exact Diffusion with Adaptive Stepsizes) and perform a non-asymptotic analysis. In addition to showing the asymptotic network independence property of EDAS, we carefully identify its non-asymptotic convergence rate as a function of characteristics of the objective functions and the underlying network. Moreover, we derive the transient time needed for EDAS to achieve the network independent convergence rate, which behaves as $\mo\left(\frac{n}{1-\lambda_2}\right)$. This is the shortest transient time for minimizing strongly convex and smooth objective functions in a distributed fashion to the best of our knowledge. Numerical experiments corroborate and strengthen the theoretical findings.
	
	\subsection{Related Work}
	
	Our work is closely related to the concept of ``asymptotic network independence'' in distributed optimization over networks, that is, asymptotically, the distributed
	algorithm converges to the optimal solution at a comparable rate to a centralized algorithm with the same computational power. In addition to \cite{chen2012limiting,chen2015learning,chen2015learning2} discussed before, the work in \cite{morral2014success,morral2017success} proved that distributed stochastic approximation methods perform asymptotically
	as well as centralized schemes by means of a central limit theorem. The papers \cite{pu2017flocking,pu2018swarming} demonstrated the advantages of distributively implementing a stochastic gradient method assuming that sampling times are random and non-negligible.  
	The first ``asymptotic network independence" result  appeared in \cite{towfic2016excess}, where it assumed that all the local functions $f_i$ have the same minimum. 
	A recent paper \cite{spiridonoff2020robust} discussed an algorithm that asymptotically performs as well as the best bounds on centralized stochastic gradient descent subject to possible message losses, delays, and asynchrony. In a parallel recent work \cite{koloskova2019decentralized}, a similar result was presented with a further compression technique which allowed agents to save on communication. The work in \cite{pu2020distributed} considered a distributed stochastic gradient tracking method (DSGT) which asymptotically performs as well as centralized stochastic gradient descent. The result was generalized to the setting of directed communication graphs in \cite{qureshi2020s}.
	For nonconvex objective functions, the paper \cite{lian2017can,tang2018d} proved that decentralized algorithms can achieve a linear speedup similar to a centralized algorithm when the number of iterates $k$ is large enough. The work in \cite{assran2018stochastic} further extended the theory to directed communication networks  for training deep neural networks.
	A recent survey \cite{pu2020asymptotic} provided a detailed discussion on this topic.
	
	When restricted to considering strongly convex and smooth objective functions, several recent works have not only shown asymptotic network independence  for the proposed algorithms, but also provided the transient times for the algorithms to reach the network independent convergence rate. \cref{tab:time} compares the transient times for some existing algorithms, which are functions of the network characteristics. In particular, the paper \cite{pu2019sharp} showed that for the distributed stochastic gradient descent(DSGD) method, it takes $\mo\left(\frac{n}{(1-\lambda_2)^2}\right)$ ($1-\lambda_2$ stands for the spectral gap of the mixing matrix) to achieve network independence. This result was also proved to be sharp. For another distributed stochastic gradient tracking (DSGT) method, the paper \cite{xin2020improved} showed the transient time is $\mo\left(\frac{n}{(1-\lambda_2)^3}\right)$.
	
	
	
	
	\begin{table}[H]\label{tab:time}
		\begin{center}
			\begin{tabular}{m{1.6 cm} m{2.2 cm} m{2.1cm} m{2cm}c}
				\hline
				& DSGD with compression & DSGD \& Diffusion & DSGT & EDAS (this paper) \\ \hline
				Transient Time & $\mathcal{O}\left(\frac{n}{(1-\lambda_2)^4}\right)$ \cite{koloskova2019decentralized} &
				$\mathcal{O}\left(\frac{n}{(1-\lambda_2)^2}\right)$ \cite{pu2019sharp} &
				$\mathcal{O}\left(\frac{n}{(1-\lambda_2)^3}\right)$ \cite{xin2020improved} &
				$\mathcal{O}\left(\frac{n}{1-\lambda_2}\right)$ \\ \hline
			\end{tabular}
		\end{center}
		\caption{Transient times of distributed stochastic gradient methods for minimizing strongly convex and smooth objective functions.}
	\end{table}
	
	\subsection{Main Contribution}
	\label{sec:contri}
	We summarize the main contribution of the paper as follows. First, we consider a distributed stochastic gradient algorithm (EDAS) adapted from NIDS/Exact Diffusion and perform a non-asymptotic analysis. For strongly convex and smooth objective functions, we show in \cref{thm:mag_extra} that EDAS asymptotically achieves the same network independent convergence rate of centralized stochastic gradient descent (SGD). We also carefully characterize the convergence rate of EDAS as a function of the characteristics of the objective functions and the communication network. 
	
	Second, by comparing the convergence result of EDAS to that of SGD, we derive the transient time needed for EDAS to achieve the network independent rate, which is upper bounded by $K_T=\mo\left(\frac{n}{1-\lambda_2}\right)$ under mild conditions. This is the shortest transient time so far to the best of our knowledge, and significantly improves upon the  shortest transient time previously obtained, i.e., $\mo\left(\frac{n}{(1-\lambda_2)^2}\right)$. This is the main contribution of the paper.
	
	Finally, we provide two simulation examples that corroborate and strengthen the obtained theoretical results. For solving a ``hard" problem with EDAS, the observed transient times are consistent with the theoretical upper bounds for both the ring network topology and the square grid network topology, which verifies the sharpness of the theoretical results. For the problem of classifying handwritten digits with logistic regression, the observed transient times are even better than the upper bounds. For both problems, we compare the performance of EDAS with other existing algorithms and show the superiority of EDAS.
	
	\subsection{Notation}\label{sec:notation}
	Throughout the paper, we use column vectors if not otherwise specified. Let each agent $i$ hold a local copy $x_{i,k}\in\R^p$ of the decision variable at the $k$-th iteration. We denote
	\begin{align*}
		\x_k & := (x_{1,k}, x_{2,k}, ..., x_{n,k})^{\T}\in \R^{n\times p},\\
		\bar x_k & : = \frac{1}{n}\sum_{i=1}^k x_{i,k}.
	\end{align*} 
	We also define the following aggregative functions, 
	\begin{align*}
		F(\x_k)&:= \sum_{i=1}^n f_i(x_{i,k}),\\
		\nabla F(\x_k) &:= \left(\nabla f_1(x_{1,k}), \nabla f_2(x_{2,k}), ..., \nabla f_n(x_{n,k})\right)^{\T}\in\R^{n\times p},\\
		\g_k &:= \left(g_1(x_{1,k}, \xi_{1,k}), g_2(x_{2,k}, \xi_{2,k}\right), ..., g_n(x_{n,k},\xi_{n,k}))^{\T}\in\R^{n\times p},
	\end{align*}
	and write
	\begin{align*}
		\bar\nabla  F(\x_k) &:= \frac{1}{n}\sum_{i=1}^n \nabla f_i(x_{i,k}),\\
		\bar\sigma^2 & := \frac{1}{n}\sum\limits_{i=1}^n \sigma_i^2,\\
		\bar g_k  &:= \frac{1}{n}\sum_{i=1}^n g_i(x_{i,k},\xi_{i,k}).
	\end{align*}
	For two vectors $a,b$ of the same dimension, $\langle a, b\rangle$ denotes the inner product. For two matrices $A,B\in\R^{n\times p}$, we define 
	\begin{equation*}
		\left\langle A,B\right\rangle := \sum_{i=1}^n \langle A_i, B_i\rangle,
	\end{equation*}
	where $A_i$ (respectively, $B_i$) represents the $i$-th row of $A$ (respectively, $B$). $\Vert\cdot \Vert$ denotes $\ell_2$-norm for vectors and Frobenius norm for matrices. 
	
	Next we specify the conditions concerning the local cost functions and the communication network.
	\begin{assumption}\label{edfi}
		Each $f_i:\R^p\rightarrow \R$ has $L$-Lipschitz continuous gradient, and $f:\R^p\rightarrow \R$ is $\mu$-strongly convex, i.e., $\forall x, x'\in\R^p$,
		\begin{align*}
			\Vert \nabla f_i(x)-\nabla f_i(x')\Vert\leq L\Vert x-x'\Vert,\quad \forall i,
		\end{align*}
		\begin{align*}
			\langle\nabla f(x)-\nabla f(x'), x-x'\rangle\geq \mu \Vert x-x'\Vert^2.\\
		\end{align*}
	\end{assumption}
	Under Assumption \ref{edfi}, problem \eqref{eq:eloss} has a unique optimal solution $x^*\in\mathbb{R}^p$.
	
	The agents are connected through a graph $\mathcal{G}=(\mathcal{N}, \mathcal{E})$, where $\mathcal{N}$ is the set of nodes and $\mathcal{E}\subseteq \mathcal{N}\times \mathcal{N}$ denotes the set of edges. Regarding the graph $\mathcal{G}$ and the corresponding mixing matrix $W=[w_{ij}]\in \R^{n\times n}$, we make the following standing assumption.
	\begin{assumption}\label{edw}
		The graph $\mathcal{G}$ is undirected and strongly connected. There exists a link from $i$ and $j$ ($i\neq j$) in $\mathcal{G}$ if and only if $w_{ij},w_{ji}>0$; otherwise, $w_{ij}=w_{ji}=0$. The mixing matrix $W$ is nonnegative, symmetric and stochastic, i.e., $W\1 =\1$. In addition, suppose the smallest eigenvalue of $W$ satisfies $\lambda_n\geq\underline{\lambda}>0$ for some fixed $\bar{\lambda}$ independent of $n$.\footnote{To satisfy such an assumption, it is convenient to let the mixing matrix $W$ be chosen as $W = \beta I_n + (1-\beta)\tilde{W}$ for some $\beta\in(\frac{1}{2},1)$ with $\tilde{W}$ being nonnegative, symmetric and stochastic. Therefore, the condition is not restrictive.}
	\end{assumption}
	
	The rest of this paper is organized as follows. In \cref{sec:extra} we present the EDAS algorithm  and its transformed error dynamics in the form of a primal-dual-like algorithm. Next we provide some supporting lemmas and prove the sublinear convergence rate of EDAS in \cref{sec:ana}. The main convergence results of the algorithm are presented in \cref{sec:main_res}, where we derive the convergence rate of EDAS as a function of the characteristics of the objective functions and the mixing matrix, and subsequently, obtain the transient time of EDAS. Numerical experiments are provided in \cref{sec:experiments}, and we conclude the paper in \cref{sec:conclusions}.
	
	\section{Algorithm}
	\label{sec:extra}
	We study the following algorithm adapted from EXTRA \cite{shi2015extra}, Exact Diffusion \cite{yuan2018exact}, NIDS \cite{li2019decentralized} and $D^2$ \cite{tang2018d}, where the main characteristics of EDAS lie in the diminishing stepsize policy $\alpha_{k}=\mathcal{O}_k (\frac{1}{k})$: 
	
	\begin{algorithm}[H]
		\begin{algorithmic}[1]
			\STATE Initialize $x_{i,0}\in\R^p$ for all agent $i\in\mathcal{N}$, determine $W = [w_{ij}]\in\R^{n\times n}$, stepsize $\alpha_{k}=\frac{\theta}{\mu(k+m)}$ for some $\theta>0$ and $m>0$.
			\FOR{$k=0, 1, 2, ...$}
			\FOR{$i = 1, 2, ..., n$}
			\STATE Agent $i$ acquires a stochastic gradient $g_i(x_{i,k}, \xi_{i,k})\in\R^p$.
			\IF{k = 0}
			\STATE $x_{i,k+1} = \sum\limits_{j=1}^n w_{ij}\left(x_{j,k} - \alpha_k g_j(x_{j, k}, \xi_{j,k})\right)$ 
			\ELSE 
			\STATE $x_{i, k+1} = \sum\limits_{j=1}^n w_{ij}\left(2x_{j, k} - x_{j, k-1} - \alpha_k g_j(x_{j, k}, \xi_{j,k}) + \alpha_{k-1} g_{j}(x_{j, k-1}, \xi_{j, k-1})\right)$
			\ENDIF
			\ENDFOR
			\ENDFOR
		\end{algorithmic}
		\caption{EDAS (Exact Diffusion with Adaptive Stepsizes)}
		\label{alg:extra}
	\end{algorithm}
	
	Using the notation in \cref{sec:notation}, EDAS can be written in the following compact form:
	\begin{subequations}\label{extra_ak}
		\begin{align}
			\x_1 &= W(\x_0 - \alpha_0 \g_0), \\
			\x_{k+1} &= W(2\x_k-\x_{k-1}-\alpha_k\g_k+\alpha_{k-1}\g_{k-1}),\;\forall k\ge 1.
		\end{align}
	\end{subequations}
	
	For further investigation, note that since $W$ is nonnegative, symmetric and stochastic, 
	it has a spectral decomposition given by 
	\begin{equation}
		W = Q\Lambda Q^{\T},
	\end{equation} 
	where 
	\begin{equation*}
		\begin{aligned}
			\Lambda := \mathrm{diag}(\lambda_1,\lambda_2,...,\lambda_n),
		\end{aligned}
	\end{equation*}
	with $1=\lambda_1 >\lambda_2\ge \cdots \ge \lambda_n\ge -1$.
	Thus $I_n-W = Q(I_n-\Lambda)Q^{\T}$. 
	
	Let $V:= (I_n-W)^{1/2}$. EDAS can be rewritten in the following equivalent form that resembles a primal-dual algorithm (see \cite{yuan2018exact}): let $\y_0 = \mathbf{0}$, for $k\ge 0$,
	\begin{subequations}\label{pdextra}
		\begin{align}
			\x_{k+1} &= W(\x_k - \alpha_k\g_k) - V\y_k,\\
			\y_{k+1} &= \y_k + V\x_{k+1}.
		\end{align}
	\end{subequations}
	To verify the equivalence relationship, note that from \eqref{pdextra}, when $k\ge 1$, we have
	\begin{equation*}
		\begin{aligned}
			\x_{k+1}-\x_{k}&=W\left(\x_{k}-\x_{k-1}-\alpha_k\g_k
			+\alpha_{k-1}\g_{k-1}\right)-V\y_k+V\y_{k-1}\\
			&=W\left(\x_{k}-\x_{k-1}-\alpha_k\g_k
			+\alpha_{k-1}\g_{k-1}\right)-V(\y_k - \y_{k-1})\\
			&=W (2 \x_k-\x_{k-1}-\alpha_k\g_k+\alpha_{k-1}\g_{k-1})-\x_{k}.
		\end{aligned}
	\end{equation*}
	
	\subsection{Optimality Condition}
	
	To facilitate the analysis of EDAS, we consider the optimality condition for problem \eqref{eq:eloss}. From the condition we are able to construct the error dynamics for \cref{alg:extra}.
	
	\begin{lemma}[Optimality Condition \cite{yuan2019performance}]\label{lem:optcond}
		Under \cref{edfi,edw}, if there exists some $(\x^*, \y_k^*)$ that satisfies:\footnote{Note that $\y_k^*$ depends on $k$.}
		\begin{subequations}
			\label{eq:optcond}
			\begin{align}
				\alpha_kW\nabla F(\x^*) + V\y_k^* &= \mathbf{0}\label{eq:yk_star}\\
				V\x^* &= \mathbf{0},
			\end{align}
		\end{subequations}
		where $\x^* : = (x_1^*, x_2^*, ..., x_n^*)^{\T}\in\R^{n\times p}$, then it holds that 
		$$
		x_1^* = x_2^* = ... = x_n^* = x^*.
		$$
	\end{lemma}
	
	\begin{proof}
		We follow the proof of Lemma 2 in \cite{yuan2018exact}, where the stepsize $\alpha_k$ was assumed to be a constant.
		First, since $\mathrm{null}(V)=\mathrm{span}\{\mathbf{1}_n\}$, we have
		\begin{equation}
			\0 = V\x^*\Longleftrightarrow \x^*\in\mathrm{null}(V),
		\end{equation}
		and hence $x^*_1=x^*_2=...=x^*_n$.
		
		Additionally, multiplying $\mathbf{1}^{\T}$ on both sides of \cref{eq:yk_star} 
		and noticing that $\mathbf{1}^{\T}V=\0$, we have
		\begin{equation}
			\begin{aligned}
				\alpha_k\mathbf{1}^{\T}W\nabla F(\x^*) = \alpha_k\mathbf{1}^{\T}\nabla F(\x^*)
				= \alpha_k \sum_{i=1}^n \nabla f_i(x_i^*)
				= 0,
			\end{aligned}
		\end{equation}
		where we use the column stochastic property of $W$, i.e., $\1^{\T} W = \1^{\T}$. Relation $\sum\limits_{i=1}^n \nabla f_i(x_i^*)=0$ then implies that $x_1^*=x_2^*=...=x_n^*$ must coincide with the minimizer $x^*$ of problem \eqref{eq:eloss}.
	\end{proof}
	
	We next show the existence of $(\x^*,\y_k^*)$  and how to select a unique pair of $(\x^*,\y_k^*)$ in light of the singularity of $V$.
	\begin{lemma}
		\label{lem:existence}
		Under \cref{edfi,edw}, there exists $(\x^*,\y_k^*)$ that satisfies \eqref{eq:optcond}. In particular, 
		we can choose
		\begin{subequations}
			\label{xy_unique}
			\begin{align}
				& \x^* = \1 (x^*)^{\T},\\
				& \y_k^* = -V^{-} (\alpha_k W \nabla F(\x^*)),
			\end{align}
		\end{subequations}
		where $V^{-}$ is the pseudoinverse of $V$ given by 
		\begin{equation*}
			V^{-} = Q \left(\sqrt{I_n-\Lambda}\right)^{-} Q^{\T},
		\end{equation*}
		in which $Q$ is the same orthogonal matrix as in decomposing $W$, and $\left(\sqrt{I_n-\Lambda}\right)^{-} : = \mathrm{\diag}(d_1, d_2, ..., d_n)$ is the pseudoinverse of $\sqrt{I_n-\Lambda}$ with
		\begin{equation*}
			d_i = \begin{cases}
				\frac{1}{\sqrt{1-\lambda_i}} & \text{if }1-\lambda_i\neq 0,\\
				0 & \text{otherwise.}
			\end{cases}
		\end{equation*}
	\end{lemma} 
	
	\begin{proof}
		The proof for the existence of $(\x^*,\y_k^*)$ is inspired by Lemma 3 in \cite{yuan2018exact}.
		Let $\x^* = \1 (x^*)^{\T}$. Then noticing that $V\mathbf{1}_n=\0$, we have $V\x^*=0$.
		To prove that $\y_k^*$ exists for all $k\geq 0$ is equivalent to showing that the following linear system is consistent w.r.t. $\y_k^*$:
		\begin{equation}
			\label{eq:lrs}
			V\y_k^* = -\alpha_kW \nabla F(\x^*).
		\end{equation}
		Hence we need only show that $-\alpha_k W \nabla F(\x^*)$ lies in $\range(V)$ for all $k\geq0$. 
		By \cref{edw},
		\begin{equation*}
			\1^{\T}(-\alpha_kW\nabla F(\x^*)) = -\alpha_k\sum_{i=1}^n f_i(x^*) = 0,
		\end{equation*}
		we have
		\begin{align*}
			(-\alpha_kW\nabla F(\x^*)) &\perp \mathrm{span}\{\1\} = \mathrm{null}(V),\\
			\Leftrightarrow (-\alpha_kW\nabla F(\x^*)) &\in \range(V^{\T})=\range(V),
		\end{align*}
		which proves the existence of $\y_k^*$. 
		
		We next verify that $\y_k^* = -V^{-} (\alpha_k W \nabla F(\x^*))$ satisfies \eqref{eq:lrs}. First of all, it is easy to see that the choice of $V^{-}$ is indeed a pseudoinverse of $V$.
		Then since the linear system \eqref{eq:lrs} is consistent, we have that $V^{-}(-\alpha_k W \nabla F(\x^*))$ is a solution of \eqref{eq:lrs}, or
		$V\y_k^* = -\alpha_kW \nabla F(\x^*)$.
		
	\end{proof}
	\begin{remark}
		Throughout the paper, we let  $(\x^*,\y_k^*)$ be the unique pair given in \eqref{xy_unique}.
	\end{remark}

	\subsection{Preliminary Analysis}
	
	In light of the optimality condition of problem \eqref{eq:eloss} given in \Cref{lem:optcond}, we consider the following error dynamics as the starting point of the convergence analysis for EDAS.
	\begin{lemma}
		Denote  $\tx_k:=\x_k-\x^*$, $\ty_{k} := \y_k-\y_k^*$ and $\s_k:=\nabla F(\x_k) - \g_k$. We have
		\begin{equation}\label{edym}
			\begin{aligned}
				\left(
				\begin{array}{c}
					\tx_{k+1}\\
					\ty_{k+1}
				\end{array}
				\right)
				&
				=B
				\left(
				\begin{array}{c}
					\tilde{\x}_{k} + \alpha_k(\nabla F(\x^*) - \nabla F(\x_k))\\
					\tilde{\y}_{k}\\
				\end{array}
				\right)
				+\alpha_k B
				\left(
				\begin{array}{c}
					\s_k\\
					\0
				\end{array}
				\right)\\ 
				&\quad + \left(\begin{array}{c}
					\0  \\
					\y_{k}^*-\y_{k+1}^* 
				\end{array}\right),
			\end{aligned}
		\end{equation}
		where
		\begin{equation}
			B := \left(
			\begin{array}{cc}
				W & -V\\
				VW & W\\
			\end{array}
			\right).
		\end{equation}
	\end{lemma}
	\begin{proof}
		From relation \eqref{pdextra} and the definitions of $\tx_k$, $\ty_{k}$ and $\s_k$, we have
		\begin{equation}
			\label{edx}
			\begin{aligned}
				\tx_{k+1} &= \x_{k+1} - \x^*\\
				&= W(\x_k - \alpha_k\g_k) - V\y_k - \x^*\\
				&= W\left[\tx_k + \alpha_k\left(\nabla F(\x^*)-\nabla F(\x_k)\right)\right] + \alpha_k W \s_k - \alpha_k W\nabla F(\x^*) - V\y_k\\
				&= W\left[\tx_k + \alpha_k\left(\nabla F(\x^*)-\nabla F(\x_k)\right)\right] - V\ty_k + \alpha_k W \s_k,
			\end{aligned}
		\end{equation}
		and
		\begin{equation}\label{edy}
			\ty_{k+1} = \y_{k+1}-\y_{k+1}^* = \y_k+V\x_{k+1}-\y_{k+1}^*  =  \ty_k + V\tx_{k+1} +(\y_{k}^*-\y_{k+1}^*).
		\end{equation}
		
		Substituting \eqref{edx} into relation \eqref{edy}, we obtain
		\begin{equation}
			\begin{aligned}
				\left(
				\begin{array}{c}
					\tx_{k+1}\\
					\ty_{k+1}
				\end{array}
				\right)
				&=
				\left(\begin{array}{cc}
					W&-V\\
					VW&W
				\end{array}
				\right)
				\left(
				\begin{array}{c}
					\tilde{\x}_{k} + \alpha_k(\nabla F(\x^*) - \nabla F(\x_k))\\
					\tilde{\y}_{k}\\
				\end{array}
				\right)\\
				&\quad +\alpha_k
				\left(
				\begin{array}{cc}
					W & -V\\
					VW & W\\
				\end{array}
				\right)
				\left(
				\begin{array}{c}
					\s_k\\
					\0
				\end{array}
				\right)
				+ \left(\begin{array}{c}
					\0  \\
					\y_{k}^*-\y_{k+1}^* 
				\end{array}\right).
			\end{aligned}
		\end{equation}
	\end{proof}
	
	In the next lemma, we show the matrix $B\in\R^{2n\times 2n}$ has an eigendecomposition in the form of  $B = UDU^{-1}$, inspired by the arguments in \cite{yuan2018exact}. In light of the decomposition, we are able to derive the transformed error dynamics for EDAS, which is then used for performing the non-asymptotic analysis for EDAS in Section \ref{sec:ana}. 
	\begin{lemma}\label{lem:B_decomposition}
		Under \cref{edw}, the matrix $B$ has an eigendecomposition given by
		\begin{equation*}
			B = U D U^{-1} := \begin{pmatrix}
				\1_n & \0 & c U_{R,u}\\
				\0 & \1_n & c U_{R,l}
			\end{pmatrix} \mathrm{diag}(1, 1, D_1)
			\begin{pmatrix}
				\frac{1}{n}\1_n^{\T} & \0\\
				\0 & \frac{1}{n}\1_n^{\T}\\
				\frac{1}{c} U_{L,l} & \frac{1}{c} U_{L,r} 
			\end{pmatrix}
		\end{equation*}
		where $D_1$ is a diagonal matrix with complex entries,  matrices $U_{R,u},U_{R,l}\in\mathbb{R}^{n\times (n-2)}$, $U_{L,l},U_{L,r}\in\mathbb{R}^{(n-2)\times n}$, and $c>0$ is an arbitrary scaling parameter. 
		
		Denote 
		\begin{equation*}
			U_R :=\begin{pmatrix}
				U_{R,u}\\
				U_{R,l}
			\end{pmatrix},\quad 
			U_L := \begin{pmatrix}
				U_{L,l} & U_{L,r} 
			\end{pmatrix}.
		\end{equation*}
		Then the product $\Vert U_L\Vert_2^2\Vert U_R\Vert_2^2$ has an upper bound that is independent of $n$. More specifically, 
		\begin{equation*}
			\Vert U_L\Vert_2^2\Vert U_R\Vert_2^2\leq \frac{1}{\lambda_n}\leq \frac{1}{\underline{\lambda}}.
		\end{equation*}
	\end{lemma}
	\begin{proof}
		see \cref{sec:prf_KL_KR}.
	\end{proof}
	
	Based on \cref{lem:B_decomposition}, we derive the transformed error dynamics for EDAS by multiplying  $U^{-1}$ on both sides of the error dynamics \eqref{edym}. In this way, we decompose the error $(\tx_k, \ty_k)$ into two parts: $\bar z_k$ that measures the difference between the average iterate and $x^*$, and $\check z_k$ that represents the remaining error. Then we further study the relationship between $\bar z_k$ and $\check z_k$ in \cref{sec:ana}.
	
	\begin{lemma}\label{lem:tedy}
		Under \cref{edfi,edw}, the transformed error dynamics for \cref{alg:extra} is given by
		\begin{equation}
			\begin{aligned}
				\left(\begin{array}{c}
					\bar z_{k+1}\\
					\check z_{k+1}
				\end{array}
				\right)
				&=\left(\begin{array}{c}
					\bar z_k + \frac{\alpha_k}{n}\1_n^{\T}\left(\nabla F(\x^*) - \nabla F(\x_k)\right)\\
					D_1\check{z}_k + \frac{\alpha_k}{c} D_1 U_{L,l}\left(\nabla F(\x^*) - \nabla F(\x_k)\right)
				\end{array}
				\right)\\
				&\quad + \alpha_k \left(
				\begin{array}{c}
					\frac{1}{n}\mathbf{1}_n^{\T} \s_k\\
					\frac{1}{c}D_1U_{L,l}\s_k
				\end{array}
				\right)
				+ 
				\left(\begin{array}{c}
					0  \\
					\frac{1}{c}U_{L,r}(\y_{k}^*-\y_{k+1}^*)
				\end{array}\right)
				,
			\end{aligned}
		\end{equation}
		where
		\begin{equation}\label{rel}
			\left(
			\begin{array}{c}
				\bar z_k\\
				\check z_k
			\end{array}
			\right) 
			:= 
			\left(
			\begin{array}{c}
				\frac{1}{n}\1_n^{\T}\tx_k\\
				\frac{1}{c}(U_{L,l}\tx_k+U_{L,r}\ty_k)
			\end{array}
			\right)\; \forall k.
		\end{equation}
	\end{lemma}
	
	\begin{proof}
		Denote 
		\begin{equation}
			\begin{pmatrix}
				\bar{z}_{k}\\
				\hat{z}_{k}\\
				\check{z}_{k}
			\end{pmatrix}:=
			U^{-1}
			\begin{pmatrix}
				\tx_{k}\\
				\ty_{k}
			\end{pmatrix}
			= \begin{pmatrix}
				\frac{1}{n}\1_n^{\T}\tx_k\\
				\frac{1}{n}\1_n^{\T}\ty_k\\
				\frac{1}{c}U_{L,l}\tx_{k} + \frac{1}{c}U_{L,r}\ty_{k}
			\end{pmatrix},
		\end{equation}
		
		Multiplying $U^{-1}$ on both sides of \eqref{edym} leads to 
		\begin{align*}
			\begin{pmatrix}
				\bar{z}_{k+1}\\
				\hat{z}_{k+1}\\
				\check{z}_{k+1}
			\end{pmatrix}
			&
			= DU^{-1} 
			\begin{pmatrix}
				\tx_k + \alpha_k \left(\nabla F(\x^*) - \nabla F(\x_k)\right)\\
				\ty_k
			\end{pmatrix} + \alpha_k DU^{-1} 
			\begin{pmatrix}
				\s_k\\
				\mathbf{0}
			\end{pmatrix}\\
			&\quad   + U^{-1} \begin{pmatrix}
				\mathbf{0}\\
				\y_{k}^*-\y_{k+1}^*
			\end{pmatrix}\\
			&= \begin{pmatrix}
				\bar{z}_{k} + \frac{\alpha_k}{n}\1_n^{\T}\left(\nabla F(\x^*) - \nabla F(\x_k)\right)\\
				\hat{z}_{k}\\
				D_1\check{z}_k + \frac{1}{c}D_1U_{L,l}\alpha_k\left(\nabla F(\x^*) - \nabla F(\x_k)\right)
			\end{pmatrix}
			+ \begin{pmatrix}
				\frac{\alpha_k}{n}\1_n^{\T}\s_k\\
				\mathbf{0}\\
				\frac{\alpha_k}{c}D_1U_{L,l}\s_k
			\end{pmatrix}\\
			&\quad + \begin{pmatrix}
				\mathbf{0}\\
				\frac{1}{n}\1_n^{\T}(\y_{k}^*-\y_{k+1}^*)\\
				\frac{1}{c}U_{L,r}(\y_{k}^*-\y_{k + 1}^*)
			\end{pmatrix}.
		\end{align*}
	\end{proof}
	
	\section{Analysis}
	\label{sec:ana}
	In this section, we provide some useful supporting lemmas based on \cref{lem:tedy}, and then prove the sublinear convergence rate of EDAS in \cref{lem:hat_W}. For the proof of sublinear convergence, we first propose \cref{barzk,lem:checkz} to derive two coupled recursions for $\E\left[\norm{\bar z_k}^2\right]$ and $\E\left[\norm{\check z_k}^2\right]$ respectively according to \cref{lem:tedy}. Then we introduce a Lyapunov function $H_k:= \E\left[\norm{\bar z_k}^2\right] + \omega_k \E\left[\norm{\check z_k}^2\right]$, and by constructing the recursion for $H_k$, we derive the sublinear convergence rate.
	
	\subsection{Supporting Lemmas}
	\label{sec:s_lems}
	In this part, we provide two important coupled recursions regarding $\E\left[\norm{\bar z_k}^2\right]$ and $\E\left[\norm{\check z_k}^2\right]$, respectively. Before introducing these results, we state some preliminary lemmas.
	\begin{lemma}
		\label{lem:bar_s}
		Suppose \cref{ngrad} holds, then
		\begin{subequations}
			\begin{align}
				\E\left[\bar{s}_k| \mF_{k}\right] & = \0,\label{eq:unbias}\\
				\E\left[\norm{\bar{s}_k}^2| \mathcal{F}_{k}\right] 
				& \leq \frac{\bar\sigma^2}{n}. \label{eq:bound_var}
			\end{align}
		\end{subequations}
	\end{lemma} 
	
	The results of \cref{lem:bar_s} directly come from \cref{ngrad}.
	
	\begin{lemma}\label{lem:from_L}
		Under \cref{edfi}, there holds
		\begin{equation*}
			\norm{\bar{\nabla} F(\1\bar{x}_k^{\T}) - \bar{\nabla}F(\x_k)} \leq \frac{L}{\sqrt{n}}\norm{cU_{R,u}\check{z}_k}.
		\end{equation*}
	\end{lemma}
	\begin{proof}
		By definition,
		\begin{align}
			\norm{\bar{\nabla} F(\1\bar{x}_k^{\T}) - \bar{\nabla}F(\x_k)} 
			&= \norm{\frac{1}{n}\sum_{i=1}^n\left(\nabla f_i(\bar{x}_k) - \nabla f_i(x_{i,k})\right)}\nonumber\\
			&\leq \frac{L}{n}\sum_{i=1}^n\norm{\bar{x}_k-x_{i,k}}\label{eq:bar_L}\\
			&= \frac{L}{n}\sum_{i=1}^n\norm{\bar{x}_k-x^* + x^*-x_{i,k}}\nonumber\\
			& \le \frac{L}{\sqrt{n}}\norm{\tx_k - \1\bar{z}_k}\label{eq:bar_trans_tx}\\
			&= \frac{L}{\sqrt{n}}\norm{cU_{R,u}\check{z}_k},\label{eq:bar_res}
		\end{align}
		where \eqref{eq:bar_L} holds because of \cref{edfi}, inequality \eqref{eq:bar_trans_tx} results from Cauchy-Schwarz's inequality, and \cref{eq:bar_res} comes from
		\begin{equation}
			\begin{pmatrix}
				\tx_{k}\\
				\ty_{k}
			\end{pmatrix}
			=U 
			\begin{pmatrix}
				\bar{z}_{k}\\
				\hat{z}_{k}\\
				\check{z}_{k}
			\end{pmatrix},
		\end{equation}
		which implies $\tx_{k} = \1\bar{z}_k + cU_{R,u}\check{z}_k$.
	\end{proof}
	
	In light of Lemma 10 in \cite{qu2017harnessing}, we have the following contraction result: 
	\begin{lemma}
		Under \cref{edfi}, suppose $\alpha_k\leq \frac{2}{\mu+L}$. There holds
		\begin{equation}
			\label{eq:fix_rate}
			\norm{\bar{z}_k - \alpha_k \bar{\nabla}F(\1\bar{x}_k)}\leq (1-\alpha_k\mu)\norm{\bar{z}_k}.
		\end{equation}
	\end{lemma}
	
	With the above results in hand, we provide an upper bound for $\E\left[\norm{\bar z_k}^2\right]$ in \cref{barzk}, following \cref{lem:tedy}. 
	
	
	
	\begin{lemma}\label{barzk}
		Under \cref{alg:extra} with \cref{ngrad,edfi,edw,}, suppose $\alpha_k\leq \min\{\frac{1}{3\mu}, \frac{2}{\mu + L}\}$, we have
		\begin{equation}
			\E\left[\norm{\bar{z}_{k+1}}^2\right] \leq (1-\frac{3}{2}\alpha_k\mu)\E\left[\norm{\bar{z}_{k}}^2\right] + \frac{3\alpha_kL^2c^2\norm{U_R}^2}{\mu n}\E\left[\norm{\check{z}_k}\right] + \frac{\alpha_k^2\bar{\sigma}^2}{n}.
		\end{equation}
	\end{lemma}
	
	\begin{proof}
		Denote $\bar{s}_k:=\frac{1}{n}\sum_{i=1}^{n}\s_{i, k}^{\T}$, where $s_{i, k}$ is the $i$-th row of $\mathbf{s}_k$. We have from \cref{lem:tedy} that
		\begin{align*}
			\bar{z}_{k+1} &= \bar z_k + \alpha_k\left(  \bar{\nabla} F(\x^*) - \bar{\nabla} F(\x_k) +  \bar{s}_k\right)   \\
			&= \bar z_k - \alpha_k \bar{\nabla} F(\x_k)+ \alpha_k\bar{s}_k,
		\end{align*}
		where the second equality follows from $\bar{\nabla} F(\x^*) = \sum_{i=1}^n\nabla f_i(x^*) = 0$ in light of the optimality of $x^*$. Hence
		\begin{align}
			\E\left[\norm{\bar{z}_{k+1}}^2|\mF_k\right] &= \E\left[\norm{\bar z_k - \alpha_k \bar{\nabla} F(\x_k)+ \alpha_k\bar{s}_k}^2|\mF_k\right]\nonumber\\
			&= \norm{\bar{z}_k - \alpha_k\bar{\nabla}F(\x_k)}^2 + \alpha_k^2\E\left[\norm{\bar{s}_k}^2|\mF_k\right]\label{eq:cond_barz}.
		\end{align}
		
		For the first term on the right hand side of \eqref{eq:cond_barz}, we have
		\begin{align}
			& \norm{\bar{z}_k - \alpha_k\bar{\nabla}F(\x_k)}^2 = \norm{\bar{z}_k - \alpha_k \bar{\nabla}F(\1\bar{x}_k) + \alpha_k\left(\bar{\nabla}F(\1\bar{x}_k) - \bar{\nabla}F(\x_k)\right)}\nonumber\\
			&= \norm{\bar{z}_k - \alpha_k \bar{\nabla}F(\1\bar{x}_k)}^2 + \alpha_k^2\norm{\bar{\nabla}F(\1\bar{x}_k) - \bar{\nabla}F(\x_k)}^2\nonumber\\
			&\quad + 2\alpha_k\left\langle \bar{z}_k - \alpha_k \bar{\nabla}F(\1\bar{x}_k), \bar{\nabla}F(\1\bar{x}_k) - \bar{\nabla}F(\x_k)\right\rangle\nonumber\\
			&\leq \norm{\bar{z}_k - \alpha_k \bar{\nabla}F(\1\bar{x}_k)}^2 + \alpha_k^2\norm{\bar{\nabla}F(\1\bar{x}_k) - \bar{\nabla}F(\x_k)}^2\label{eq:fix_rate_cauchy}\\
			&\quad + 2 \alpha_k \norm{\bar{z}_k - \alpha_k \bar{\nabla}F(\1\bar{x}_k)}\cdot\norm{\bar{\nabla}F(\1\bar{x}_k) - \bar{\nabla}F(\x_k)}\nonumber\\
			&\leq (1+\gamma_1)\norm{\bar{z}_k - \alpha_k \bar{\nabla}F(\1\bar{x}_k)}^2 + (1+\frac{1}{\gamma_1})\alpha_k^2\norm{\bar{\nabla}F(\1\bar{x}_k) - \bar{\nabla}F(\x_k)}^2,\label{eq:fix_rata_inter}
		\end{align}
		where \eqref{eq:fix_rate_cauchy} comes from Cauchy-Schwarz's inequality and \eqref{eq:fix_rata_inter} holds because $2\norm{a}\norm{b} \leq \gamma_1 \norm{a}^2 + \frac{1}{\gamma_1}\norm{b}^2, \forall \gamma_1>0$. 
		Substituting \eqref{eq:fix_rata_inter} into \eqref{eq:cond_barz} and taking full expectation on both sides of \eqref{eq:cond_barz}, we get
		\begin{equation}
			\label{eq:lem1_inter}
			\E\left[\norm{\bar{z}_{k+1}}^2\right] \leq (1+\gamma_1)(1-\alpha_k\mu)^2\E\left[\norm{\bar{z}_{k}}^2\right] + (1+\gamma_1^{-1})\frac{\alpha_k^2L^2c^2\norm{U_R}^2}{n}\E\left[\norm{\check{z}_k}^2\right] + \frac{\alpha_k^2\bar{\sigma}^2}{n},
		\end{equation}
		where we invoked relations \eqref{eq:fix_rate}, \eqref{eq:bound_var},and \cref{lem:from_L}.
		Let $\gamma_1 = \frac{3}{8}\alpha_k\mu$ and notice that $\alpha_k\le \frac{1}{3\mu}$. We conclude from \eqref{eq:lem1_inter} that
		\begin{equation*}
			E\left[\norm{\bar{z}_{k+1}}^2\right] \leq (1-\frac{3}{2}\alpha_k\mu)\E\left[\norm{\bar{z}_{k}}^2\right] + \frac{3\alpha_kL^2c^2\norm{U_R}^2}{\mu n}\E\left[\norm{\check{z}_k}\right] + \frac{\alpha_k^2\bar{\sigma}^2}{n}.
		\end{equation*} 
		This finishes the proof. 
	\end{proof}

	We obtain the recursion for $\E\left[\Vert \check z_{k+1}\Vert^2\right]$ stated in \cref{lem:checkz} below. 
	\begin{lemma}\label{lem:checkz}
		Under \cref{alg:extra} with \cref{ngrad,edfi,edw,}, suppose
		\begin{equation*}
			\alpha_k \leq \frac{1-\sqrt{\lambda_2}}{2\norm{U_L}\norm{U_R}L}
		\end{equation*}
		where $\lambda_2$ is the second largest eigenvalue of $W$.
		Then for all $k\geq 1$, we have
		\begin{equation*}
			\begin{aligned}
				\E\left[\Vert \check{z}_{k+1}\Vert^2\right]&\leq 
				\frac{3+\sqrt{\lambda_2}}{4} \E\left[\Vert \check{z}_k\Vert^2\right] 
				+ \frac{\alpha_k^2\lambda_2n\Vert U_L\Vert^2}{c^2}\left(\frac{4L^2}{1-\lambda_2}\E\left[\Vert \bar z_k\Vert^2\right] + \bar\sigma^2\right)\\
				&\quad + \frac{4\Vert U_L\Vert^2\Vert V^{-}\Vert^2\Vert\nabla F(\x^*)\Vert^2}{c^2(1-\sqrt{\lambda_2})}\Vert\alpha_{k+1}-\alpha_k\Vert^2
			\end{aligned}
		\end{equation*}
	\end{lemma}
	\begin{proof}
		By squaring and taking conditional expectation on both sides of the second recursion in \cref{lem:tedy}, we obtain
		\begin{align}
			\E\left[\Vert \check z_{k+1}\Vert^2|\mathcal{F}_{k}\right]&=\Vert D_1\check z_k+\frac{\alpha_k}{c}D_1U_{L,l}(\nabla F(\x^*)-\nabla F(\x_k))\Vert^2\nonumber\\
			&\quad + \frac{\alpha_k^2\Vert D_1U_{L,l}\Vert^2}{c^2}\E\left[\Vert \s_k\Vert^2|\mathcal{F}_k\right] + \frac{1}{c^2}\Vert U_{L,r}\Vert^2\Vert \y_{k+1}^*-\y_k^*\Vert^2 \nonumber\\ 
			&\quad + 2\left\langle D_1\check z_k+\frac{\alpha_k}{c}D_1U_{L,l}(\nabla F(\x^*)-\nabla F(\x_k)), 
			\frac{1}{c}U_{L,r}(\y_{k}^*-\y_{k+1}^*)\right\rangle\nonumber\\
			&\leq (1 + \gamma)\Vert D_1\check z_k+\frac{\alpha_k}{c}D_1U_{L,l}(\nabla F(\x^*)-\nabla F(\x_k))\Vert^2 \nonumber\\
			&\quad + \frac{\alpha_k^2\norm{D_1}^2\Vert U_L\Vert^2}{c^2}\E\left[\norm{\mathbf{s}_k}^2|\mathcal{F}_k\right]
			+ \left(1 + \frac{1}{\gamma}\right)\frac{\Vert U_L\Vert^2}{c^2} \Vert \y_{k+1}^*-\y_k^*\Vert^2,\label{eq:checkz}
		\end{align}
		where $\gamma>0$ is arbitrary.
		To further bound the right hand side of \eqref{eq:checkz} in terms of $\norm{\bar z_k}^2$ and $\norm{\check{z}_k}^2$, first consider
		\begin{align}
			&\Vert D_1\check z_k+\frac{\alpha_k}{c}D_1U_{L,l}(\nabla F(\x^*)-\nabla F(\x_k))\Vert^2\nonumber\\
			&\leq (1+\gamma_2)\Vert D_1\check z_k\Vert^2 + (1+\gamma_2^{-1})\Vert \frac{\alpha_k}{c}D_1U_{L,l}(\nabla F(\x^*)-\nabla F(\x_k))\Vert^2\label{eq:checkz_rhs_split}
		\end{align}
		where $\gamma_2>0$ is arbitrary. Since $\lambda_2 := |\lambda_2(W)|$, it can be seen from \eqref{D_1} that
		\begin{equation*}
			\Vert D_1\Vert_2 = \sqrt{|\lambda_2(W)|} = \sqrt{\lambda_2}.
		\end{equation*}
		
		Take $\gamma_2 = 1/\sqrt{\lambda_2} - 1$, and notice that
		\begin{equation}
			\begin{aligned}
				\Vert \nabla F(\x^*) - \nabla F(\x_k)\Vert^2 &= \sum_{i=1}^n\Vert \nabla f_i(x_{i,k}) - \nabla f_i(x_i^*)\Vert^2\\
				&\leq L^2\sum_{i=1}^n \Vert x_{i,k} - x^*\Vert^2\\ 
				&= L^2 \Vert \tx_k\Vert^2 = L^2 \Vert \1 \bar{z}_k+c U_{R,u}\check{z}_k\Vert^2 \\
				&\leq L^2(2n\Vert \bar z_k\Vert^2 + 2c^2\Vert U_R\Vert^2\Vert \check{z}_k\Vert^2).
			\end{aligned}
		\end{equation}
		
		Relation \eqref{eq:checkz_rhs_split} leads to
		\begin{equation}\label{eq:checkz_gamma}
			\begin{aligned}
				&\Vert D_1\check z_k+\frac{\alpha_k}{c}D_1U_{L,l}(\nabla F(\x^*)-\nabla F(\x_k))\Vert^2\\
				& \leq \sqrt{\lambda_2}\Vert \check{z}_k\Vert^2 + \frac{\alpha_k^2 \lambda_2}{c^2(1-\sqrt{\lambda_2})}\|U_{L,l}\|^2\|\nabla F(\x^*)-\nabla F(\x_k)\|^2\\
				& \leq \sqrt{\lambda_2}\Vert \check{z}_k\Vert^2 + \frac{\alpha_k^2 \lambda_2}{c^2(1-\sqrt{\lambda_2})}\|U_{L,l}\|^2 L^2(2n\Vert \bar z_k\Vert^2 + 2c^2\Vert U_R\Vert^2\Vert \check{z}_k\Vert^2)\\
				&\leq \left(\sqrt{\lambda_2} + \frac{2\alpha_k^2\lambda_2\Vert U_L\Vert^2\Vert U_R\Vert^2L^2}{1-\sqrt{\lambda_2}}\right) \Vert \check{z}_k\Vert^2 
				+ \frac{2n\alpha_k^2\lambda_2\Vert U_L\Vert^2L^2}{c^2(1-\sqrt{\lambda_2})}\Vert \bar z\Vert^2.
			\end{aligned}
		\end{equation}
		
		We then bound the term $\Vert \y_{k}^* - \y_{k+1}^*\Vert$ on the right hand side of \eqref{eq:checkz}. As stated in \cref{lem:existence}, $\y_k^* = -V^{-}\alpha_kW \nabla F(\x^*)$ where $V^{-}$ is the pseudoinverse of $V$. Then, 
		\begin{align}
			\norm{\y_{k+1}^*-\y_k^*} &= \norm{V^{-}W \nabla F(\x^*) \left(\alpha_{k+1}-\alpha_k\right)}\nonumber\\
			&\leq \norm{V^{-}}\norm{\nabla F(\x^*)}\norm{\alpha_{k+1}-\alpha_k}\label{eq:ystar_res}
		\end{align}
		where \eqref{eq:ystar_res} holds since $\norm{W}= 1$.
		
		Substituting \eqref{eq:checkz_gamma} and \eqref{eq:ystar_res} into \eqref{eq:checkz}, and taking full expectation on both sides of the inequality, we obtain
		\begin{equation}
			\label{eq:checkz_inter}
			\begin{aligned}
				\E\left[\Vert \check{z}_{k+1}\Vert^2 \right]
				&\leq (1+\gamma) \left(\sqrt{\lambda_2} + \frac{2\alpha_k^2\lambda_2\Vert U_L\Vert^2\Vert U_R\Vert^2L^2}{1-\sqrt{\lambda_2}}\right) \E\left[\Vert \check{z}_k\Vert^2\right] 
				+  \frac{\alpha_k^2\lambda_2\Vert U_L\Vert^2n\bar\sigma^2}{c^2} \\
				&\quad + (1+\gamma)\frac{2n\alpha_k^2\lambda_2\Vert U_L\Vert^2L^2}{c^2(1-\sqrt{\lambda_2})}\E\left[\Vert \bar z\Vert^2\right]\\
				&\quad + (1 + \frac{1}{\gamma})\frac{\Vert U_L\Vert^2}{c^2} \Vert V^{-}\Vert^2 \Vert \nabla F(\x^*)\Vert^2 \Vert \alpha_{k+1} - \alpha_k\Vert^2.
			\end{aligned}
		\end{equation}
		
		Given that 
		$\alpha_k \leq \frac{1-\sqrt{\lambda_2}}{2\norm{U_L}\norm{U_R}L}$,
		we have
		\begin{equation*}
			\sqrt{\lambda_2}+ \frac{2\alpha_{k}^2\lambda_2\Vert U_L\Vert^2\Vert U_R\Vert^2L^2}{1-\sqrt{\lambda_2}} \leq \frac{\sqrt{\lambda_2} + 1}{2}.
		\end{equation*}
		
		Further choosing $\gamma = \frac{1-\sqrt{\lambda_2}}{2\sqrt{\lambda_2}+2}$ such that 
		\begin{equation}
			(1+\gamma)\frac{(\sqrt{\lambda_2}+1)}{2}=\frac{(\sqrt{\lambda_2}+3)}{4}.
		\end{equation}
		
		We obtain the desired result.
	\end{proof}

	\subsection{Preliminary Convergence Results}
	\label{sec:pre_r}
	We now prove the sublinear convergence rate of the EDAS algorithm under the diminishing stepsizes policy 
	\begin{equation*}
		\alpha_k = \frac{\theta}{\mu(k+m)},
	\end{equation*}
	where $\theta,m>0$ are to be determined later. For simplicity, we let 
	\begin{equation}
		\label{eq:MT}
		M_k := \E\left[\Vert \bar z_k\Vert^2\right], \quad T_k := \E\left[\Vert \check z_k\Vert^2\right], \, \forall k\ge 0.
	\end{equation}
	
	Our analysis builds upon constructing a Lyapunov function \eqref{eq:Lya}, which is inspired by the arguments in \cite{pu2019sharp}. 
	The sublinear convergence results will then be combined with \cref{lem:checkz} and \cref{barzk} respectively, so as to derive the improved convergence rates for $M_k$ and $T_k$ which will be stated in the next section.
	
	The following result comes from Lemma 4.1 in \cite{pu2019sharp} and will be used for bounding specific terms in the proofs. It will be used repeatedly in the analysis.
	\begin{lemma}\label{lem:prodineq}
		$\forall 1<a<k(a\in\mathbb{N})$ and $1<\gamma<a/2$, 
		$$
		\frac{a^{2 \gamma}}{k^{2 \gamma}} \leq \prod_{t=a}^{k-1}\left(1-\frac{\gamma}{t}\right) \leq \frac{a^{\gamma}}{k^{\gamma}}
		$$  
	\end{lemma}
	
	The Lyapunov function is defined as
	\begin{equation}
		\label{eq:Lya}
		H_k := M_k + \omega_k T_k, \,\forall k\geq 0,
	\end{equation}
	where 
	\begin{equation}
		\label{omega_k_def}
		\omega_k:= \frac{24\alpha_kc^2\Vert U_R\Vert^2L^2}{n\mu(1-\sqrt{\lambda_2})}.
	\end{equation}
	
	It will be made clear why $\omega_k$ takes this particular form in the proof of \cref{lem:hat_W}.
	\begin{lemma}\label{lem:hat_W}
		Suppose
		\begin{equation}\label{eq:K1}
			m \ge \max\left\{
			\frac{24\theta}{1-\sqrt{\lambda_2}}, \frac{24\theta L^2 \Vert U_R\Vert\Vert U_L\Vert}{\mu^2(1-\sqrt{\lambda_2})}
			\right\}, \quad \theta > 4.
		\end{equation}
		
		Under \cref{alg:extra} with \cref{ngrad,edfi,edw,}, for all $k\geq 1$, we have
		\begin{equation}
			\label{eq:Mk0}
			M_k\leq \frac{\hat H_1}{k+m}+\frac{\hat{H}_2}{(k+m)^2},
		\end{equation}
		where 
		\begin{equation*}
			\hat H_1  := \frac{1}{\theta-3}\left(p_2+\frac{p_5}{m^3}\right),\ \
			\hat H_2  := m^2 H_0+\frac{2p_3}{2\theta-3},
		\end{equation*}
		with constants $p_2,p_3,p_5$ defined in \eqref{p_form}.
		
		In addition, 
		\begin{equation}\label{eq:Tk}
			\begin{aligned}
				T_k&\leq q_0^{k}T_0 + \frac{8q_2}{(1-\sqrt{\lambda_2})(k+m)^2}+\frac{8q_3}{(1-\sqrt{\lambda_2})(k+m)^3}+ \frac{8q_4}{(1-\sqrt{\lambda_2})(k+m)^4},
			\end{aligned}
		\end{equation}
		where constants $q_0,q_2,q_3,q_4$ are defined in \eqref{q_form}.
	\end{lemma}
	\begin{proof}
		{\bf Step 1: bounding $H_k$.} From \cref{barzk,lem:checkz} and the definition of $H_{k+1}$, we have for all $k\ge 0$ that
		\begin{equation}
			\label{H_k_ineq}
			\begin{aligned}
				H_{k+1} &\leq 
				\left[(1-\frac{3}{2}\alpha_k\mu)+\omega_k\frac{4\alpha_k^2\lambda_2L^2\Vert U_L\Vert^2n}{c^2(1-\lambda_2)}\right] M_k\\
				&\quad +
				\left[\frac{3\alpha_k c^2 \Vert U_R\Vert^2 L^2}{n\mu}+\omega_k\frac{(3+\sqrt{\lambda_2})}{4}\right] T_k\\ 
				&\quad + \left(\frac{1}{n} + \frac{\omega_k \lambda_2 n\Vert U_L\Vert^2}{c^2}\right)\alpha_k^2\bar\sigma^2\\
				&\quad + \omega_k\frac{4\Vert U_L\Vert^2\Vert V^{-}\Vert^2\Vert\nabla F(\x^*)\Vert^2}{c^2(1-\sqrt{\lambda_2})}\Vert\alpha_{k+1}-\alpha_k\Vert^2.
			\end{aligned}
		\end{equation}
		
		We now show the following inequalities hold for all $k\geq 0$: 
		\begin{subequations}
			\label{UV_ineq}
			\begin{align}
				(1-\frac{3}{2}\alpha_k\mu)+\omega_k\frac{4\alpha_k^2\lambda_2L^2\Vert U_L\Vert^2n}{c^2(1-\lambda_2)} &\leq 1-\frac{4}{3}\alpha_k\mu, \label{U_ineq}\\
				\frac{3\alpha_k c^2 \Vert U_R\Vert^2 L^2}{n\mu}+\omega_k\frac{3+\sqrt{\lambda_2}}{4} &\leq (1-\frac{4}{3}\alpha_k\mu)\omega_k.\label{eq:V_ineq}
			\end{align}
		\end{subequations}
		
		First of all, condition \eqref{eq:V_ineq} is equivalent to 
		\[\left(1-\frac{4}{3}\alpha_k\mu-\frac{3+\sqrt{\lambda_2}}{4}\right)\omega_k\ge \frac{3\alpha_k c^2 \Vert U_R\Vert^2 L^2}{n\mu}.\]
		
		Since $m\geq \frac{24\theta}{1-\sqrt{\lambda_2}}$, we have $\alpha_k\le \frac{1-\sqrt{\lambda_2}}{24\mu}$ for all $k\ge 0$, and so that
		\begin{equation*}
			1-\frac{4}{3}\alpha_k\mu-\frac{3+\sqrt{\lambda_2}}{4}\geq \frac{1-\sqrt{\lambda_2}}{4}-\frac{1-\sqrt{\lambda_2}}{8}=\frac{1-\sqrt{\lambda_2}}{8}.
		\end{equation*}
		
		Hence for \eqref{eq:V_ineq} to hold, it is sufficient that
		\begin{equation}\label{eq:omegak_1}
			\omega_k\geq \frac{24\alpha_kc^2\Vert U_R\Vert^2L^2}{n\mu(1-\sqrt{\lambda_2})},
		\end{equation} 
		which is satisfied given the definition of $\omega_k$ in \eqref{omega_k_def}.
		
		Second, condition \eqref{U_ineq} requires
		\begin{equation}\label{eq:omegak_2}
			\omega_k\leq \frac{(1-\lambda_2)c^2\mu}{24\lambda_2L^2 n\Vert U_L\Vert^2}\frac{1}{\alpha_k},
		\end{equation}
		or
		\begin{equation*}
			\frac{24\alpha_kc^2\Vert U_R\Vert^2L^2}{n\mu(1-\sqrt{\lambda_2})}\leq \frac{(1-\lambda_2)c^2\mu}{24\lambda_2L^2 n\Vert U_L\Vert^2}\frac{1}{\alpha_k}.
		\end{equation*}
		
		It is sufficient that
		\begin{equation}
			\alpha_k\le \frac{\mu(1-\sqrt{\lambda_2})}{24 L^2\Vert U_R\Vert\Vert U_L\Vert},\,\forall k\ge 0\quad \text{or}\quad m\ge \frac{24\theta L^2 \Vert U_R\Vert\Vert U_L\Vert}{\mu^2(1-\sqrt{\lambda_2})}.
		\end{equation}
		
		We have verified condition \eqref{UV_ineq}, and thus obtain the following recursion for $H_{k+1}$ from \eqref{H_k_ineq}:
		\begin{equation*}
			\begin{aligned}
				H_{k+1}&\leq \left(1-\frac{4}{3}\alpha_k\mu\right)H_k +  \left(\frac{1}{n} + \frac{\omega_k \lambda_2 n\Vert U_L\Vert^2}{c^2}\right)\alpha_k^2\bar\sigma^2\\
				&\quad + \omega_k\frac{4\Vert U_L\Vert^2\Vert V^{-}\Vert^2\Vert\nabla F(\x^*)\Vert^2}{c^2(1-\sqrt{\lambda_2})}\Vert\alpha_{k+1}-\alpha_k\Vert^2\\
				&\leq \left(1-\frac{4\theta}{3(k+m)}\right)H_k + \left[\frac{24\theta L^2 \Vert U_L\Vert^2 \Vert U_R\Vert^2\lambda_2}{\mu^2(1-\sqrt{\lambda_2})(k+m)}+\frac{1}{n}\right]\frac{\theta^2\bar\sigma^2}{\mu^2(k+m)^2}\\
				&\quad + \frac{24\theta c^2\Vert U_R\Vert^2L^2}{n\mu^2(1-\sqrt{\lambda_2})(k+m)}\frac{4\Vert U_L\Vert^2\Vert V^{-}\Vert^2\Vert\nabla F(\x^*)\Vert^2}{c^2(1-\sqrt{\lambda_2})}\frac{\theta^2}{\mu^2(k+m)^4}\\
				&= \left(1-\frac{4\theta}{3(k+m)}\right)H_k+\frac{p_2}{(k+m)^2}+\frac{p_3}{(k+m)^3}+\frac{p_5}{(k+m)^5},
			\end{aligned}
		\end{equation*}
		where 
		\begin{equation}\label{p_form}
			\begin{aligned}
				p_2 & := \frac{\theta^2\bar\sigma^2}{n\mu^2},\\
				p_3 & := \frac{24\lambda_2\theta^3\bar\sigma^2\Vert U_L\Vert^2\Vert U_R\Vert^2 L^2}{\mu^4(1-\sqrt{\lambda_2})},\\
				p_5 & := \frac{96\theta^3\Vert U_L\Vert^2\Vert U_R\Vert^2L^2\Vert V^{-}\Vert^2\Vert\nabla F(\x^*)\Vert^2}{n\mu^4(1-\sqrt{\lambda_2})^2}.
			\end{aligned}
		\end{equation}
		
		Then for all $k\geq 1$, we have
		\begin{align*}
			H_k&\leq\prod_{t=0}^{k-1}\left(1-\frac{4\theta}{3(t+m)}\right)H_0\nonumber\\
			&\quad +\sum_{t=0}^{k-1}\left[\prod_{j=t+1}^{k-1}\left(1-\frac{4\theta}{3(j+m)}\right)\right]\left[\frac{p_2}{(t+m)^2}+\frac{p_3}{(t+m)^3}+\frac{p_5}{(t+m)^5}\right]\nonumber\\
			&\leq \frac{m^{\frac{4\theta}{3}}}{(k+m)^{\frac{4\theta}{3}}}H_0 + \sum_{t=0}^{k-1}\frac{(m+t+1)^{\frac{4\theta}{3}}}{(k+m)^{\frac{4\theta}{3}}}\left[\frac{p_2}{(t+m)^2}+\frac{p_3}{(t+m)^3}+\frac{p_5}{(t+m)^5}\right] \\
			&\nonumber\leq \frac{m^{\frac{4\theta}{3}}}{(k+m)^{\frac{4\theta}{3}}}H_0 \\
			&\nonumber\quad + \frac{4}{3(k+m)^{\frac{4\theta}{3}}}\sum_{t=0}^{k-1}\left[p_2(m+t)^{\frac{4\theta}{3}-2} + p_3(m+t)^{\frac{4\theta}{3}-3} + p_5(m+t)^{\frac{4\theta}{3}-5}\right],
		\end{align*}
		where we invoked \cref{lem:prodineq} for the second inequality.
		Additionally, since $\theta>4$, we have
		\begin{equation}\label{eq:sum_t}
			\sum_{t=0}^{k-1}(m+t)^{\frac{4\theta}{3}-i}\leq \int_{-1}^k (m+t)^{\frac{4\theta}{3}-i}\mathrm{dt}\leq \frac{3}{4\theta-3i+3}(m+k)^{\frac{4\theta}{3}-i+1},\, i=2,3,5.
		\end{equation}
		
		Further noticing that $\frac{m^{\frac{4\theta}{3}}}{(k+m)^{\frac{4\theta}{3}}}\leq \frac{m^2}{(k+m)^2}$ since $\frac{m}{k+m}\leq 1$, we have
		\begin{align*}
			H_k\leq & \frac{m^2}{(k+m)^2} H_0 + \frac{4p_2}{4\theta-3}\frac{1}{k+m} 
			+\frac{4p_3}{4\theta-6}\frac{1}{(m+k)^2}+\frac{4p_5}{4\theta-12}\frac{1}{(m+k)^4}\\
			= & \frac{1}{k+m}\left[\frac{4p_2}{4\theta-3}+\frac{p_5}{\theta-3}\frac{1}{(m+k)^3}\right] +\frac{1}{(k+m)^2}\left[m^2 H_0+\frac{2p_3}{2\theta-3}\right]\\
			\le & \frac{1}{k+m}\left[\frac{1}{\theta-3}\left(p_2+\frac{p_5}{m^3}\right)\right]+\frac{1}{(k+m)^2}\left[m^2 H_0+\frac{2p_3}{2\theta-3}\right] \\\
			= & \frac{\hat{H}_1}{k+m}+\frac{\hat{H}_2}{(k+m)^2}.
		\end{align*}
		
		Since $H_k=M_k+\omega_kT_k\ge M_k$, we get
		\begin{equation}
			\label{eq:hatW_Mk_bound}
			M_k\leq\frac{\hat H_1}{k+m}+\frac{\hat{H}_2}{(k+m)^2},\; \forall k.
		\end{equation}
		
		{\bf Step 2: bounding $T_k$.} To bound $T_k$, we substitute \eqref{eq:hatW_Mk_bound} into \cref{lem:checkz} and get
		\begin{equation*}
			\begin{aligned}
				T_{k+1}&\leq 
				\frac{4\lambda_2n \Vert U_L\Vert^2 \theta^2L^2}{\mu^2c^2(1-\lambda_2)}\left(\frac{\hat{H}_1}{(k+m)^3}+\frac{\hat{H}_2}{(k+m)^4}\right) + \frac{\lambda_2n\Vert U_L\Vert^2\bar\sigma^2\theta^2}{c^2\mu^2}\frac{1}{(k+m)^2}\\
				&\quad + \frac{(3+\sqrt{\lambda_2})}{4} T_k +  \frac{4\Vert U_L\Vert^2\Vert V^{-}\Vert^2\Vert\nabla F(\x^*)\Vert^2\theta^2}{\mu^2c^2(1-\sqrt{\lambda_2})}\frac{1}{(k+m)^4}\\
				& = q_0T_k+\frac{q_2}{(k+m)^2}+\frac{q_3}{(k+m)^3}+\frac{q_4}{(k+m)^4},
			\end{aligned}
		\end{equation*}
		where 
		\begin{equation}\label{q_form}
			\begin{alignedat}{2}
				q_0 & := \frac{3+\sqrt{\lambda_2}}{4},&\quad
				q_2 & := \frac{\lambda_2n\Vert U_L\Vert^2\bar\sigma^2\theta^2}{c^2\mu^2},\\
				q_3 & := \frac{4\lambda_2n \Vert U_L\Vert^2 \theta^2L^2\hat H_1}{\mu^2c^2(1-\lambda_2)},&\quad
				q_4 & := \frac{4\lambda_2n \Vert U_L\Vert^2 \theta^2L^2\hat H_2}{\mu^2c^2(1-\lambda_2)}+\frac{4\Vert U_L\Vert^2\Vert V^{-}\Vert^2\Vert\nabla F(\x^*)\Vert^2\theta^2}{\mu^2c^2(1-\sqrt{\lambda_2})}.
			\end{alignedat}
		\end{equation}
		
		Hence for all $k\ge 0$, we have
		\begin{equation*}
			\begin{aligned}
				T_k
				&\leq q_0^{k}T_0 + \sum_{t=0}^{k-1} q_0^{k-1-t}\left[\frac{q_2}{(m+t)^2}+\frac{q_3}{(m+t)^3}+\frac{q_4}{(m+t)^4}\right]\\
				&:= q_0^{k}T_0 + A_2(k) + A_3(k)+A_4(k),
			\end{aligned}
		\end{equation*}
		with
		\begin{equation}\label{A_form}
			\begin{alignedat}{2}
				A_2(k) & := \sum_{t=0}^{k-1}\frac{q_0^{k-t-1}q_2}{(m+t)^2},&\quad
				A_3(k) & := \sum_{t=0}^{k-1}\frac{q_0^{k-t-1}q_3}{(m+t)^3},\\
				A_4(k) & := \sum_{t=0}^{k-1}\frac{q_0^{k-t-1}q_4}{(m+t)^4}.
			\end{alignedat}
		\end{equation}
		
		Note that for $i=2,3,4$, 
		\begin{equation*}
			A_i(k+1) = q_0\left[\sum_{t=0}^{k-1}q_0^{k-t-1}\frac{q_i}{(m+t)^i}+\frac{q_i}{q_0(m+k)^i}\right]= q_0A_i(k)+\frac{q_i}{(m+k)^i},\; \forall k\ge0.
		\end{equation*}
		
		By induction we obtain
		\begin{equation}\label{eq:A_k_ineq}
			A_i(k)\leq \frac{1}{(k+m)^i}\frac{q_i}{(1-\frac{1}{m+1})^i-q_0},\; i=2,3,4.
		\end{equation}
		
		Since $m\geq \frac{24\theta}{1-\sqrt{\lambda_2}}$, it can be verified that
		\begin{equation}\label{eq:A_k}
			\left(1-\frac{1}{m+1}\right)^i-q_0\geq \frac{1-q_0}{2}=\frac{1-\sqrt{\lambda_2}}{8},\; i=2,3,4.
		\end{equation}
		
		We conclude that for all $k \ge0$, there holds
		\begin{equation}
			\begin{aligned}
				T_k&\leq q_0^{k}T_0 + \frac{8q_2}{(1-\sqrt{\lambda_2})(k+m)^2}+\frac{8q_3}{(1-\sqrt{\lambda_2})(k+m)^3} +\frac{8q_4}{(1-\sqrt{\lambda_2})(k+m)^4}.
			\end{aligned}
		\end{equation}
	\end{proof}
	
	\section{Main Results}
	\label{sec:main_res}
	In this section, we present the main convergence results for EDAS. First, we show EDAS performs asymptotically as well as the centralized stochastic gradient descent (SGD) method in \cref{thm:Uk}. In other words, EDAS enjoys the so-called `` asymptotic network independence'' property. Then we refine the convergence results by iteratively combining the existing results with \cref{barzk,lem:checkz}. 
	Finally, we derive the transient time for EDAS to approach the network independent rate in \cref{thm:transient_time} by comparing with the performance of centralized SGD.
	
	\begin{theorem}\label{thm:Uk}
		Under \cref{alg:extra} with \cref{ngrad,edfi,edw,}, suppose $\theta>4$ and $m$ is chosen according to \eqref{eq:K1}. We have 
		\begin{equation*}
			\begin{aligned}
				M_k&\leq  \frac{4\theta^2\bar\sigma^2}{(3\theta-2)n\mu^2(k+m)} + \frac{4c_0 T_0m^{\frac{3\theta}{2}-1}}{(1-q_0)(k+m)^{\frac{3\theta}{2}}} + \frac{m^{\frac{3\theta}{2}}M_0}{(k+m)^{\frac{3\theta}{2}}} \\
				&\quad + \frac{32c_0}{1-\sqrt{\lambda_2}}\left[\frac{q_2}{(3\theta-4)(m+k)^2} + 
				\frac{q_3}{(3\theta-6)(k+m)^3}+ \frac{q_4}{(3\theta-8)(k+m)^4}\right],
			\end{aligned}
		\end{equation*}
		where $c_0:= \frac{3\theta L^2 c^2 \Vert U_R\Vert^2}{n\mu^2}$
	\end{theorem}
	\begin{proof}
		From \cref{barzk}, we have
		\begin{equation*}
			M_{k+1}\leq 
			\left(1-\frac{3\alpha_k\mu}{2}\right)M_k+ \frac{3\alpha_kL^2c^2\Vert U_R\Vert^2}{n\mu}T_k+ \frac{\alpha_k^2\bar\sigma^2}{n}.
		\end{equation*}
		
		Thus for all $k\ge 0$,
		\begin{equation}\label{eq:rec_M_1}
			\begin{aligned}
				M_k&\leq \prod_{t=0}^{k-1}\left(1-\frac{3\theta}{2(t+m)}\right)M_0\\
				&\quad +\sum_{t=0}^{k-1}\left[\prod_{j=t+1}^{k-1}\left(1-\frac{3\theta}{2(j+m)}\right)\right]\left[\frac{3\theta L^2 c^2 \Vert U_R\Vert^2}{n\mu^2(t+m)}T(t)
				+\frac{\theta^2\bar\sigma^2}{n\mu^2(t+m)^2}\right]\\
				&\leq \left(\frac{m}{k+m}\right)^{\frac{3\theta}{2}}M_0 + \sum_{t=0}^{k-1}\left(\frac{t+m+1}{k+m}\right)^{\frac{3\theta}{2}}\left[\frac{c_0 T(t)}{t+m}
				+\frac{\theta^2\bar\sigma^2}{n\mu^2(t+m)^2}\right],
			\end{aligned}
		\end{equation}
		where \cref{lem:prodineq} was invoked for obtaining the second inequality. In light of \eqref{eq:Tk},
		\begin{equation*}
			\begin{aligned}
				&\quad \sum_{t=0}^{k-1}\left(\frac{t+m+1}{k+m}\right)^{\frac{3\theta}{2}}\left[\frac{c_0 T(t)}{t+m}+ \frac{\theta^2\bar\sigma^2}{n\mu^2(t+m)^2}\right]\\
				&\leq 
				c_0T_0\sum_{t=0}^{k-1}\frac{(m+t+1)^{\frac{3\theta}{2}}q_0^{t}}{(m+t)(k+m)^{\frac{3\theta}{2}}} + \frac{\theta^2\bar\sigma^2}{n\mu^2(k+m)^{\frac{3\theta}{2}}}\sum_{t=0}^{k-1}\frac{(m+t+1)^{\frac{3\theta}{2}}}{(m+t)^2}\\
				&\quad + \frac{8c_0}{(1-\sqrt{\lambda_2})(k+m)^{\frac{3\theta}{2}}}\sum_{t=0}^{k-1}\biggl[\frac{q_2(m+t+1)^{\frac{3\theta}{2}}}{(m+t)^3} + \frac{q_3(m+t+1)^{\frac{3\theta}{2}}}{(m+t)^4}\\
				&\quad  +\frac{q_4(m+t+1)^{\frac{3\theta}{2}}}{(m+t)^5}\biggl]\\
				&\leq 
				2c_0T_0\sum_{t=0}^{k-1}\frac{(m+t)^{\frac{3\theta}{2}-1}q_0^{t}}{(k+m)^{\frac{3\theta}{2}}} +\frac{2\theta^2\bar\sigma^2}{n\mu^2(k+m)^{\frac{3\theta}{2}}}\sum_{t=0}^{k-1}(m+t)^{\frac{3\theta}{2}-2}\\
				&\quad + \frac{16c_0}{(1-\sqrt{\lambda_2})(k+m)^{\frac{3\theta}{2}}}\sum_{t=0}^{k-1}\left[q_2(m+t)^{\frac{3\theta}{2}-3} 
				+q_3(m+t)^{\frac{3\theta}{2}-4} + q_4(m+t)^{\frac{3\theta}{2}-5}\right].
			\end{aligned}
		\end{equation*}
		
		Since $(m+t)^{\frac{3\theta}{2}-1}q_0^{t}$ is decreasing in $t$ and
		\begin{align*}
			\frac{1}{\ln q_0}d((m+t)^{\frac{3\theta}{2}-1}q_0^{t})& = \frac{1}{\ln q_0}\left(\frac{3\theta}{2}-1\right)(m+t)^{\frac{3\theta}{2}-2}q_0^{t}dt+(m+t)^{\frac{3\theta}{2}-1}q_0^{t} dt\\
			& \ge \frac{1}{2} (m+t)^{\frac{3\theta}{2}-1}q_0^{t} dt,
		\end{align*}
		we have
		\begin{equation*}
			\begin{aligned}
				\sum_{t=0}^{k-1}(m+t)^{\frac{3\theta}{2}-1}q_0^{t}&\leq 
				\int_{0}^{+\infty}(m+t)^{\frac{3\theta}{2}-1}q_0^{t}dt\\
				&\leq \frac{2}{\ln q_0}\int_{m}^{+\infty}\frac{d(t^{\frac{3\theta}{2}-1}q_0^{t-m})}{dt}\\
				&=-\frac{2m^{\frac{3\theta}{2}-1}}{\ln q_0}\leq \frac{2m^{\frac{3\theta}{2}-1}}{1-q_0}.
			\end{aligned}
		\end{equation*}
		
		In addition, from an argument similar to \eqref{eq:sum_t}, we obtain
		\begin{equation*}
			\sum_{t=0}^{k-1}(m+t)^{\frac{3\theta}{2}-i}\leq \frac{1}{\frac{3\theta}{2}-i+1}(m+k)^{\frac{3\theta}{2}-i+1},\; i=2,3,4,5.
		\end{equation*}
		
		Hence for all $k\ge 0$, 
		\begin{equation*}
			\begin{aligned}
				M_k&\leq  \frac{4\theta^2\bar\sigma^2}{(3\theta-2)n\mu^2(k+m)} + \frac{4c_0 T_0m^{\frac{3\theta}{2}-1}}{(1-q_0)(k+m)^{\frac{3\theta}{2}}} + \frac{m^{\frac{3\theta}{2}}M_0}{(k+m)^{\frac{3\theta}{2}}} \\
				&\quad + \frac{32c_0}{1-\sqrt{\lambda_2}}\left[\frac{q_2}{(3\theta-4)(m+k)^2} + 
				\frac{q_3}{(3\theta-6)(k+m)^3}+ \frac{q_4}{(3\theta-8)(k+m)^4}\right].
			\end{aligned}
		\end{equation*}
	\end{proof}
	
	From \cref{thm:Uk}, we can see that the asymptotic convergence rate of EDAS, which behaves as $\frac{4\theta^2\bar\sigma^2}{(3\theta-2)n\mu^2(k+m)}$, is of the same order as that of the centralized SGD method \cite{pu2019sharp}. Our next goal is to derive the transient time needed for EDAS to reach this asymptotic convergence rate. In particular, we would like to characterize the dependence relationship between the transient time and the network characteristics as well as function properties. For this purpose, we first simplify the presentation of \cref{thm:Uk} with $\mo{(\cdot)}$ notation in  \cref{cor:rec_MT_2} and then improve the convergence results by iteratively substituting the existing upper bounds into \cref{barzk,lem:checkz}.
	
	In the following lemma, we estimate the constants $q_i, i = 2, 3, 4$, $M_0, T_0, H_0$, and $c_0$ appearing in the statement of \cref{thm:Uk} with $\mo{(\cdot)}$ notation. 
	\begin{lemma}\label{lem:bigO}
		Let the free scaling parameter $c$ in \cref{lem:B_decomposition} be  chosen as
		\begin{equation*}
			c^2 = n\Vert U_L\Vert^2.
		\end{equation*}
		Then we have
		\begin{align*}
			& M_0 = \mbb{\frac{\norm{\x_0-\x^*}^2}{n}},\quad T_0 = \mbb{\frac{\norm{\x_0-\x^*}^2 + (1-\lambda_2)\norm{\nabla F(\x^*)}^2}{n}}\\
			& H_0 = \mbb{\frac{\norm{\x_0-\x^*}^2 + (1-\lambda_2)\norm{\nabla F(\x^*)}^2}{n}}, \quad c_0 = \mbb{1}\\
			& \frac{1}{1-q_0} = \mbb{\frac{1}{1-\lambda_2}}, \quad q_2 = \mbb{1}\\
			& q_3 = \mbb{\frac{1+\norm{\nabla F(\x^*)}^2}{n(1-\lambda_2)}},\quad 
			q_4 =\mbb{\frac{\norm{\x_0-\x^*}^2}{n(1-\lambda_2)^3}+\frac{\norm{\nabla F(\x^*)}^2}{n(1-\lambda_2)^2}+\frac{1}{(1-\lambda_2)^2}}.
		\end{align*}
	\end{lemma}
	\begin{proof}
		See \cref{sec:prf_bigO}.
	\end{proof}

	In light of \cref{lem:bigO}, we again iteratively apply \cref{barzk,lem:checkz} to derive \cref{cor:rec_MT_2} based on \cref{thm:Uk}.
	\begin{corollary}\label{cor:rec_MT_2}
		Under \cref{alg:extra} with \cref{ngrad,edfi,edw,}, suppose $\theta>4$ and $m$ be chosen  according to \eqref{eq:K1}. We have for all $k\geq 0$, 
		\begin{equation*}
			\begin{aligned}
				M_k&\leq \frac{4\theta^2\bar\sigma^2}{(3\theta-2)n\mu^2(k+m)} + \mo\left(\frac{1}{1-\lambda_2}\right)\frac{1}{(k+m)^2}\\
				&\quad + \mbb{\frac{\norm{\nabla F(\x^*)}^2}{n(1-\lambda_2)^3}}\frac{1}{(k+m)^4} + \mbb{\frac{\norm{\x_0-\x^*}^2}{n(1-\lambda_2)^6}}\frac{1}{(k+m)^6}.
			\end{aligned}
		\end{equation*}
		and
		\begin{equation*}
			\begin{aligned}
				T_k &\leq   \mbb{\frac{\norm{\x_0-\x^*}^2+\norm{\nabla F(\x^*)}^2}{n}}\left(\frac{3+\sqrt{\lambda_2}}{4}\right)^k + \mbb{\frac{1}{1-\lambda_2}}\frac{1}{(k+m)^2}\\
				& \quad + \mbb{\frac{\norm{\nabla F(\x^*)}^2}{n(1-\lambda_2)^3}}\frac{1}{(k+m)^4} + \mbb{\frac{\norm{\x_0-\x^*}^2}{n(1-\lambda_2)^6}}\frac{1}{(k+m)^6}.
			\end{aligned}
		\end{equation*}
	\end{corollary}
	\begin{proof}
		From \cref{thm:Uk}, noticing that $m=\Omega({\frac{1}{1-\lambda_2}})$, we have when $k\geq 0$,
		\begin{equation}
			\label{eq:Mk_sim}
			\begin{aligned}
				M_k&\leq  \frac{4\theta^2\bar\sigma^2}{(3\theta-2)n\mu^2(k+m)} + \frac{4c_0 T_0m^{\frac{3\theta}{2}-1}}{(1-q_0)(k+m)^{\frac{3\theta}{2}}} + \frac{m^{\frac{3\theta}{2}}M_0}{(k+m)^{\frac{3\theta}{2}}} \\
				&\quad + \frac{32c_0}{1-\sqrt{\lambda_2}}\left[\frac{q_2}{(3\theta-4)(m+k)^2} + 
				\frac{q_3}{(3\theta-6)(k+m)^3}+ \frac{q_4}{(3\theta-8)(k+m)^4}\right]\\
				&= \frac{4\theta^2\bar\sigma^2}{(3\theta-2)n\mu^2(k+m)} + \mbb{\frac{\norm{\x_0-\x^*}^2 + (1-\lambda_2)\norm{\nabla F(\x^*)}^2}{n(1-\lambda_2)^{\frac{3\theta}{2}}}}\frac{1}{(k+m)^{\frac{3\theta}{2}}}\\
				&\quad + \mbb{\frac{1}{1-\lambda_2}}\frac{1}{(k+m)^2} + \mbb{\frac{1+\norm{\nabla F(\x^*)}^2}{n(1-\lambda_2)^{2}}}\frac{1}{(k+m)^3}\\
				&\quad + \mbb{\frac{\norm{\x_0-\x^*}^2 }{n(1-\lambda_2)^4}+\frac{\norm{\nabla F(\x^*)}^2}{n(1-\lambda_2)^3}+\frac{1}{(1-\lambda_2)^3}}\frac{1}{(k+m)^4}\\
				&= \frac{4\theta^2\bar\sigma^2}{(3\theta-2)n\mu^2(k+m)} + \mbb{\frac{1}{1-\lambda_2}}\frac{1}{(k+m)^2}\\
				&\quad + \mbb{\frac{\norm{\nabla F(\x^*)}^2}{n(1-\lambda_2)^{2}}}\frac{1}{(k+m)^3} + \mbb{\frac{\norm{\x_0-\x^*}^2 }{n(1-\lambda_2)^4}}\frac{1}{(k+m)^4}\\
				& = \frac{4\theta^2\bar\sigma^2}{(3\theta-2)n\mu^2(k+m)} + \frac{c_2}{(k+m)^2} + \frac{c_3}{(k+m)^3} + \frac{c_4}{(k+m)^4},
			\end{aligned}
		\end{equation}
		where 
		\begin{equation*}
			c_2 = \mo\left(\frac{1}{1-\lambda_2}\right),\quad c_3 = \mbb{\frac{\norm{\nabla F(\x^*)}^2}{n(1-\lambda_2)^{2}}},\quad c_4 = \mbb{\frac{\norm{\x_0-\x^*}^2 }{n(1-\lambda_2)^4}}.
		\end{equation*}
		
		Substituting \eqref{eq:Mk_sim} into \cref{lem:checkz}, we have
		\begin{equation*}
			\begin{aligned}
				T_{k+1}&\leq q_0 T_k + \frac{\theta^2\lambda_2}{\mu^2(k+m)^2}\biggl[\frac{4L^2}{1-\lambda_2}\biggl(\frac{4\theta^2\bar\sigma^2}{(3\theta-2)n\mu^2(k+m)} + \frac{c_2}{(k+m)^2}\\ 
				&\quad + \frac{c_3}{(k+m)^3} + \frac{c_4}{(k+m)^4}\biggr) +\bar\sigma^2\biggr] + \frac{4\theta^2\Vert V^{-}\Vert^2\Vert\nabla F(\x^*)\Vert^2}{n\mu^2(1-\sqrt{\lambda_2})}\frac{1}{(k+m)^4}\\
				&=  q_0 T_k +  \frac{\theta^2\lambda_2\bar\sigma^2}{\mu^2}\frac{1}{(k+m)^2} + \frac{16\theta^4\bar\sigma^2\lambda_2 L^2}{(3\theta-2)n\mu^4(1-\lambda_2)}\frac{1}{(k+m)^3}\\
				&\quad  + \biggl[\frac{4\theta^2\Vert V^{-}\Vert^2\Vert\nabla F(\x^*)\Vert^2}{n\mu^2(1-\sqrt{\lambda_2})}+ \frac{4\theta^2\lambda_2c_2 L^2}{\mu^2(1-\lambda_2)}\biggr]\frac{1}{(k+m)^4}\\ 
				&\quad +  \frac{4\theta^2\lambda_2L^2c_3}{\mu^2(1-\lambda_2)}\frac{1}{(k+m)^5} + \frac{4\theta^2\lambda_2L^2c_4}{\mu^2(1-\lambda_2)}\frac{1}{(k+m)^6}\\
				&= q_0 T_k + \frac{b_2}{(k+m)^2} 
				+ \frac{b_4}{(k+m)^4} + \frac{b_6}{(k+m)^6},
			\end{aligned}
		\end{equation*}
		where 
		\begin{equation*}
			b_2 = \mo(1), \quad
			b_4 = \mbb{\frac{\norm{\nabla F(\x^*)}^2}{n(1-\lambda_2)^2}}, \quad
			b_6 = \mbb{\frac{\norm{\x_0-\x^*}^2}{n(1-\lambda_2)^5}}.
		\end{equation*}
		
		Then, 
		\begin{equation*}
			\begin{aligned}
				T_k&\leq q_0^{k}T_0
				+ \sum_{t=0}^{k-1}q_0^{k-t-1}\left[\sum_{i=2,4,6}\frac{b_i}{(k+m)^i}\right].
			\end{aligned}
		\end{equation*}
		
		Similar to the derivation of \eqref{eq:A_k_ineq}, we obtain that for all $k\ge 0$,
		\begin{equation}
			\label{eq:Tk_sim}
			\begin{aligned}
				T_k &\leq q_0^{k}T_0 + \frac{8}{1-\sqrt{\lambda_2}}\left[\frac{b_2}{(k+m)^2} + \frac{b_4}{(k+m)^4} +  \frac{b_6}{(k+m)^6}\right]\\
				&\leq q_0^{k}T_0 + \mbb{\frac{1}{1-\lambda_2}}\frac{1}{(k+m)^2}+ \mbb{\frac{\norm{\nabla F(\x^*)}^2}{n(1-\lambda_2)^3}}\frac{1}{(k+m)^4}\\
				& \quad + \mbb{\frac{\norm{\x_0-\x^*}^2}{n(1-\lambda_2)^6}}\frac{1}{(k+m)^6}\\
				& = q_0^{k}T_0 + \sum_{i=2,4,6}\frac{e_i}{(k+m)^i},
			\end{aligned}
		\end{equation}
		where each $e_i$ is of the same order of $\frac{b_i}{(1-\lambda_2)}$ for $i=2,4,6$.
		
		Substituting \eqref{eq:Tk_sim} into \cref{barzk}, we obtain
		\begin{equation}
			M_{k+1}\leq \left(1-\frac{3\theta}{2(k+m)}\right)M_k + \frac{\theta^2\bar\sigma^2}{n\mu^2(k+m)^2}
			+ \frac{c_0}{k+m}\biggl[q_0^{k}T_0 + \sum_{i=2,4,6}\frac{e_i}{(k+m)^i}\biggl].
		\end{equation}
		
		Similar to relation \eqref{eq:rec_M_1}, we derive
		\begin{equation*}
			\begin{aligned}
				M_k&\leq \frac{m^{\frac{3\theta}{2}}}{(k+m)^{\frac{3\theta}{2}}}M_0 + \sum_{t=0}^{k-1}\frac{(m+t+1)^{\frac{3\theta}{2}}}{(k+m)^{\frac{3\theta}{2}}}
				\left[\frac{c_0q_0^{t}T_0}{t+m}+ \sum_{i=2,4,6}\frac{c_0e_i}{(t+m)^{(i+1)}}\right]\\
				&\quad + \frac{\theta^2\bar\sigma^2}{n\mu^2(k+m)^{\frac{3\theta}{2}}}\sum_{t=0}^{k-1}\frac{(m+t+1)^{\frac{3\theta}{2}}}{(t+m)^2}\\
				&\leq \frac{4\theta^2\bar\sigma^2}{(3\theta-2)n\mu^2(k+m)}+ \frac{m^{\frac{3\theta}{2}}}{(k+m)^{\frac{3\theta}{2}}}M_0 + \frac{2c_0 T_0m^{\frac{3\theta}{2}-1}}{(1-q_0)(k+m)^{\frac{3\theta}{2}}}\\
				&\quad + \frac{4c_0e_2}{(3\theta-4)(k+m)^2} + \frac{4c_0e_4}{(3\theta-8)(k+m)^4} + \frac{4c_0e_6}{(3\theta-12)(k+m)^6}\\
				&= \frac{4\theta^2\bar\sigma^2}{(3\theta-2)n\mu^2(k+m)} + \mo\left(\frac{1}{1-\lambda_2}\right)\frac{1}{(k+m)^2}\\
				&\quad + \mbb{\frac{\norm{\nabla F(\x^*)}^2}{n(1-\lambda_2)^3}}\frac{1}{(k+m)^4} + \mbb{\frac{\norm{\x_0-\x^*}^2}{n(1-\lambda_2)^6}}\frac{1}{(k+m)^6}.
			\end{aligned}
		\end{equation*}
		Again we considered the fact that $(k+m)\geq \Omega{(\frac{1}{1-\lambda_2})}$ with $m$ chosen according to \eqref{eq:K1} for all $k\geq1$.
	\end{proof}
	
	We now state the convergence rate for the EDAS algorithm under the measure $\frac{1}{n}\sum_{i=1}^n\E\left[\Vert x_{i,k}-x^*\Vert^2\right]=\E\left[\Vert \tx_{k}\Vert^2 \right]$, i.e., the average expected optimization error over all the agents in the network.
	Note that from \cref{cor:rec_MT_2}, the error $\tx_k$ can be decomposed into $\tx_k=\mathbf{1}\bar z_k+cU_{R,u}\check{z}_k$. Therefore,
	\begin{equation}\label{eq:rel_SX}
		\begin{aligned}
			\Vert \tx_k\Vert^2 &\leq 2n\Vert \bar z_k\Vert^2+2c^2\Vert U_{R}\Vert^2\Vert\check{z}_k\Vert^2.
		\end{aligned}
	\end{equation} 
	
	Taking full expectation on both sides of \eqref{eq:rel_SX} and noting the choice $c^2= n\norm{U_L}^2$, we obtain
	\begin{equation}\label{eq:rel}
		\begin{aligned}
			\E\left[\Vert \tx_{k}\Vert^2 \right]\leq 2n M_k+2n\Vert U_L\Vert^2\Vert U_R \Vert^2 T_k.\\
		\end{aligned}
	\end{equation}
	
	Combining \cref{cor:rec_MT_2} and \eqref{eq:rel} leads to the improved convergence rate of EDAS stated in \cref{thm:mag_extra}.
	\begin{theorem}\label{thm:mag_extra}
		Under \cref{alg:extra} with \cref{ngrad,edfi,edw,}, suppose $\theta>5$, and $m$ is chosen as in \eqref{eq:K1}. We have for all $k\geq 1$, 
		\begin{multline*}
			\frac{1}{n}\sum_{i=1}^n\E\left[\Vert x_{i,k}-x^*\Vert^2\right]
			\leq \frac{4\theta^2\bar\sigma^2}{(3\theta-2)n\mu^2(k+m)} + \mo\left(\frac{1}{1-\lambda_2}\right)\frac{1}{(k+m)^2}\\
			\quad + \mbb{\frac{\norm{\x_0-\x^*}^2 + \norm{\nabla F(\x^*)}^2}{n}}\left(\frac{3+\lambda_2}{4}\right)^k\\
			\quad +\mbb{\frac{\norm{\nabla F(\x^*)}^2}{n(1-\lambda_2)^3}}\frac{1}{(k+m)^4}
			+ \mbb{\frac{\norm{\x_0-\x^*}^2}{n(1-\lambda_2)^6}}\frac{1}{(k+m)^6}.
		\end{multline*}
	\end{theorem}
	
	\begin{proof}
		Substitute the bounds on $M_k$ and $T_k$ in \cref{cor:rec_MT_2} into \eqref{eq:rel}, and notice that the product $\norm{U_R}^2\norm{U_L}^2$ is independent of $n$ according to \cref{lem:B_decomposition}. We obtain the result. 
	\end{proof}
	
	\subsection{Transient Time}
	In this part, we estimate how long it takes for EDAS to achieve the convergence rate of centralized stochastic gradient descent (SGD) method, i.e., \emph{transient time} of the algorithm. First, we state the convergence rate of SGD \cite{pu2019sharp}. 
	\begin{theorem}\label{thm:mag_cen}
		Under the centralized stochastic gradient descent method (SGD) with stepsize policy $\alpha_k= \frac{\theta}{\mu(k+m)}$, suppose $m\geq \left\lceil \frac{\theta L}{\mu}\right\rceil$. We have
		\begin{equation*}
			\E\left[\Vert x_k- x^*\Vert^2\right]\leq \frac{\theta^2\bar\sigma^2}{(2\theta-1)n\mu^2(k+m)}+\mathcal{O}\left(\frac{1}{n}\right)\frac{1}{(k+m)^2}.
		\end{equation*}
	\end{theorem}
	
	In the next theorem, we derive the transient time for EDAS.
	\begin{theorem}\label{thm:transient_time}
		Under \cref{alg:extra} with \cref{ngrad,edfi,edw,}, suppose $\theta>5$ and $m$ is chosen according to \eqref{eq:K1}. Then it takes 
		\begin{multline*}
			K_T = \max\left\{\mbb{\frac{n}{1-\lambda_2}}, \mbb{\left(\frac{\norm{\nabla F(\x^*)}^2}{(1-\lambda_2)^3}\right)^{\frac{1}{3}}},\mbb{\left(\frac{\norm{\x_0-\x^*}^2}{(1-\lambda_2)^6}\right)^{\frac{1}{5}}},\right.\\
			\left.\mbb{\frac{\max\left\{\log \left(\norm{\x_0-\x^*}^2 + \norm{\nabla F(\x^*)}^2\right),-\log(1-\lambda_2)\right\}}{1-\lambda_2}}\right\}
		\end{multline*}
		time for \cref{alg:extra} to reach the asymptotic, network independent convergence rate, that is, when $k\geq K_T$, we have $\frac{1}{n}\sum\limits_{i=1}^n\E\left[\Vert x_{i,k}-x^*\Vert^2\right]\leq \frac{\theta^2\bar\sigma^2}{(2\theta-1)n\mu^2k}\mathcal{O}(1)$. 
	\end{theorem}
	\begin{proof}
		From \cref{thm:mag_extra}, we have
		\begin{equation*}
			\begin{aligned}
				\frac{1}{n}\sum_{i=1}^n\E\left[\Vert x_{i,k}-x^*\Vert^2\right]
				&\leq \frac{\theta^2\bar\sigma^2}{(2\theta-1)n\mu^2 k}\biggl[ \frac{4(2\theta-1)}{3\theta-2}
				+ \mbb{\frac{n}{1-\lambda_2}}\frac{1}{k}\\
				&\quad + \mbb{\norm{\x_0-\x^*}^2 + \norm{\nabla F(\x^*)}^2}\left(\frac{3+\lambda_2}{4}\right)^k k\\
				&\quad + \mbb{\frac{\norm{\nabla F(\x^*)}^2}{(1-\lambda_2)^3}}\frac{1}{k^3}+\mbb{\frac{\norm{\x_0-\x^*}^2}{(1-\lambda_2)^6}}\frac{1}{k^5}
				\biggl].
			\end{aligned}
		\end{equation*}
		
		Let $K_T$ be such that 
		\begin{equation*}
			\begin{aligned}
				\mo(1) &= \frac{4(2\theta-1)}{3\theta-2}
				+ \mbb{\frac{n}{1-\lambda_2}}\frac{1}{K_T} \\
				&\quad + \mbb{\norm{\x_0-\x^*}^2 + \norm{\nabla F(\x^*)}^2}\left(\frac{3+\lambda_2}{4}\right)^{K_T}K_T\\
				&\quad +\mbb{\frac{\norm{\nabla F(\x^*)}^2}{(1-\lambda_2)^3}}\frac{1}{K_T^3}+ \mbb{\frac{\norm{\x_0-\x^*}^2}{(1-\lambda_2)^6}}\frac{1}{K_T^5}.
			\end{aligned}
		\end{equation*}
		
		We obtain
		\begin{multline*}
			K_T = \max\left\{\mbb{\frac{n}{1-\lambda_2}}, \mbb{\left(\frac{\norm{\nabla F(\x^*)}^2}{(1-\lambda_2)^3}\right)^{\frac{1}{3}}},\mbb{\left(\frac{\norm{\x_0-\x^*}^2}{(1-\lambda_2)^6}\right)^{\frac{1}{5}}},\right.\\
			\left.\mbb{\frac{\max\left\{\log \left(\norm{\x_0-\x^*}^2 + \norm{\nabla F(\x^*)}^2\right),-\log(1-\lambda_2)\right\}}{1-\lambda_2}}\right\},
		\end{multline*}
		which finishes the proof.
	\end{proof}
	
	Under mild additional conditions, we can obtain a cleaner expression for the transient time of EDAS in the following corollary.
	\begin{corollary}\label{cor:transient_time}
		Under \cref{alg:extra} with \cref{ngrad,edfi,edw,}, suppose $\theta>5$ and $m$ is chosen according to \eqref{eq:K1}. Assume in addition that $\norm{\x_0-\x^*}^2=\mbb{n^{3}}$, $\norm{\nabla F(\x^*)}^2=\mbb{n^{3}}$ and $\frac{1}{1-\lambda_2}=\mbb{n^q}$ for some $q>0$. Then it takes 
		$$
		K_T = \mbb{\frac{n}{1-\lambda_2}}
		$$
		steps for \cref{alg:extra} to reach the asymptotic, network independent convergence rate. 
	\end{corollary}
	\begin{remark}
		Assuming that $\norm{\x_0-\x^*}^2=\mbb{n^{3}}$ and $\norm{\nabla F(\x^*)}^2=\mbb{n^{3}}$ is mild. 
		This condition can be satisfied for many problem settings including the ones we consider in Section \ref{sec:experiments}.
		In addition, $\frac{1}{1-\lambda_2}=\mbb{n^q}$ is generally satisfied since $\frac{1}{1-\lambda_2}=\mbb{n^2}$ holds under common choices of the mixing weights for undirected graphs \cite{nedic2018network}.
	\end{remark}
	
	\section{Numerical Examples}
	\label{sec:experiments}
	
	In this section, we present two numerical examples to verify and complement our theoretical results. To begin with, we define the transient time for the algorithms of interest as follows for practical consideration: 
	\begin{equation}\label{eq:exp_time}
		\inf_k\left\{k>100: \frac{1}{n}\sum_{i=1}^n\E\left[\norm{x_{i,k} - x^*}^2\right] \leq 2\E\left[\norm{x_k-x^*}^2\right] \right\},
	\end{equation}
	where $x_k\in\R^p$ stands for the $k$-th iterate for the centralized SGD algorithm. 
	
	In the first experiment, we construct a ``hard" optimization problem for which the experimental transient time of EDAS is of the order of $\mbb{\frac{n}{1-\lambda_2}}$, which agrees with the upper bound given in \cref{thm:transient_time}. This verifies the sharpness of the obtained theoretical results.
	Then we consider the problem of logistic regression for classifying handwritten digits from the MNIST dataset \cite{MNIST}.
	In this case, EDAS achieves a transient time close to $\mbb{\frac{n}{(1-\lambda_2)^{0.5}}}$, which is better than the upper bound in \cref{cor:transient_time}. Hence in practice, the performance of EDAS depends on the specific problem instances and can be better than the worse-case scenario.
	
\begin{figure}[htbp]
		\centering
		\subfloat[Ring network topology with $8$ nodes]{\label{fig:circle_graph}\includegraphics[width=6.5cm]{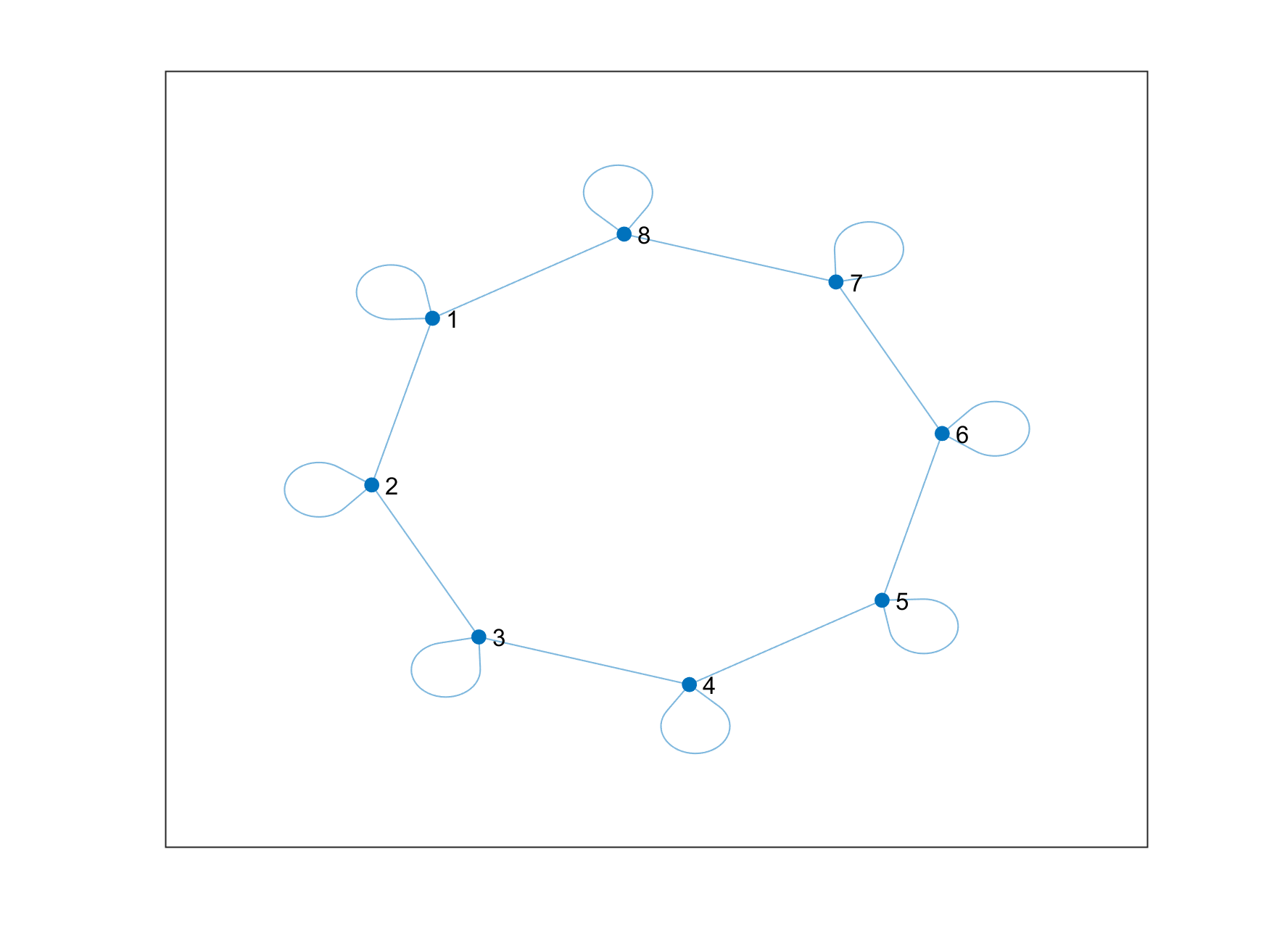}}
		\subfloat[Grid network topology with $16$ nodes]{\label{fig:grid_graph}\includegraphics[width = 6.5cm]{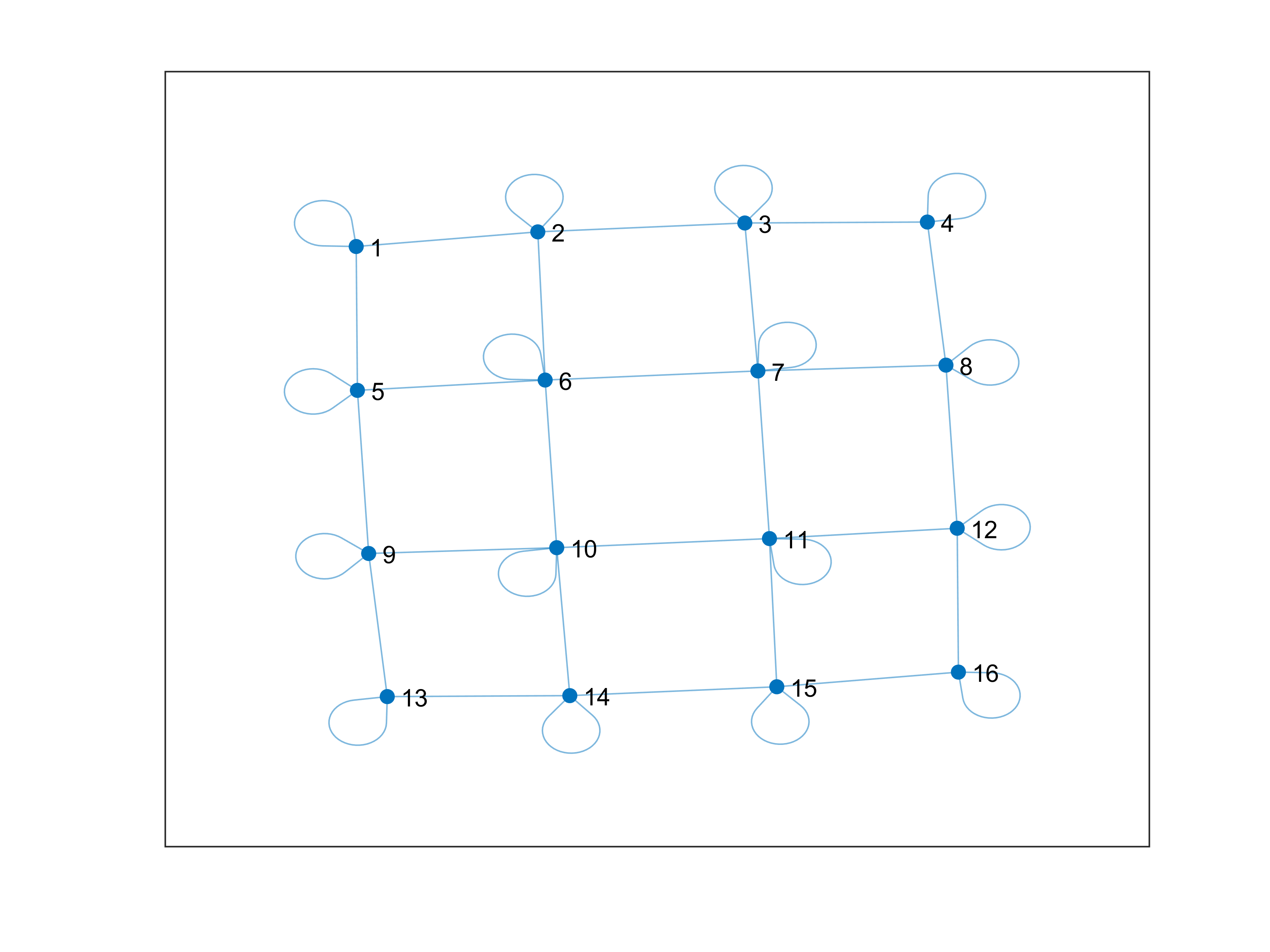}}
		\caption{Two network topologies.}
		\label{fig:graph}
	\end{figure}

	For both problems, we consider two different types of network topologies, i.e., ring network shown in \cref{fig:circle_graph} and grid network shown in \cref{fig:grid_graph}. The mixing matrices compliant with the networks are constructed under the Lazy Metropolis rule \cite{nedic2018network}. In addition to presenting the transient times for EDAS, we also compare its performance with other algorithms enjoying the asymptotic network independent property, i.e., distributed stochastic gradient descent (DSGD) \cite{pu2019sharp,yuan2018exact} and distributed stochastic gradient tracking method (DSGT) \cite{pu2020distributed}.
	
	\subsection{A Near-Sharp Quadratic Problem}\label{sec:sharp}
	Consider the following problem:
	\begin{equation}\label{eq:mean}
		\underset{x\in\R^p}{\min} f(x) = \frac{1}{n}\sum_{i=1}^n\frac{\sqrt{n}}{2}\norm{x-d_i}^2,
	\end{equation}
	where $d_i = a_i\sqrt{i}\1$ and $a_i$ is the $i$-th smallest eigenvalue of $W$. The unique solution to problem \eqref{eq:mean} is given by $x^* = \frac{1}{n}\sum_{i=1}^n d_i$. In the experiments, we set the initial solution to $\mathbf{0}$ for each node in the network; for centralized SGD we use the same initialization. 
	Let $p = 1$, and each node $i$ is able to obtain noisy gradients in the form of $g_i(x_i) = \sqrt{n}(x_i - d_i) + \omega$, where $\omega\sim\mathcal{N}(0, 0.1^2)$ denotes the noise. Stepsizes are chosen as $\alpha_k = 20 / (k + 200),\forall k$.  
	
\begin{figure}[htbp]
		\centering
		\subfloat[Error terms for centralized and decentralized methods.]{\label{fig:errors}
			\includegraphics[width = 6cm]{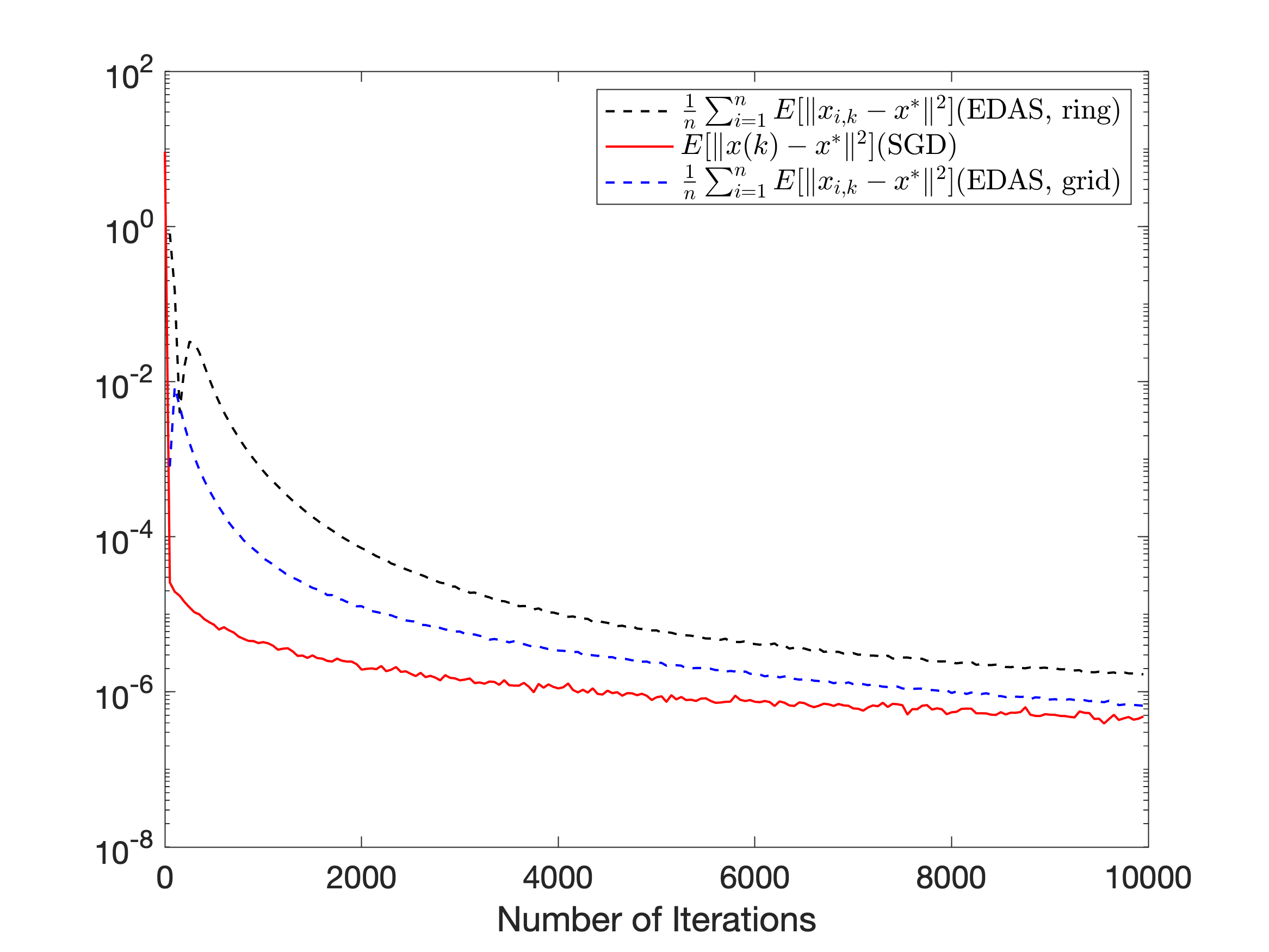}}
		\caption{An illustration of asymptotic network independence.}
		\label{fig:illu}
	\end{figure}

	Fist, we illustrate the property of ``asymptotic network independence'' by comparing the performance of EDAS and SGD  in \cref{fig:errors}. We can see that the error of EDAS gets close to that of centralized SGD as the number of iteration grows. 
	Eventually, it is as if the network were not there and the convergence rate for EDAS is exactly the same as SGD. Hence the phenomenon is called ``asymptotic network independence''. 
	
\begin{figure}[htbp]
		\centering
		\subfloat[Transient time for the ring network topology.]{\label{fig:circle_sharp}\includegraphics[width=6.5cm]{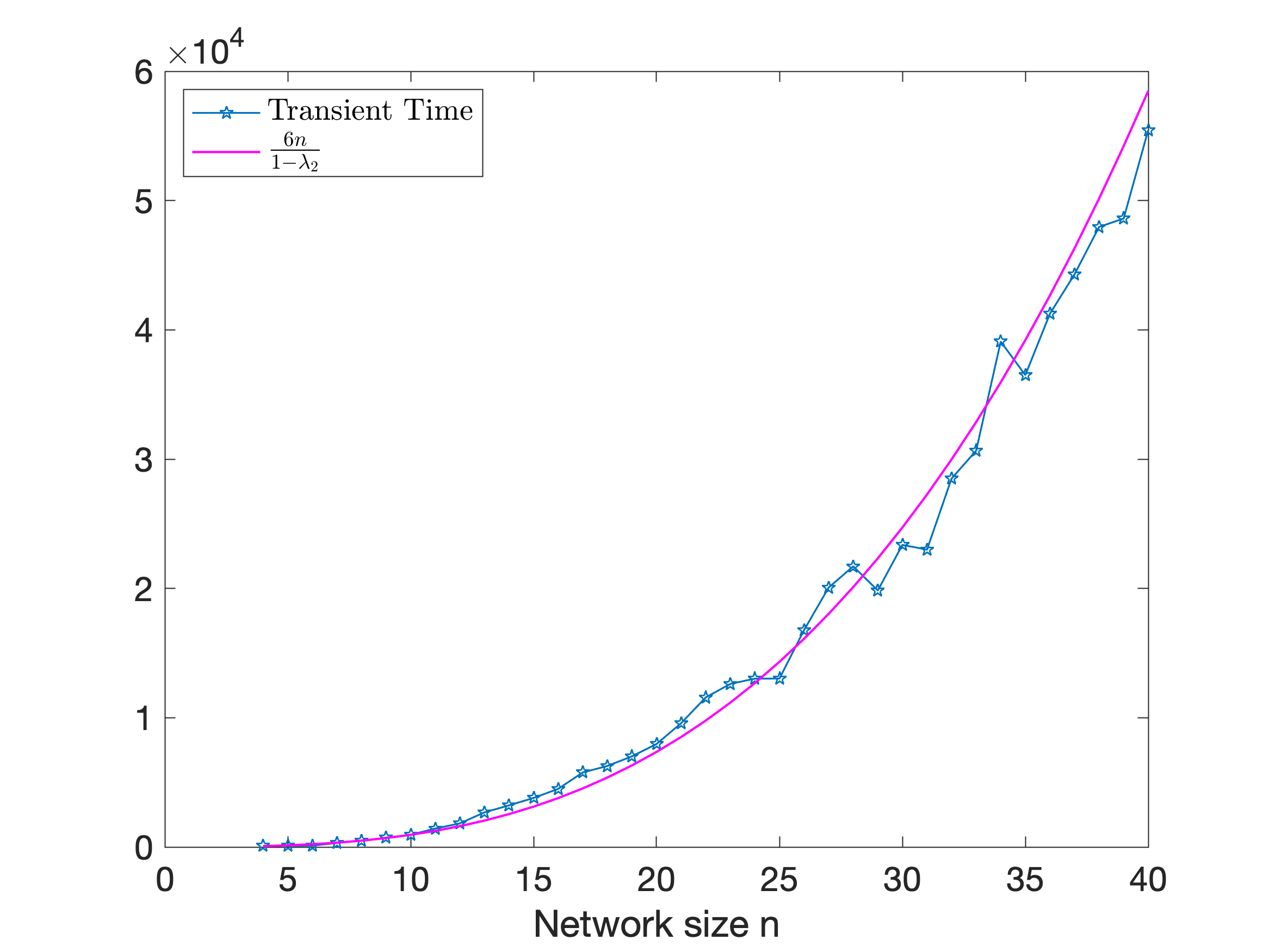}}
		\subfloat[Transient time for the grid network topology.]{\label{fig:grid_sharp}\includegraphics[width = 6.5cm]{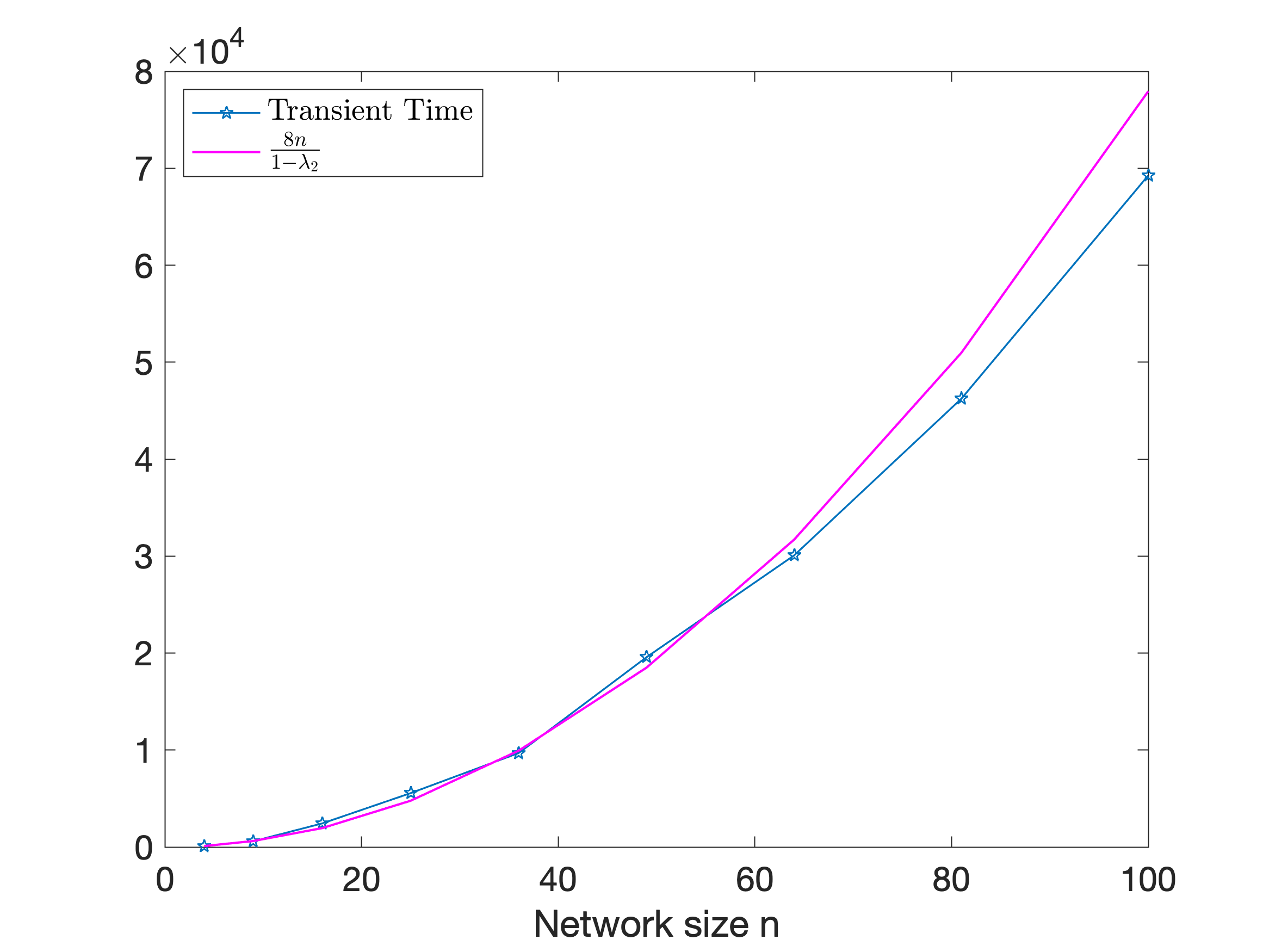}}
		\caption{Comparison of the transient time for \cref{alg:extra} and multiples of $\frac{n}{1-\lambda_2}$. The expected errors are approximated by averaging $500$ simulation results.}
		\label{fig:sharp}
	\end{figure}

	
	In \cref{fig:circle_sharp}, we plot the observed transient times of EDAS and $\frac{6n}{1-\lambda_2}$ against the network size $n$ for the ring network topology.\footnote{For problem \eqref{eq:mean}, the upper bound on the transient time is 
		$K_T = \mbb{\frac{n}{1-\lambda_2}}$ from \cref{cor:transient_time} since we have $\norm{\x_0 - \x^*}^2 = \norm{\nabla F(\x^*)}^2 = \mbb{n^3}$.} We then consider the grid network topology in \cref{fig:grid_sharp}. It can be seen that the two curves are close in both cases, which suggests that the theoretical upper bound $K_T=\mbb{\frac{n}{1-\lambda_2}}$ given in \cref{cor:transient_time} is sharp.
	
\begin{figure}[htbp]
		\centering
		\subfloat[Error terms over a ring network with 20 nodes.]{
			\label{fig:circle20_comp_error}
			\includegraphics[width = 0.3\textwidth]{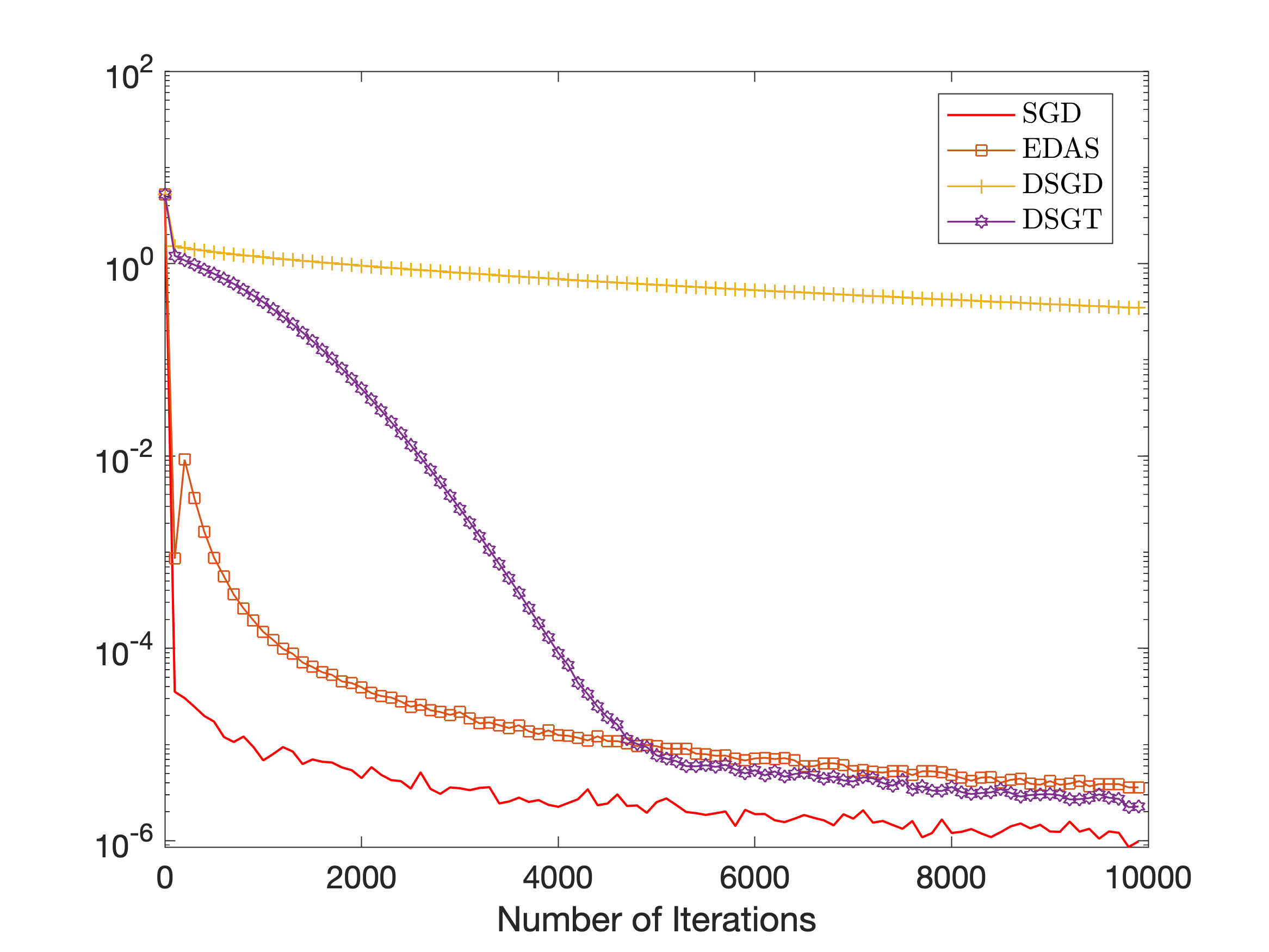}
		}
		\subfloat[Error terms over a ring network with 40 nodes.]{
			\label{fig:circle40_comp_error}
			\includegraphics[width = 0.3\textwidth]{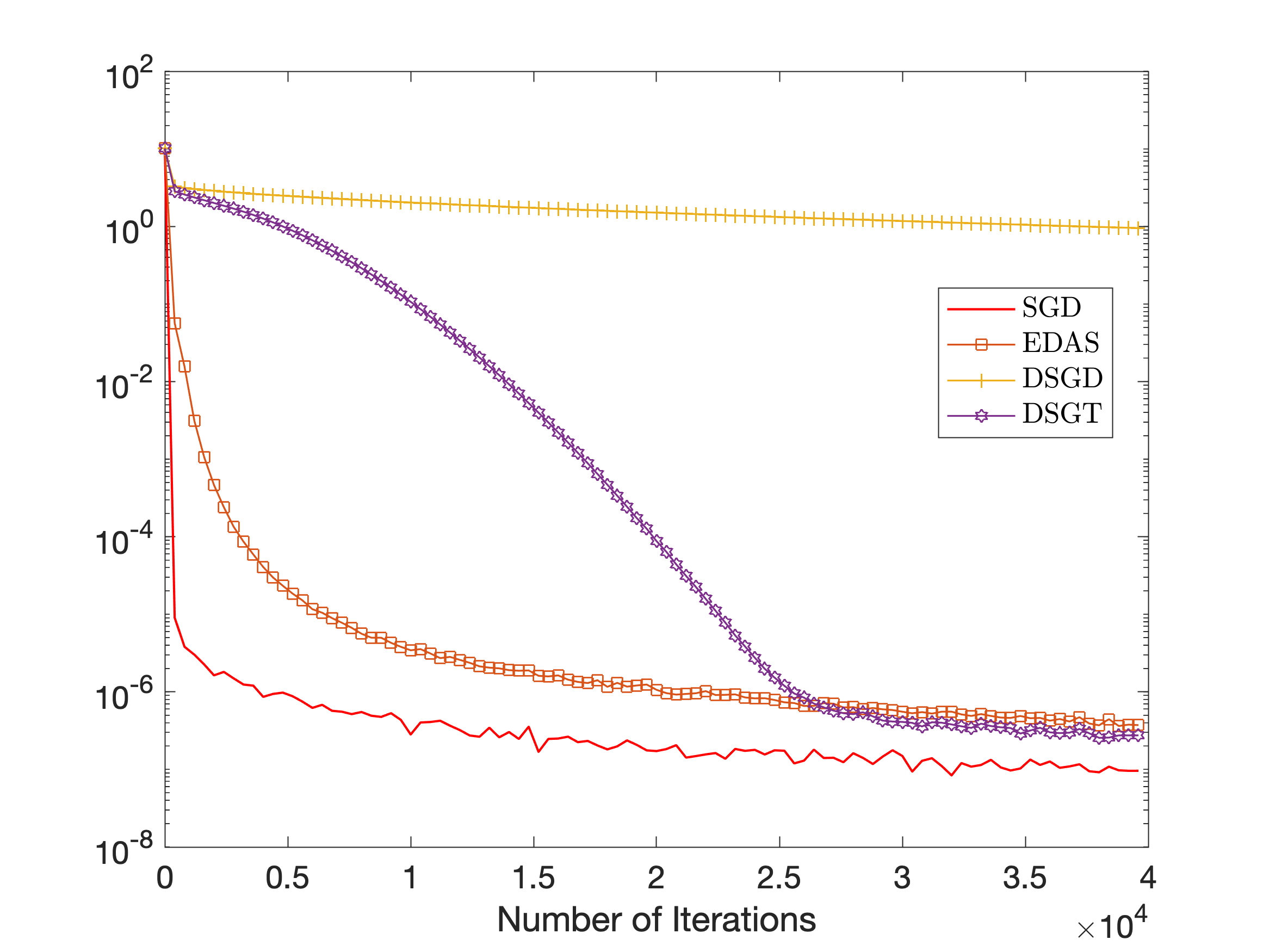}
		}
		\subfloat[Error terms over a ring network with 60 nodes.]{
			\label{fig:circle60_comp_error}
			\includegraphics[width = 0.3\textwidth]{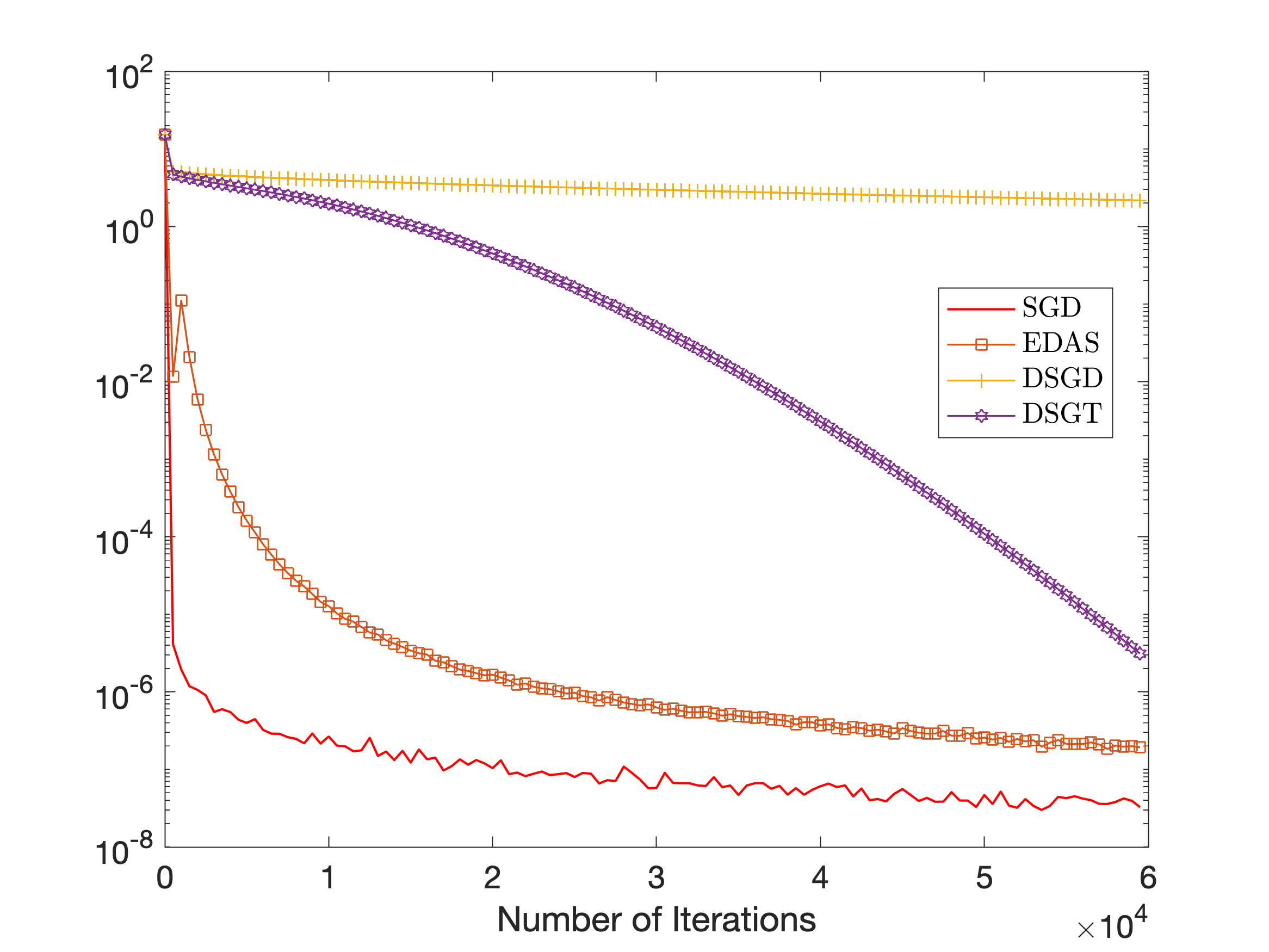}
		}
		\\
		\subfloat[Error terms over a grid network with 49 nodes.]{
			\label{fig:grid49_comp_error}
			\includegraphics[width = 0.3\textwidth]{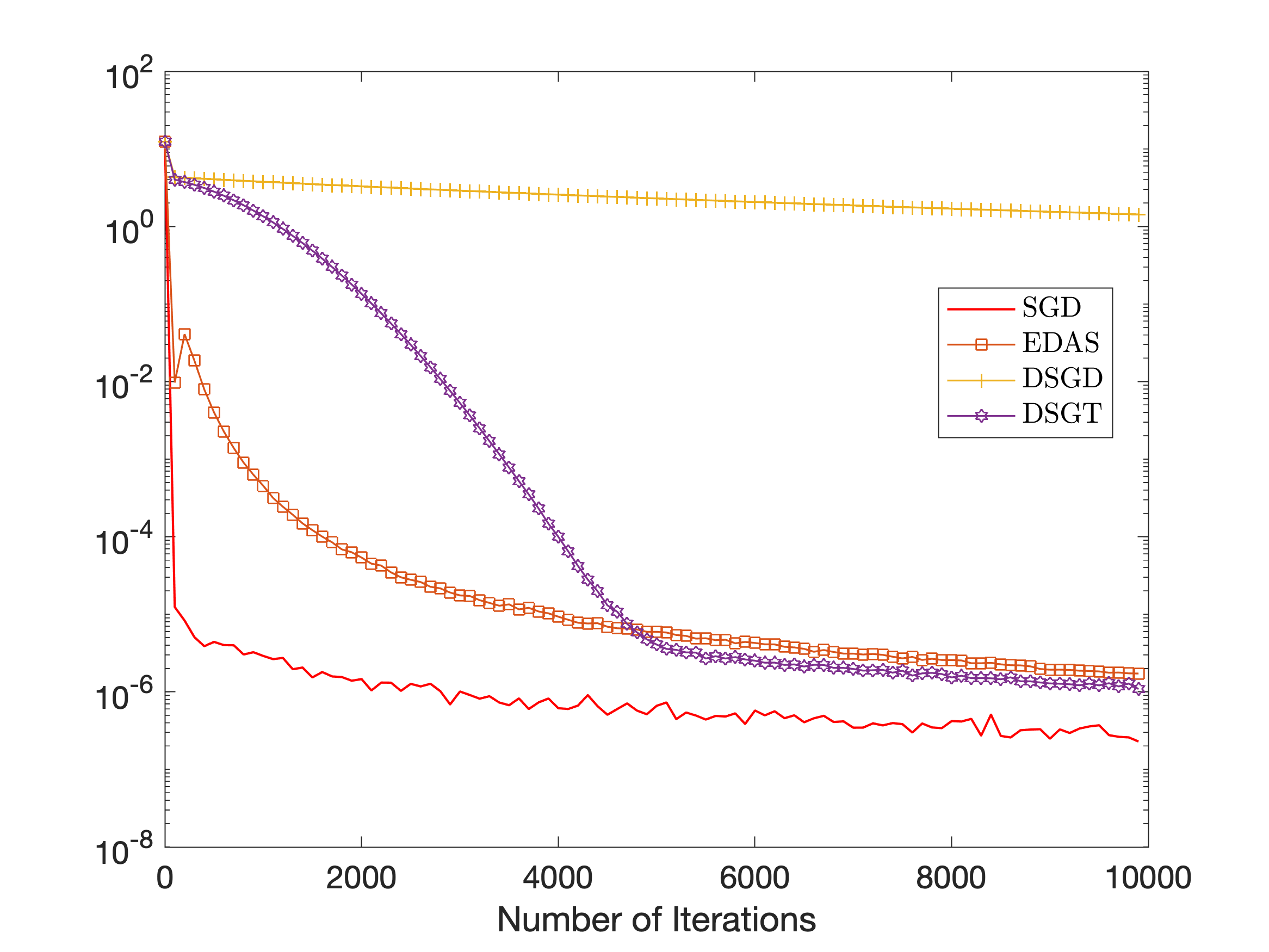}
		}
		\subfloat[Error terms over a grid network with 81 nodes.]{
			\label{fig:grid81_comp_error}
			\includegraphics[width = 0.3\textwidth]{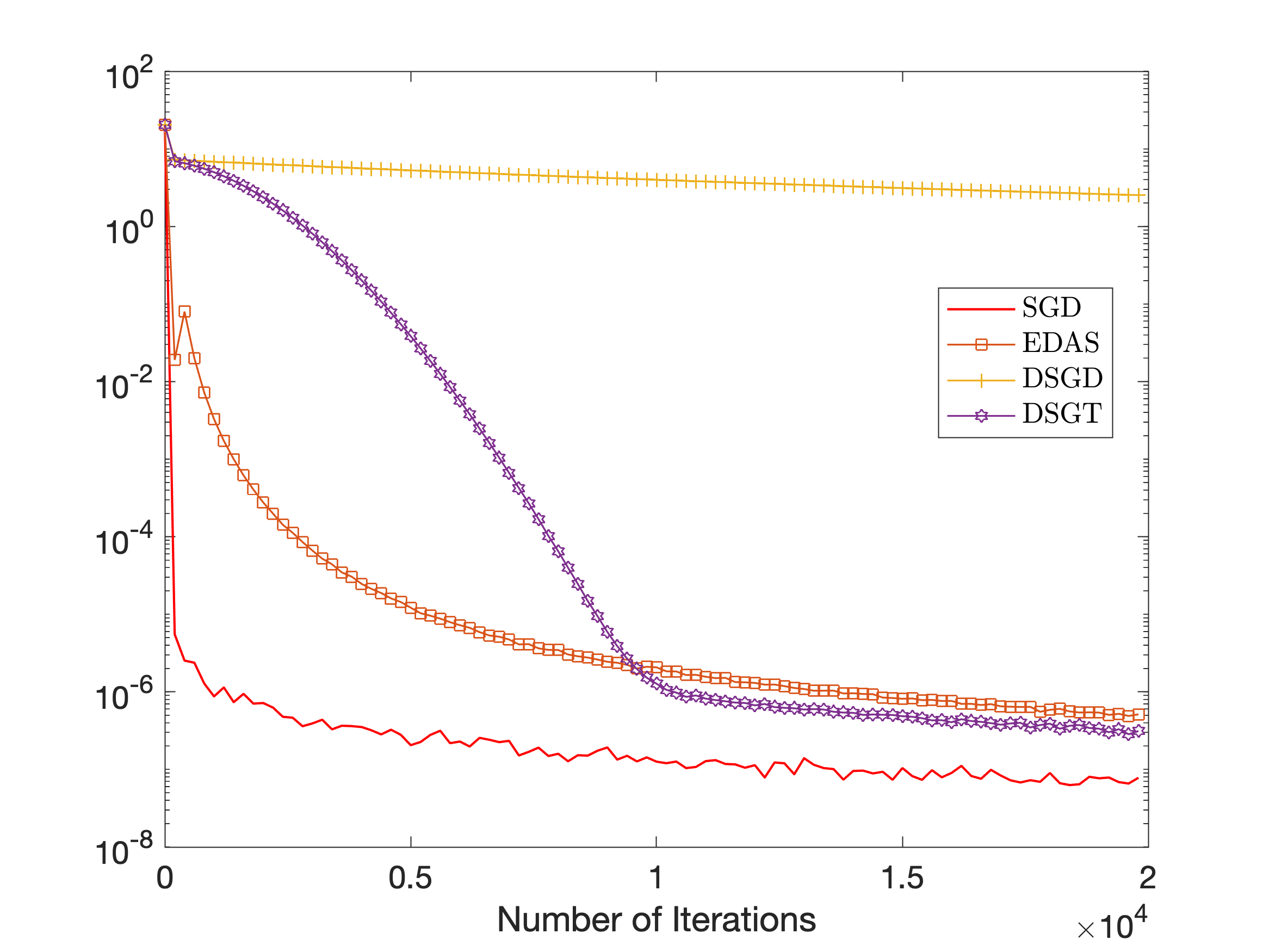}
		}
		\subfloat[Error terms over a grid network with 121 nodes.]{
			\label{fig:grid121_comp_error}
			\includegraphics[width = 0.3\textwidth]{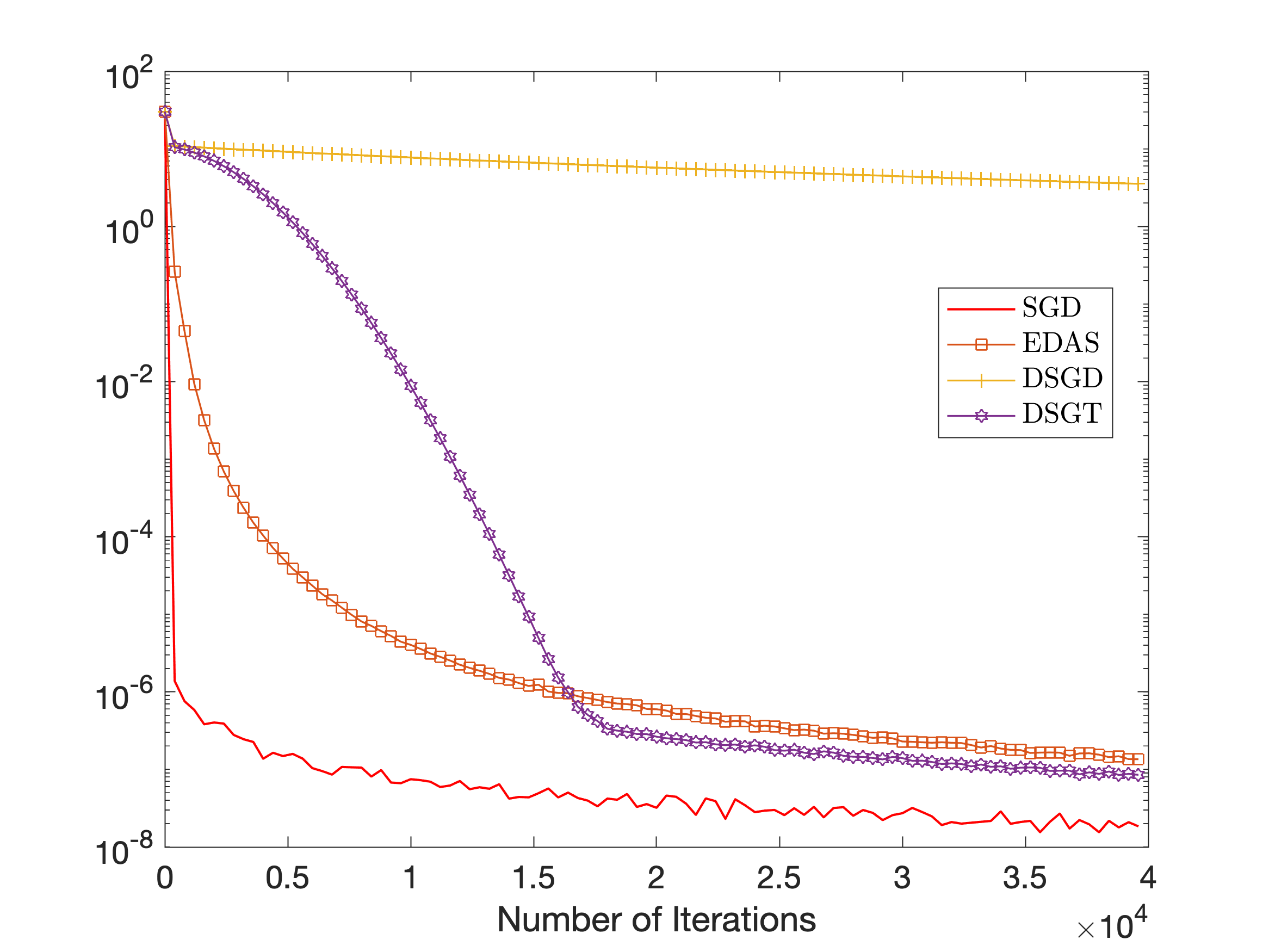}
		}
		\caption{Comparison among DSGT, DSGD, and EDAS over ring and grid network topologies for problem \eqref{eq:mean}. The curves stand for $\E(\norm{x(k)-x^*}^2)$ for SGD and $\frac{1}{n}\sum_{i=1}^n\E(\norm{x_{i,k}-x^*}^2)$ for decentralized methods. The expected errors are approximated by averaging over 100 simulation results.}
		\label{fig:comp}
	\end{figure}

	In \cref{fig:comp}, we further compare the early stage performance of DSGT, DSGD, and EDAS with the same initialization. It can be seen that the error term of EDAS decreases faster than that of the other two decentralized methods and quickly gets close to the error of SGD. This is particularly true when the network is not well connected, e.g., when the network has a ring topology with a large number of nodes (see \cref{fig:circle60_comp_error}).
	
	\subsection{Logistic Regression}
	We consider classifying handwritten digits $1$ and $2$ in the MNIST dataset using logistic regression. There are $12700$ data points $\mc{S} = \{(u, v)\}$ with $u\in\R^{785}$ denoting the image input and $v\in\{-1, 1\}$ being the label.\footnote{Label $-1$ represents handwritten digit $2$ and label $1$ represents digit $1$.} Each agent possesses a distinct local dataset $\mc{S}_i$ randomly selected from $\mc{S}$. The classifier can then be obtained by solving the following optimization problem using all the agents' local datasets $\mc{S}_i, i=1,2,...,n$:
	\begin{subequations}
		\label{eq:logistic}
		\begin{align}
			\min_{x\in\R^{785}} f(x) &= \frac{1}{n}\sum_{i=1}^n f_i(x),\\
			f_i(x) &:= \frac{1}{|\mc{S}_i|} \sum_{j\in\mc{S}_i} \log\left[1 + \exp(-x^{\T}u_jv_j)\right] + \frac{\rho}{2}\norm{x}^2,
		\end{align}
	\end{subequations}
	where $\rho >0$ is a regularization parameter (hyperparameter). To solve the problem, each agent can obtain an unbiased estimate of $\nabla f_i(x)$ using a minibatch of randomly selected data points from its local dataset $\mc{S}_i$. In the experiment, we let $\rho = 1$, $\alpha_k = 6 / (k + 20), \forall k$ for all the considered algorithms with the same initial solutions: $x_{i, 0} = \0, \forall i$ (distributed methods) and $x_0 = \0$ (centralized SGD). Each local dataset has $|\mc{S}_i| = 100$ data points and the mini-batch size is set to be $1$. 
	
\begin{figure}[htbp]
		\centering
		\subfloat[Transient time for the ring network topology.]{\label{fig:circle_mnist}\includegraphics[width=6.5cm]{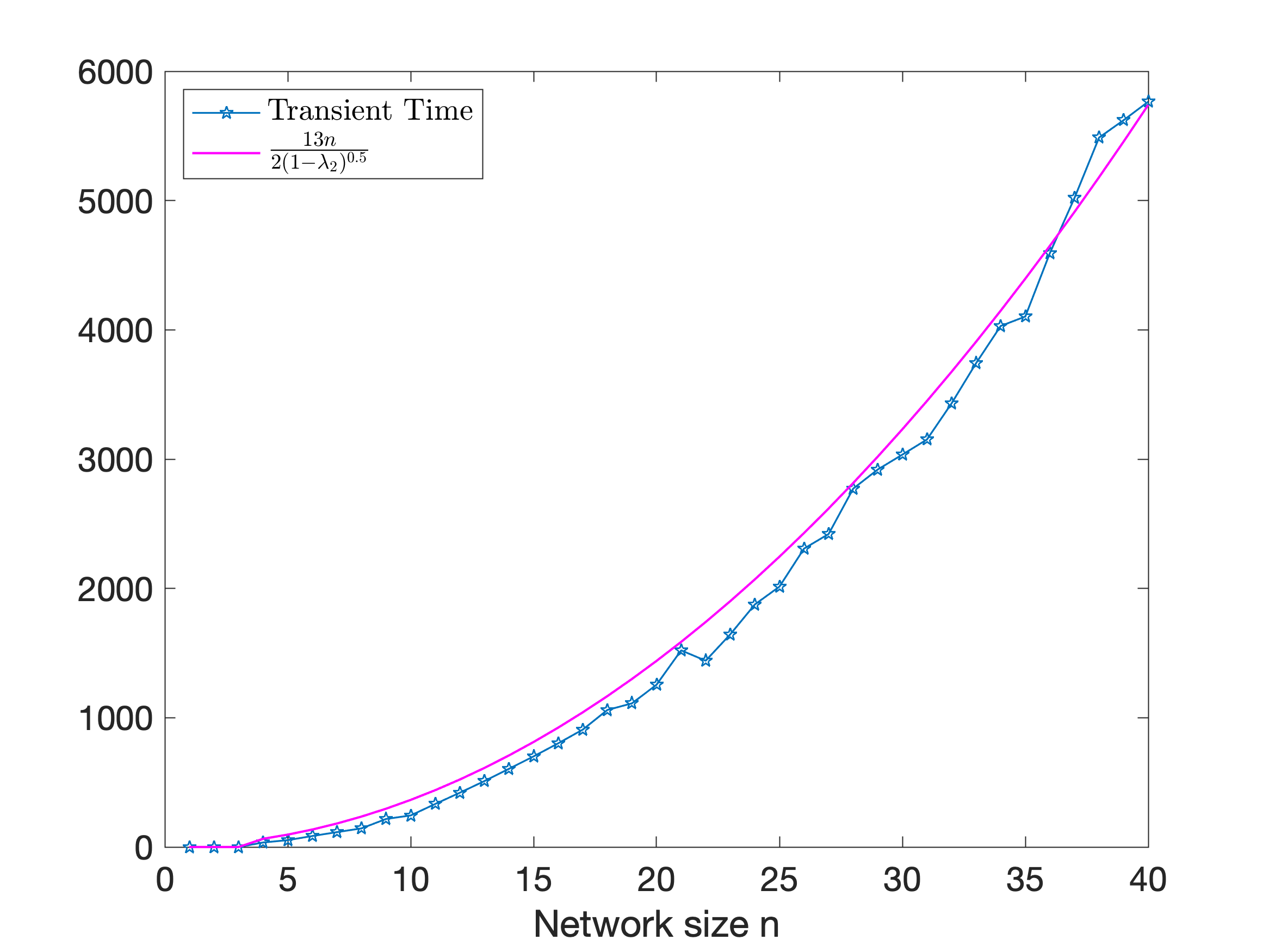}}
		\subfloat[Transient time for the grid network topology.]{\label{fig:grid_mnist}\includegraphics[width = 6.5cm]{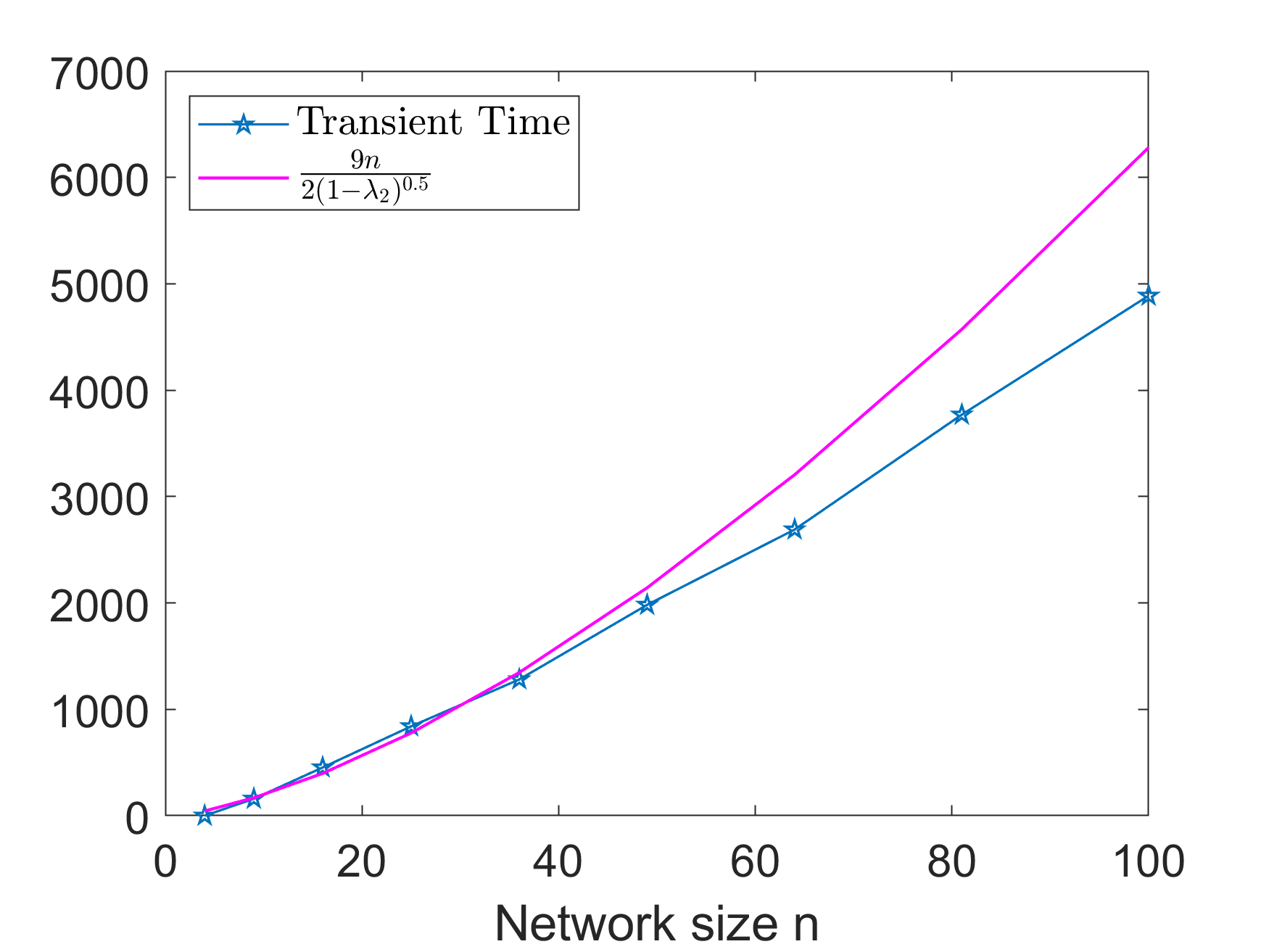}}
		\caption{Comparison between the transient times for EDAS and multiples of $\frac{n}{(1-\lambda_2)^{0.5}}$. The expected errors are approximated by averaging over $200$ simulation results.}
		\label{fig:mnist}
	\end{figure}

	In \cref{fig:mnist}, we can see that the observed transient times of EDAS are close to the multiple of $\frac{n}{(1-\lambda_2)^{0.5}}$. Hence in practice EDAS may achieve a better transient time than the derived worst-case upper bound.   
	
\begin{figure}[htbp]
		\centering
		\subfloat[Error terms over a ring network with 20 nodes.]{
			\label{fig:circle20_comp_error_mnist}
			\includegraphics[width = 0.3\textwidth]{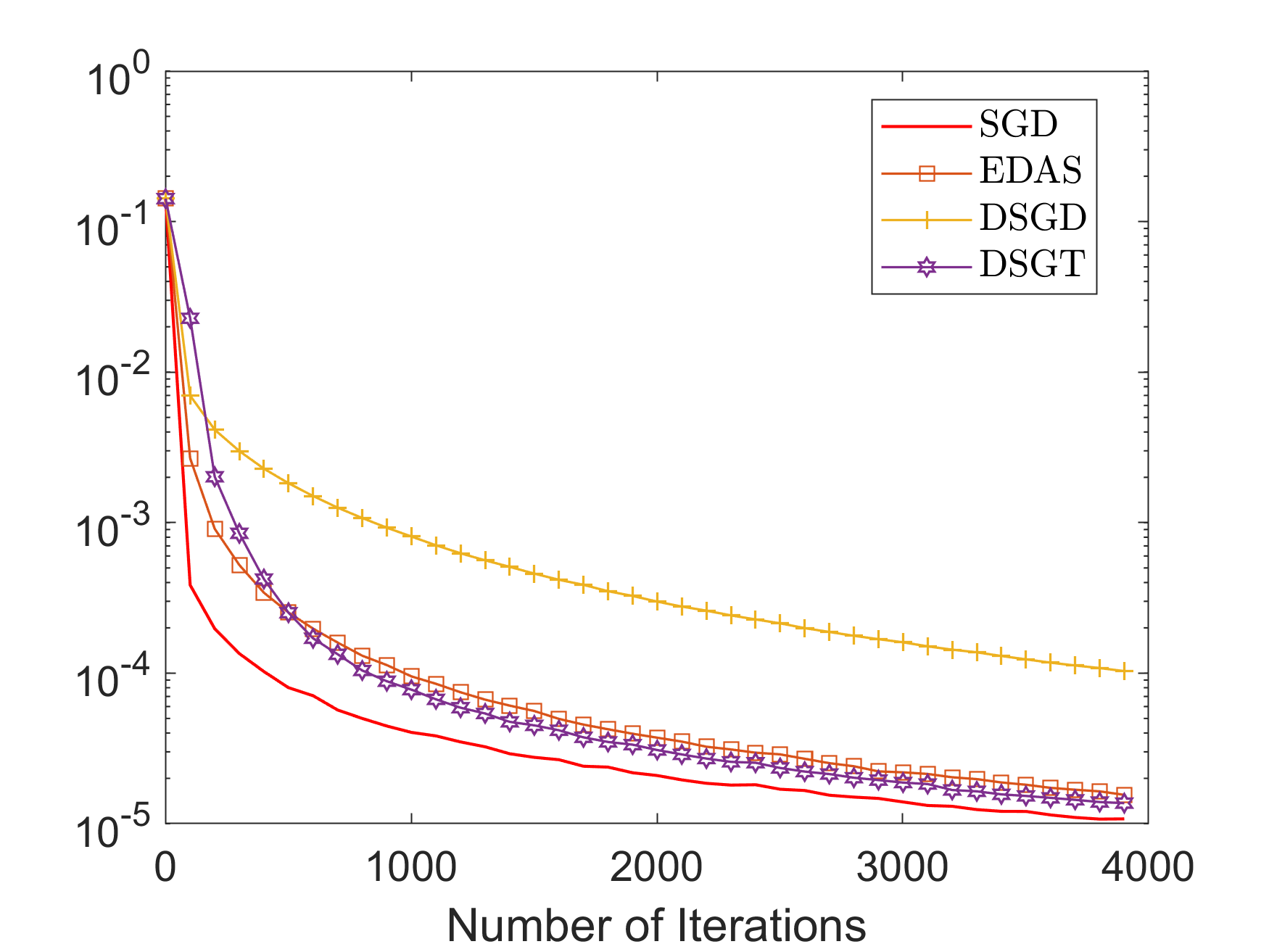}
		}
		\subfloat[Error terms over a ring network with 40 nodes.]{
			\label{fig:circle40_comp_error_mnist}
			\includegraphics[width = 0.3\textwidth]{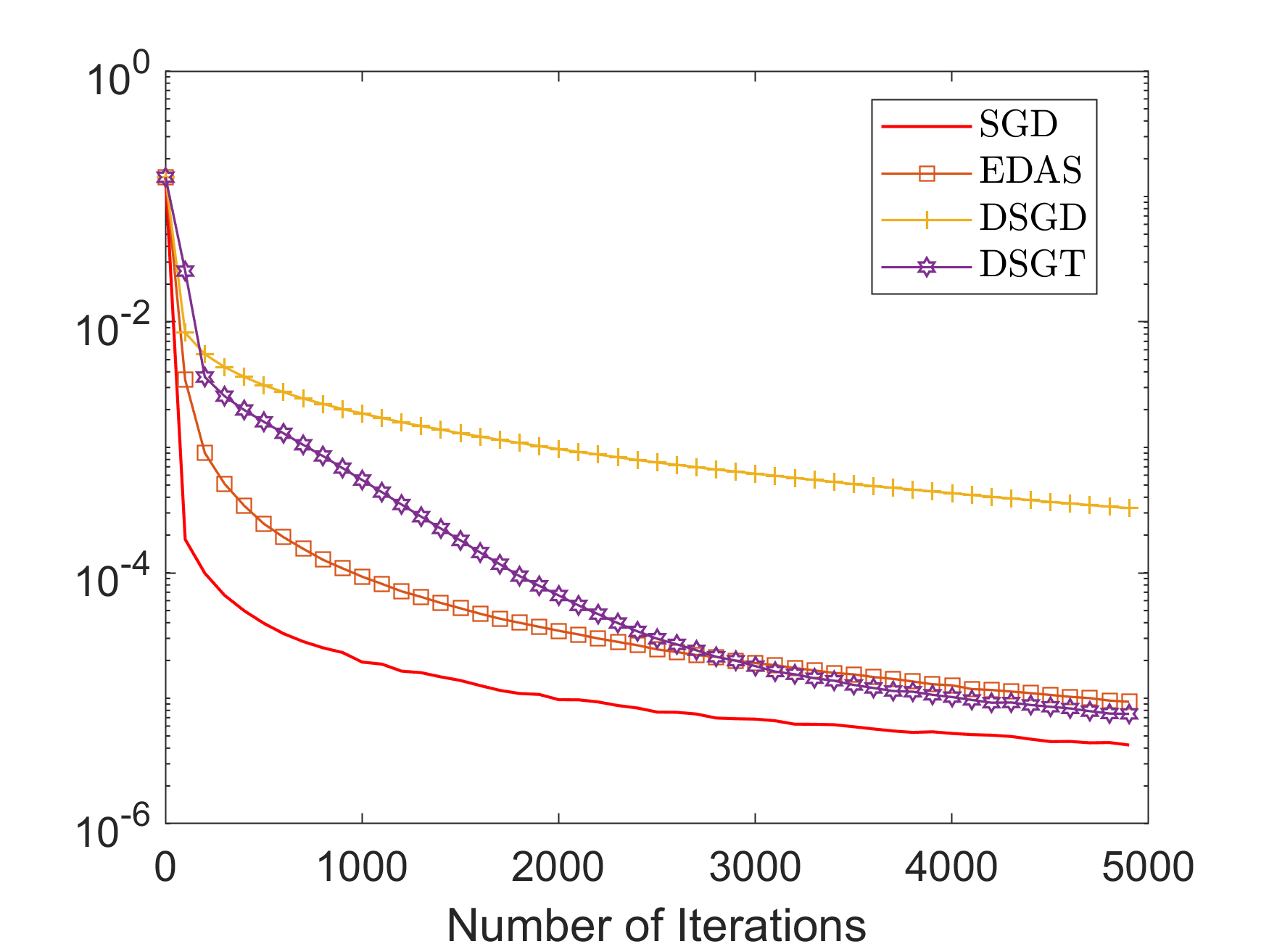}
		}
		\subfloat[Error terms over a ring network with 60 nodes.]{
			\label{fig:circle60_comp_error_mnist}
			\includegraphics[width = 0.3\textwidth]{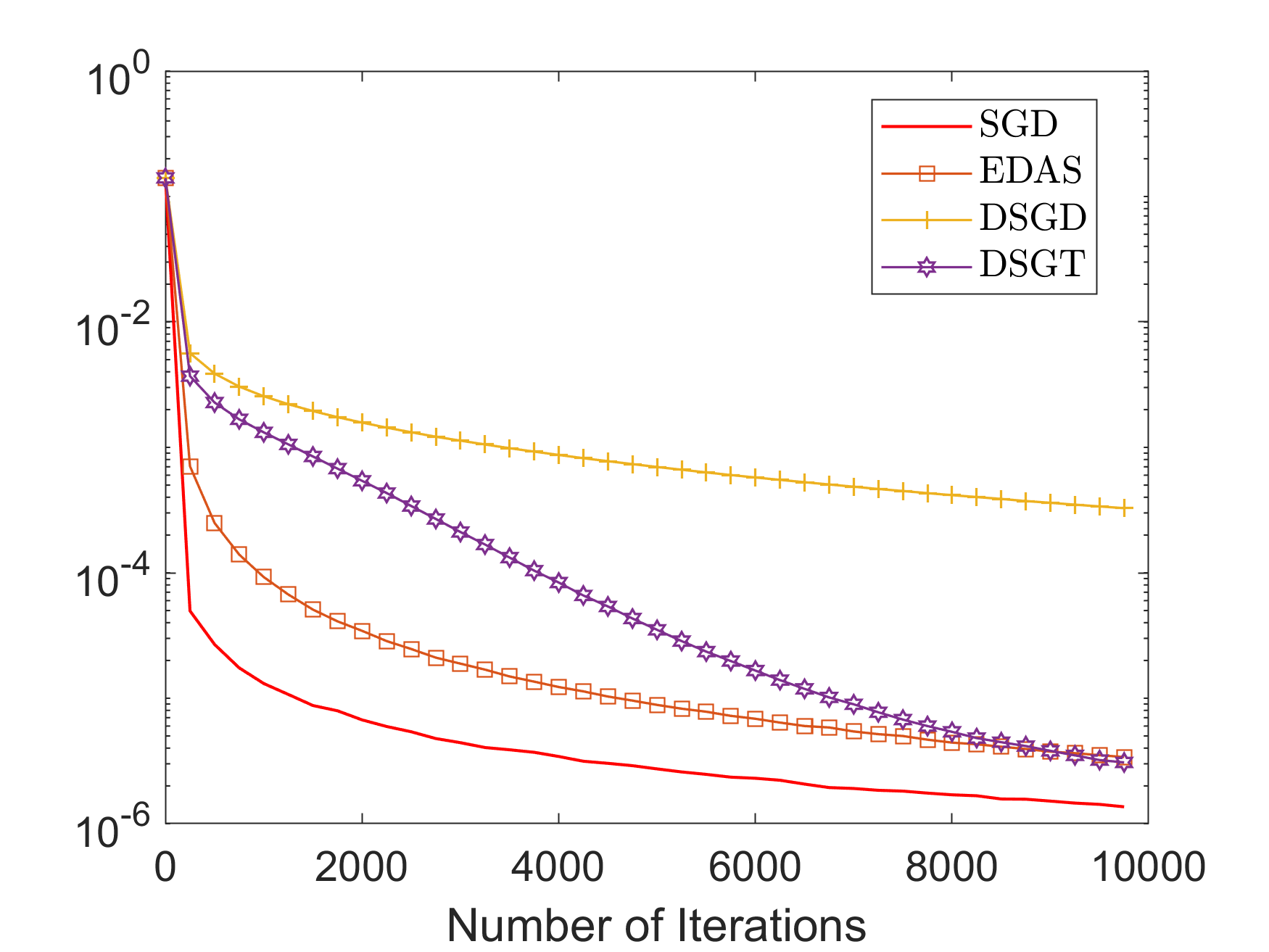}
		}
		\\
		\subfloat[Error terms over a grid network with 49 nodes.]{
			\label{fig:grid49_comp_error_mnist}
			\includegraphics[width = 0.3\textwidth]{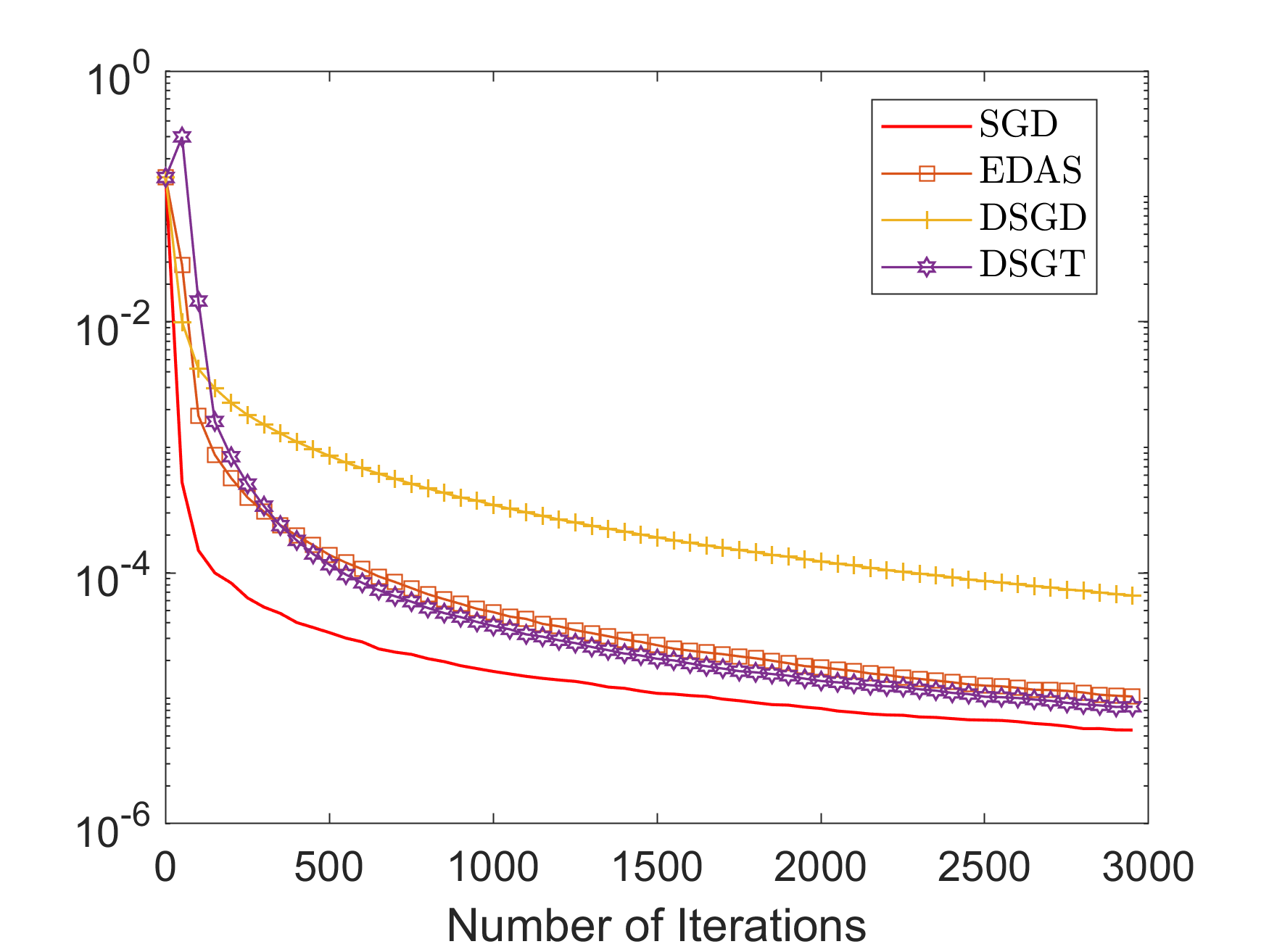}
		}
		\subfloat[Error terms over a grid network with 81 nodes.]{
			\label{fig:grid81_comp_error_mnist}
			\includegraphics[width = 0.3\textwidth]{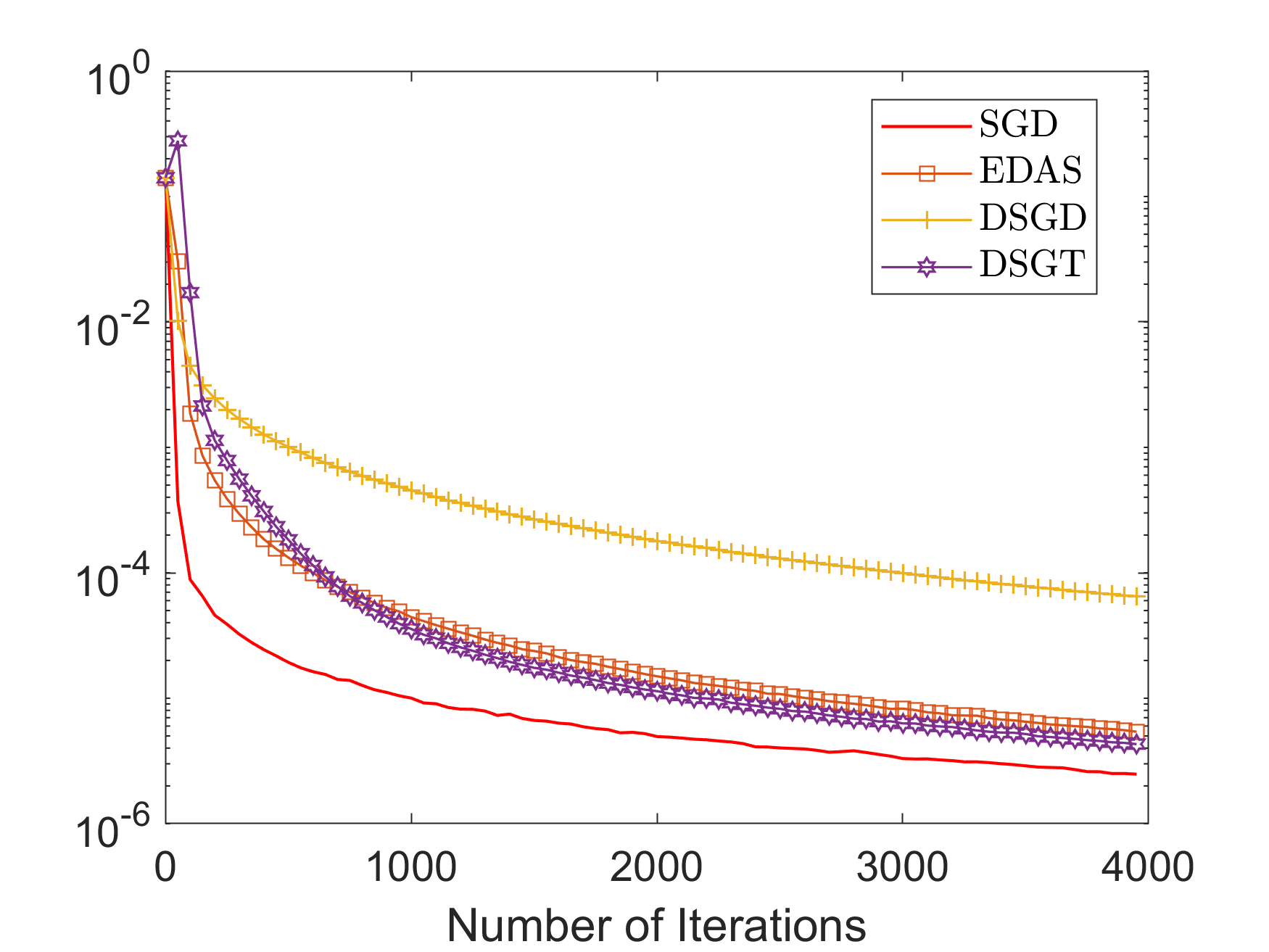}
		}
		\subfloat[Error terms over a grid network with 121 nodes.]{
			\label{fig:grid121_comp_error_mnist}
			\includegraphics[width = 0.3\textwidth]{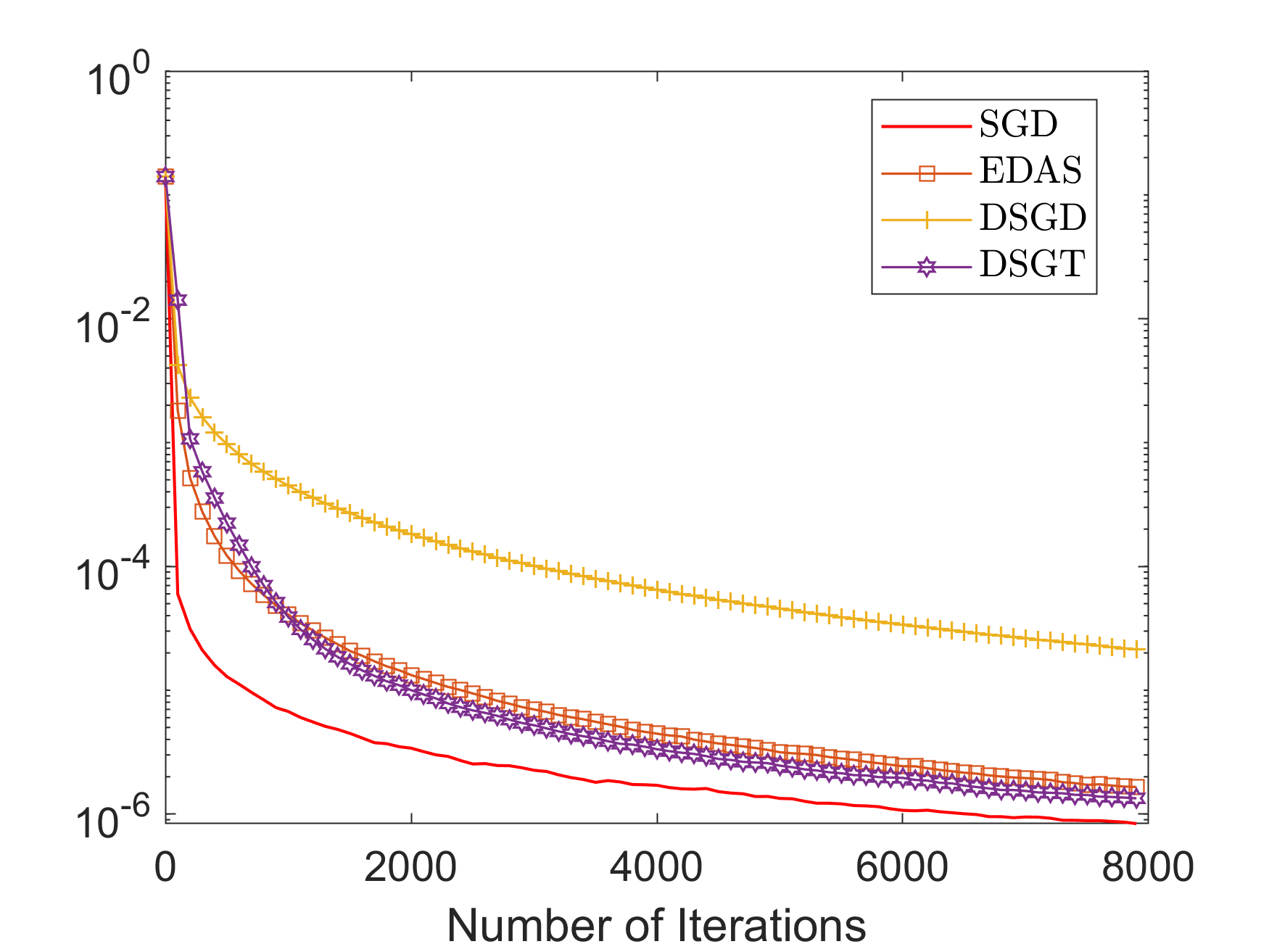}
		}
		\caption{Comparison among DSGT, DSGD, and EDAS over ring and grid network topologies for problem \eqref{eq:logistic}. The y-axis stands for $\E(\norm{x(k)-x^*}^2)$ for SGD and $\frac{1}{n}\sum_{i=1}^n\E(\norm{x_{i,k}-x^*}^2)$ for decentralized methods. The expected errors are approximated by averaging 100 simulation results.}
		\label{fig:comp_mnist}
	\end{figure}

	We also compare the performance of EDAS with DSGD and DSGT for the logistic regression problem in \cref{fig:comp_mnist}. On one hand, similar to problem \eqref{eq:mean}, EDAS achieves significantly shorter transient times for ring networks when the network size is moderately large; see \cref{fig:circle20_comp_error_mnist,fig:circle40_comp_error_mnist,fig:circle60_comp_error_mnist}. On the other hand, for grid networks, EDAS and DSGT both achieve better transient times than DSGD; see \cref{fig:grid49_comp_error_mnist,fig:grid81_comp_error_mnist,fig:grid121_comp_error_mnist}.
	
	\section{Conclusions}
	\label{sec:conclusions}
	
	In this paper, we consider a distributed stochastic gradient algorithm (EDAS) and perform a non-asymptotic analysis. In addition to showing that EDAS achieves asymptotic network independence, we identify its non-asymptotic convergence rate as a function of characteristic of the objective functions and the network. We show that EDAS achieves the fastest transient time so far to the best of our knowledge. Future work will explore the lower bound of transient times for a general class of distributed stochastic gradient methods.
	
	\appendix
	\section[Proofs]{Proofs}
	\subsection[Proof of Lemma 2.5]{Proof of \cref{lem:B_decomposition}}
	\label{sec:prf_KL_KR}
	\begin{proof}
		Since $W$ is a stochastic symmetric matrix, the matrices $W$ and $V$ have the following spectral decomposition:
		\begin{equation*}
			W = Q\Lambda Q^{\T}, \quad V = Q\sqrt{I_n-\Lambda}Q^{\T}
		\end{equation*}
		where $Q$ is an orthogonal matrix, and
		$$\Lambda := \mathrm{diag}\left(1, \lambda_2,\ldots, \lambda_n\right).$$ 
		
		Therefore, 
		\begin{equation}
			\label{B_dec1}
			B = \begin{aligned}
				\left(
				\begin{array}{cc}
					W & -V\\
					VW & W
				\end{array}
				\right) 
				&= 
				\left(
				\begin{array}{cc}
					Q & \\
					& Q
				\end{array}
				\right)
				\left(
				\begin{array}{cc}
					\Lambda & -\sqrt{I_n-\Lambda}\\
					\sqrt{I_n-\Lambda}\Lambda & \Lambda
				\end{array}
				\right)
				\left(
				\begin{array}{cc}
					Q^{\T} & \\
					& Q^{\T}
				\end{array}
				\right).
			\end{aligned}
		\end{equation}
		
		Noticing the particular form of the middle term on the right hand side of \eqref{B_dec1},  there exists a permutation operator $\Pi_1\in\R^{2n\times 2n}$ such that
		\begin{equation*}
			\left(
			\begin{array}{cc}
				\Lambda & -\sqrt{I_n-\Lambda}\\
				\sqrt{I_n-\Lambda}\Lambda & \Lambda
			\end{array}
			\right) = 
			\Pi_1 \mathrm{diag}(E_1, E_2, ..., E_n)\Pi_1^{\T},
		\end{equation*}
		where 
		\begin{equation}
			E_1 := I_2, \;
			E_k := \left(
			\begin{array}{cc}
				\lambda_k & -\sqrt{1-\lambda_k}\\
				\lambda_k\sqrt{1-\lambda_k} & \lambda_k
			\end{array}
			\right), k=2, 3, \ldots, n.
		\end{equation}
		
		Hence
		\begin{equation}
			\label{B_dec2}
			B = \begin{aligned} 
				\left(
				\begin{array}{cc}
					Q & \\
					& Q
				\end{array}
				\right)
				\Pi_1 \mathrm{diag}(I_2, E_2, ..., E_n)\Pi_1^{\T}
				\left(
				\begin{array}{cc}
					Q^{\T} & \\
					& Q^{\T}
				\end{array}
				\right).
			\end{aligned}
		\end{equation}
		
		We then consider the eigendecomposition of $E_k$ for $k=2, 3, \ldots, n$. It can be verified that
		\begin{equation}\label{D_1}
			E_k = Z_k \left(
			\begin{array}{cc}
				\lambda_k + i\sqrt{\lambda_k-\lambda_k^2} &0 \\
				0 & \lambda_k - i\sqrt{\lambda_k-\lambda_k^2}
			\end{array}
			\right) Z_k^{-1},
		\end{equation}
		where $i^2=-1$, and 
		\begin{equation}
			\begin{aligned}
				Z_k &= \left(
				\begin{array}{cc}
					1 & 1\\
					-i\sqrt{\lambda_k} & i\sqrt{\lambda_k}
				\end{array}
				\right),\quad
				Z_k^{-1} &= \left(
				\begin{array}{cc}
					\frac{1}{2} & \frac{i}{2\sqrt{\lambda_k}}\\
					\frac{1}{2} & -\frac{i}{2\sqrt{\lambda_k}}
				\end{array}
				\right).
			\end{aligned}
		\end{equation}
		
		Therefore, from \eqref{B_dec2}, we have
		\begin{equation}
			\label{B_dec3}
			B = \begin{aligned} 
				\left(
				\begin{array}{cc}
					Q & \\
					& Q
				\end{array}
				\right)
				\Pi_1 Z\left(\mathrm{diag}(I_2,D_1)\right)Z^{-1}\Pi_1^{\T}
				\left(
				\begin{array}{cc}
					Q^{\T} & \\
					& Q^{\T}
				\end{array}
				\right),
			\end{aligned}
		\end{equation}
		where $Z=\mathrm{diag}\left(I_2, Z_1, Z_2, ..., Z_n\right)$ and $D_1$ is a complex diagonal matrix.
		
		Additionally, note that $\left(\mathbf{1}_n^{\T}, \0^{\T}\right)^{\T}$ and $\left(\0^{\T}, \mathbf{1}_n^{\T}\right)^{\T}$ are the two eigenvectors of $B$ w.r.t. eigenvalue $1$. Thus after a permutation $\Pi_2$, we reach a decomposition for $B$ given by $B = \tilde{U}D\tilde{U}^{-1}$, 
		where
		\begin{align*}
			& D:=\mathrm{diag}(I_2,D_1),\\
			& \tilde{U}: = \left(
			\begin{array}{cc}
				Q & \\
				& Q
			\end{array}
			\right)\Pi_1 Z \Pi_2 = \begin{pmatrix}
				\1_n & \0 & U_{R,u}\\
				\0 & \1_n & U_{R,l}
			\end{pmatrix},
		\end{align*}
		for some $U_{R,u},U_{R,l}\in\mathbb{R}^{n\times (n-2)}$. We then further choose $c>0$ as a free parameter, to obtain $B = UDU^{-1}$,  
		where
		\begin{align*}
			& U: = \begin{pmatrix}
				\1_n & \0 & c U_{R,u}\\
				\0 & \1_n & c U_{R,l}
			\end{pmatrix},\;U^{-1}:=\begin{pmatrix}
				\frac{1}{n}\1_n^{\T} & \0\\
				\0 & \frac{1}{n}\1_n^{\T}\\
				\frac{1}{c} U_{L,l} & \frac{1}{c} U_{L,r} 
			\end{pmatrix}.
		\end{align*}
		
		We next bound $\Vert U_R\Vert_2$ and $\Vert U_L\Vert_2$. By definition, we have
		\begin{equation*}
			\begin{aligned}
				\Vert U_R\Vert_2\le \Vert \tilde{U}\Vert_2 &\leq \left\Vert \left(
				\begin{array}{cc}
					Q & \\
					& Q
				\end{array}
				\right)\right\Vert_2 \Vert\Pi_1\Vert_2 \Vert Z\Vert_2 \Vert\Pi_2\Vert_2
				= \Vert Z\Vert_2.
			\end{aligned}
		\end{equation*}
		
		Since
		\begin{equation*}
			\begin{aligned}
				\Vert Z_k\Vert_2&= \sqrt{\lambda_{\max}(Z_kZ_k^*)} 
				= \sqrt{\lambda_{\max}\left(
					\begin{array}{cc}
						2 &0 \\
						0 & 2|\lambda_k|
					\end{array}
					\right)}
				= \sqrt{2},\;\forall k=2,3,\ldots,n,
			\end{aligned}
		\end{equation*}
		we obtain $\Vert U_R\Vert_2\leq \Vert Z\Vert_2 = \sqrt{2}$. 
		
		Similarly, 
		\begin{equation*}
			\Vert U_L\Vert_2 \le \Vert \tilde{U}^{-1}\Vert_2 \leq \Vert Z^{-1}\Vert_2.
		\end{equation*}
		
		Given that
		\begin{equation*}
			Z_k^{-1} (Z_k^{-1})^* = 
			\frac{1}{4}
			\left(
			\begin{array}{cc}
				1+\frac{1}{|\lambda_k|} & 1-\frac{1}{|\lambda_k|}\\
				1-\frac{1}{|\lambda_k|} & 1+\frac{1}{|\lambda_k|}
			\end{array}
			\right),
		\end{equation*}
		we have
		\begin{equation*}
			\Vert Z_k^{-1}\Vert_2 = \sqrt{\lambda_{\max}(Z_k^{-1}(Z_k^{-1})^*)} = \sqrt{\frac{1}{2|\lambda_k|}}.
		\end{equation*}
		
		Hence $\Vert U_L\Vert_2 \leq \Vert Z^{-1}\Vert_2 = \sqrt{\frac{1}{2|\lambda_n|}}\le \sqrt{\frac{1}{2\underline\lambda}}$ by assuming $\lambda_n \geq \underline\lambda >0$.
	\end{proof}
	



	\subsection[Proof of Lemma 4.2]{Proof of \cref{lem:bigO}}\label{sec:prf_bigO}
	\begin{proof}
		From the definitions of $M_k, T_k$ in \eqref{eq:MT}, we have
		\begin{equation*}
			\begin{aligned}
				M_0 &= \|\bar{z}_0\|^2 = \Vert\frac{1}{n}\1_n^{\T}(\x_0-\x^*)\Vert^2= \mbb{\frac{\norm{\x_0-\x^*}^2}{n}},\\
				T_0 &= \Vert\frac{1}{c}(U_{L,l}\tx_0+U_{L,r}\ty_0)\Vert^2 \leq \frac{2\Vert U_L\Vert^2}{c^2}\left(\Vert \x_0-\x^* \Vert^2 + \Vert \y_0-\y^*_0\Vert^2\right).
			\end{aligned}
		\end{equation*}
		
		For $\Vert \y_0-\y^*_0\Vert^2$, notice that $\y_0 = \mathbf{0}$ in \cref{alg:extra}, and $\y^*_0 = -V^{-}\alpha_0W\nabla F(\x^*)$  by \cref{lem:existence}, we have 
		\begin{equation*}
			\norm{\y_0-\y_0^*} = \norm{V^{-}\alpha_0W\nabla F(\x^*)}.
		\end{equation*}
		
		Given that $\norm{V^{-}} = \norm{Q\left(\sqrt{I_n - \bar{\Lambda}}\right)^{-}Q^{\T}} = \sqrt{\frac{1}{1-\lambda_2}}$ and $\alpha_0 \le \frac{1-\lambda_2}{\mu}$ from \eqref{eq:K1}, we derive
		\begin{equation*}
			T_0 = \mbb{\frac{\norm{\x_0-\x^*}^2 + (1-\lambda_2) \norm{\nabla F(\x^*)}^2}{n}}.
		\end{equation*}

		Noting that $\omega_0 = \frac{24\theta c^2\Vert U_R\Vert^2L^2}{n\mu^2 m (1-\sqrt{\lambda_2})} = \mo(1)$, it follows that
		\begin{equation}
			\label{eq:H0}
			H_0 = M_0 + \omega_0 T_0 = \mbb{\frac{\norm{\x_0-\x^*}^2 +(1-\lambda_2)  \norm{\nabla F(\x^*)}^2}{n}}.
		\end{equation}
		
		To estimate $q_0,q_2,q_3,q_4$, consider 
		\begin{subequations}
			\label{eq:p}
			\begin{align}
				p_2 &= \frac{\theta^2\bar\sigma^2}{n\mu^2} = \mo\left(\frac{1}{n}\right),\\
				p_3 &= \frac{24\lambda_2\theta^3\bar\sigma^2\Vert U_L\Vert^2\Vert U_R\Vert^2 L^2}{\mu^4(1-\sqrt{\lambda_2})} = \mbb{\frac{1}{1-\lambda_2}},\\
				p_5 &= \frac{96\theta^3\Vert U_L\Vert^2\Vert U_R\Vert^2L^2\Vert V^{-}\Vert^2\Vert\nabla F(\x^*)\Vert^2}{n\mu^4(1-\sqrt{\lambda_2})^2} = \mbb{\frac{\norm{\nabla F(\x^*)}^2}{n(1-\lambda_2)^3}}
			\end{align}
		\end{subequations}
		
		Combining \eqref{eq:p} and \eqref{eq:H0} leads to 
		\begin{align*}
			\hat H_1 &= \frac{1}{\theta-3}\left(p_2+\frac{p_5}{m^3}\right) = \mbb{\frac{1+\norm{\nabla F(\x^*)}^2}{n}},\\
			\hat H_2 &= m^2 H_0+\frac{2p_3}{2\theta-3} = \mbb{\frac{\norm{\x_0-\x^*}^2}{n(1-\lambda_2)^2}+\frac{\norm{\nabla F(\x^*)}^2}{n(1-\lambda_2)}+\frac{1}{1-\lambda_2}}.
		\end{align*}
		
		Therefore,
		\begin{equation*}
			\begin{aligned}
				c_0 &= \frac{3\theta L^2 c^2 \Vert U_R\Vert^2}{n\mu^2} = \mo(1),\quad 
				\frac{1}{1-q_0} = \frac{4}{1-\sqrt{\lambda_2}} = \mbb{\frac{1}{1-\lambda_2}},\\
				q_2 &= \frac{\lambda_2n\Vert U_L\Vert^2\bar\sigma^2\theta^2}{c^2\mu^2} = \mbb{1},\quad
				q_3 = \frac{4\lambda_2n \Vert U_L\Vert^2 \theta^2L^2\hat H_1}{\mu^2c^2(1-\lambda_2)} = \mbb{\frac{1+\norm{\nabla F(\x^*)}^2}{n(1-\lambda_2)}},\\
				q_4 &= \frac{4\lambda_2n \Vert U_L\Vert^2 \theta^2L^2\hat H_2}{\mu^2c^2(1-\lambda_2)}+\frac{4\Vert U_L\Vert^2\Vert V^{-}\Vert^2\Vert\nabla F(\x^*)\Vert^2\theta^2}{\mu^2c^2(1-\sqrt{\lambda_2})}\\
				&= \mbb{\frac{\norm{\x_0-\x^*}^2}{n(1-\lambda_2)^3}+\frac{\norm{\nabla F(\x^*)}^2}{n(1-\lambda_2)^2}+\frac{1}{(1-\lambda_2)^2}}.
			\end{aligned}
		\end{equation*}
	\end{proof}
	
	
	\section*{Acknowledgment}
	We would like to thank Bingyu Wang from The Chinese University of Hong Kong, Shenzhen for providing helpful feedback.
	
	\bibliographystyle{siamplain}
	\bibliography{references}
\end{document}


\maketitle

\section[Proof of Lemma]{Proof of \cref{lem:KL_KR}}
\label{sec:prf_KL_KR}
\begin{proof}
  $\overline{W}=\frac{I_n+W}{2}$ is symmetric, then $\overline{W}$ and $V^2=I_n-\overline{W}$ have the following spectral decomposition,
  \begin{equation*}
      \overline{W} = Q\overline{\Lambda}Q^{\T}, \quad V = Q\sqrt{I_n-\overline{\Lambda}}Q^{\T}
  \end{equation*}
  where $Q$ is an orthogonal matrix, 
  $$\overline{\Lambda}=\mathrm{diag}\left(\lambda_1(\overline{W}), \lambda_2(\overline{W}),..., \lambda_n(\overline{W})\right)$$, and 
  $$\sqrt{I_n-\overline{\Lambda}}=\mathrm{diag}\left(\sqrt{1-\lambda_1(\overline{W})}, \sqrt{1-\lambda_2(\overline{W})}, ..., \sqrt{1-\lambda_n(\overline{W})}\right)$$. 
  Therefore, $B$ can be decomposed as the following,
  \begin{equation}
      \begin{aligned}
          \left(
              \begin{array}{cc}
                  \overline{W} & -V\\
                  V\overline{W} & \overline{W}
              \end{array}
          \right) 
          &= 
          \left(
              \begin{array}{cc}
                  Q & \\
                    & Q
              \end{array}
          \right)
          \left(
              \begin{array}{cc}
                   \overline{\Lambda} & -\sqrt{I_n-\overline{\Lambda}}\\
                   \sqrt{I_n-\overline{\Lambda}}\overline{\Lambda} & \overline{\Lambda}
              \end{array}
          \right)
          \left(
              \begin{array}{cc}
                  Q^{\T} & \\
                    & Q^{\T}
              \end{array}
          \right)
      \end{aligned}
  \end{equation}
  With permutation operator $\Pi_1\in\R^{2n\times 2n}$, we have, 
  \begin{equation}
      \left(
          \begin{array}{cc}
                  \overline{\Lambda} & -\sqrt{I_n-\overline{\Lambda}}\\
                  \sqrt{I_n-\overline{\Lambda}}\overline{\Lambda} & \overline{\Lambda}
          \end{array}
      \right) = 
      \Pi_1 \mathrm{diag}(E_1, E_2, ..., E_n)\Pi_1 := \Pi_1 D \Pi_1
  \end{equation}
  where 
  \begin{equation}
      E_1 = I_2, 
      E_2 = \left(
          \begin{array}{cc}
              \lambda_k(\overline{W}) & -\sqrt{1-\lambda_k(\overline{W})}\\
              \lambda_k(\overline{W})\sqrt{1-\lambda_k(\overline{W})} & \lambda_k(\overline{W})
          \end{array}
      \right), k=2, 3, ..., n
  \end{equation}
  Next we do a spectral decomposition on $E_k$. Suppose $d_k$ is an eigenvalue of $E_k$, the characteristic polynomial of $E_k$ is, 
  \begin{equation}
      d_k^2 -2\lambda_k(\overline{W})d_k + \lambda_k(\overline{W})
  \end{equation}
  Since $\lambda_(\overline{W})<1$, when $k=2,3,...,n$, it holds that $4\lambda_k^2(\overline{W})-4\lambda_k(\overline{W})<0$. Therefore, $E_k$ can be diagnoalized in $\mathbb{C}$, i.e., 
  \begin{equation}\label{D_1}
      E_k = Z_k \left(
          \begin{array}{cc}
              \lambda_k(\overline{W})+i\sqrt{\lambda_k(\overline{W})-\lambda_k^2(\overline{W})} &0 \\
              0 & \lambda_k(\overline{W})-i\sqrt{\lambda_k(\overline{W})-\lambda_k^2(\overline{W})}
          \end{array}
      \right) Z_k^{-1}
  \end{equation}
  Where $i^2=-1$ and, 
  \begin{equation}
      \begin{aligned}
          Z_k &= \left(
              \begin{array}{cc}
                  1 & 1\\
                  -i\sqrt{\lambda_k(\overline{W})} & i\sqrt{\lambda_k(\overline{W})}
              \end{array}
          \right)\\
          Z_k^{-1} &= \left(
              \begin{array}{cc}
                  \frac{1}{2} & \frac{i}{2\sqrt{\lambda_k(\overline{W})}}\\
                  \frac{1}{2} & -\frac{i}{2\sqrt{\lambda_k(\overline{W})}}
              \end{array}
          \right)
      \end{aligned}
  \end{equation}
  Additionally, it can be verified that $\left(\mathbf{1}_n, \0\right)^{\T}$ and $\left(\0, \mathbf{1}_n\right)^{\T}$ are two eigenvectors of $B$. Thus after a permutation $\Pi_2$, we reach a decomposition for $B$ in lemma \ref{lem:B_decomposition}. 
  \begin{equation}
      K = \left(
          \begin{array}{cc}
              Q & \\
               & Q
          \end{array}
      \right)\pi_1 Z \Pi_2
  \end{equation}
  Where $Z=\mathrm{diag}\left(I_2, Z_1, Z_2, ..., Z_n\right)$. 
  Thus, 
  \begin{equation}
      \begin{aligned}
          \Vert K\Vert_2 &\leq \left\Vert \left(
              \begin{array}{cc}
                  Q & \\
                  & Q
              \end{array}
          \right)\right\Vert_2 \Vert\Pi_1\Vert_2 \Vert Z\Vert_2 \Vert\Pi_2\Vert_2
          = \Vert Z\Vert_2 
      \end{aligned}
  \end{equation}
  Next we compute $\Vert Z\Vert_2$. 
  \begin{equation}
      \begin{aligned}
          \Vert Z_k\Vert_2&= \sqrt{\lambda_{\max}(Z_kZ_k^*)} 
          = \sqrt{\lambda_{\max}\left(
              \begin{array}{cc}
              2 &0 \\
              0 & 2\lambda_k(\overline{W})
              \end{array}
              \right)}
          = \sqrt{2}
      \end{aligned}
  \end{equation}
  Therefore $\Vert Z\Vert_2 = \sqrt{2}$. 
  Similarly, 
  \begin{equation}
      \Vert K^{-1}\Vert_2 \leq \Vert Z^{-1}\Vert_2
  \end{equation}
  We first consider $\Vert Z_k^{-1}\Vert_2$.
  \begin{equation}
      Z_k^{-1} (Z_k^{-1})^* = 
      \frac{1}{4}
      \left(
          \begin{array}{cc}
              1+\frac{1}{\lambda_k(\overline{W})} & 1-\frac{1}{\lambda_k(\overline{W})}\\
              1-\frac{1}{\lambda_k(\overline{W})} & 1+\frac{1}{\lambda_k(\overline{W})}
          \end{array}
      \right)
  \end{equation}
  Thus $Z_k$ has eigenvalues $\frac{1}{2}$ and $\frac{1}{2\lambda_k(\overline{W})}$. Since $\lambda_k(\overline{W})<1$, 
  \begin{equation*}
      \Vert Z_k^{-1}\Vert = \sqrt{\frac{1}{2\lambda_k(\overline{W})}}
  \end{equation*}
  Hence $\Vert Z\Vert_2 = \sqrt{\frac{1}{2\lambda_n(\overline{W})}}= \frac{1}{1+\lambda_n}$.
  Therefore $\Vert K\Vert_2 \Vert K^{-1}\Vert_2\leq \sqrt{\frac{2}{1+\lambda_n}}<\sqrt{3}$, which reaches the desired result. 
\end{proof}

\section[Proof of Lemma]{Proof of \cref{lem:checkz}}\label{sec:prf_checkz}
\begin{proof}
    By squaring and taking conditional expectation on both sides of the recursion (\ref{edym}), we obtaion, 
    \begin{equation}\label{eq:checkz}
        \begin{aligned}
            \E\left[\Vert \check z_{k+1}\Vert^2|\mathcal{F}_{k}\right]&=\Vert D_1\check z_k+\frac{\alpha_k}{c}D_1K_{L,l}(\nabla F(\x_k)-\nabla F(\x^*))\Vert^2\\
            &\quad + \frac{\alpha_k^2\Vert D_1K_{L,l}\Vert^2}{c^2}\E\left[\Vert s_k(\x_k)\Vert^2|\mathcal{F}_k\right] + \frac{1}{c^2}\Vert K_{L,r}\Vert^2\Vert \y_{k+1}^*-\y_k^*\Vert^2 \\ 
            &\quad + 2\left\langle D_1\check z_k+\frac{\alpha_k}{c}D_1K_{L,l}(\nabla F(\x_k)-\nabla F(\x^*)), 
            \frac{1}{c}K_{L,r}(\y_{k+1}^*-\y_k^*)\right\rangle\\
            &\leq (1 + \gamma)\Vert D_1\check z_k+\frac{\alpha_k}{c}D_1K_{L,l}(\nabla F(\x_k)-\nabla F(\x^*))\Vert^2 \\
            &\quad + \frac{\alpha_k^2\lambda'\Vert K_L\Vert^2n\sigma^2}{c^2}
            + (1 + \frac{1}{\gamma})\frac{\Vert K_L\Vert^2}{c^2} \Vert \y_{k+1}^*-\y_k^*\Vert^2
        \end{aligned}
    \end{equation}
    Consider $\Vert F(\x_k)-F(\x^*)\Vert^2$, 
    \begin{equation}
        \begin{aligned}
            \Vert \nabla F(\x_k)-\nabla F(\x^*)\Vert^2 &= \sum_{i=1}^n\Vert \nabla f_i{x_{i,k}} - \nabla f_i(x_i^*)\Vert^2\\
            &\leq L^2\sum_{i=1}^n \Vert x_{i,k} - x^*\Vert^2\\ 
            &= L^2 \Vert \tx_k\Vert^2\\
            &\leq L^2(2n\Vert \bar z\Vert^2 + 2c^2\Vert K_R\Vert^2\Vert \check{z}_k\Vert^2)
        \end{aligned}
    \end{equation}
    Then, according to Jensen's Inequlity, $\forall t\in(0,1)$, 
    \begin{equation}
        \begin{aligned}
            &\Vert D_1\check z_k+\frac{\alpha_k}{c}D_1K_{L,l}(\nabla F(\x_k)-\nabla F(\x^*))\Vert^2\\
            &\leq 
            \frac{1}{t}\Vert D_1\check z_k\Vert^2 + \frac{1}{1-t}\Vert \frac{\alpha_k}{c}D_1K_{L,l}(\nabla F(\x_k)-\nabla F(\x^*))\Vert^2\\
            &\leq \frac{1}{t}\Vert D_1\Vert^2 \Vert\check z_k\Vert^2 + \frac{2n\alpha_k^2\Vert D_1\Vert^2\Vert K_L\Vert^2L^2}{c^2(1-t)}\Vert \bar z\Vert^2\\ 
            &\quad + \frac{2\alpha_k^2\Vert D_1\Vert^2\Vert K_L\Vert^2\Vert K_R\Vert^2L^2}{1-t}\Vert \check{z}_k\Vert^2
        \end{aligned}
    \end{equation}
    Let $\lambda' = \lambda_2(\overline{W})$, we can see from \eqref{D_1} that, 
    \begin{equation}
        \Vert D_1\Vert_2 = \sqrt{\lambda_2(\overline{W})}
    \end{equation}
    Choose $t = \sqrt{\lambda'}$, we derive, 
    \begin{equation}\label{eq:checkz_gamma}
        \begin{aligned}
            &\Vert D_1\check z_k+\frac{\alpha_k}{c}D_1K_{L,l}(\nabla F(\x_k)-\nabla F(\x^*))\Vert^2\\
            &\leq \left(\sqrt{\lambda'} + \frac{2\alpha_k^2\lambda'\Vert K_L\Vert^2\Vert K_R\Vert^2L^2}{1-\sqrt{\lambda'}}\right) \Vert \check{z}_k\Vert^2 
            + \frac{2n\alpha_k^2\lambda'\Vert K_L\Vert^2L^2}{c^2(1-\sqrt{\lambda'})}\Vert \bar z\Vert^2
        \end{aligned}
    \end{equation}
    Substitute \eqref{eq:checkz_gamma} and \eqref{eq:ystar} into \eqref{eq:checkz} and take the full expectation, 
    \begin{equation}
        \begin{aligned}
            &\E\left[\Vert \check{z}_{k+1}\Vert^2 \right]\\
            &\leq (1+\gamma) \left(\sqrt{\lambda'} + \frac{2\alpha_k^2\lambda'\Vert K_L\Vert^2\Vert K_R\Vert^2L^2}{1-\sqrt{\lambda'}}\right) \E\left[\Vert \check{z}_k\Vert^2\right] 
             +  \frac{\alpha_k^2\lambda'\Vert K_L\Vert^2n\sigma^2}{c^2} \\
            &\quad + (1+\gamma)\frac{2n\alpha_k^2\lambda'\Vert K_L\Vert^2L^2}{c^2(1-\sqrt{\lambda'})}\E\left[\Vert \bar z\Vert^2\right]\\
            &\quad + (1 + \frac{1}{\gamma})\frac{\Vert K_L\Vert^2}{c^2} \Vert V'^{-1}\Vert^2 \Vert \nabla F(\x^*)\Vert^2 \Vert \alpha_{k+1} - \alpha_k\Vert^2
        \end{aligned}
    \end{equation}
    Consider the matrix 
    \begin{equation}
        A_k := 
        \begin{pmatrix}
            (1-\frac{1}{2}\alpha_k\mu)  & \frac{2\alpha_k c^2\Vert K_R\Vert^2L^2}{n\mu}\\
            (1+\gamma)\frac{2\alpha_k^2 \lambda' \Vert K_L\Vert^2L^2 n}{c^2(1-\sqrt{\lambda'})} & (1+\gamma)\lambda'(\frac{1}{\sqrt{\lambda'}}+
            \frac{2\alpha_k^2\Vert K_L\Vert^2\Vert K_R\Vert^2L^2}{1-\sqrt{\lambda'}})
        \end{pmatrix}
    \end{equation}
    Under the diminishing stepsizes policy $\alpha_k = \frac{\theta}{\mu(k+m)}$, we choose proper $\gamma$ and $K_1-m$ such that when $k\geq K_1-m$, the spectral radius of $A_k$, $\rho(A_k)<1$. 
    First we let $K_1\geq \frac{2\sqrt{\lambda'}\theta\Vert K_L\Vert\Vert K_R\Vert L}{\mu(1-\sqrt{\lambda'})}$, then the following holds
    \begin{equation}
        \lambda'(\frac{1}{\sqrt{\lambda'}}+ \frac{2\alpha_{K_1-m}^2\Vert K_L\Vert^2\Vert K_R\Vert^2L^2}{1-\sqrt{\lambda'}}) < \frac{\sqrt{\lambda'} + 1}{2}
    \end{equation}
    We further choose $\gamma$ such that, 
    \begin{equation}
        (1+\gamma)\frac{\sqrt{\lambda'}+1}{2}=\frac{\sqrt{\lambda'}+3}{4}
    \end{equation}
    Therefore, 
    \begin{equation*}
        \gamma = \frac{1-\sqrt{\lambda'}}{2\sqrt{\lambda'}+2}
    \end{equation*}
    Then $A_k$ becomes, 
    \begin{equation}
        A_k = 
        \begin{pmatrix}
            1 - \frac{\alpha_k\mu}{2} & \frac{2\alpha_k c^2\Vert K_R\Vert^2 L^2}{n\mu}\\
            \frac{(\sqrt{\lambda'} + 3)\alpha_k^2\lambda'\Vert K_L\Vert^2 L^2n}{c^2(1-\lambda')} & \frac{\sqrt{\lambda'} + 3}{4}
        \end{pmatrix}
    \end{equation}
    Thus,
    \begin{equation}
        \begin{aligned}
            \mathrm{det}(I_2-A_k) &= 
            \begin{pmatrix}
                \frac{\alpha_k\mu}{2} & -\frac{2\alpha_k c^2\Vert K_R\Vert^2 L^2}{n\mu}\\ 
                - \frac{(\sqrt{\lambda'} + 3)\alpha_k^2\lambda'\Vert K_L\Vert^2 L^2n}{c^2(1-\lambda')} & \frac{1-\sqrt{\lambda'}}{4}
            \end{pmatrix}\\
            &= \frac{(1-\sqrt{\lambda'})\alpha_k\mu}{8} - \frac{2(\sqrt{\lambda'}+3)\alpha_k^3\lambda'\Vert K_L\Vert^2\Vert K_R\Vert^2L^4}{\mu(1-\lambda')}>0
        \end{aligned}
    \end{equation}
    which leads to 
    \begin{equation*}
        \alpha_k^2 < \frac{(1-\sqrt{\lambda'})^2(1+\sqrt{\lambda'})\mu^2}{16(\sqrt{\lambda'}+3)\lambda'\Vert K_L\Vert^2\Vert K_R\Vert^2L^4}
    \end{equation*}
    When $k\geq K_1-m$, $\alpha_k = \frac{\theta}{\mu(k+m)}\leq \alpha_{K_1-m}$, we choose $K_1$ such that, 
    \begin{equation*}
        K_1\geq \sqrt{\frac{(3+\sqrt{\lambda'})\lambda'}{1+\sqrt{\lambda'}}} \frac{4\theta\Vert K_L\Vert \Vert K_R\Vert L^2}{\mu^2(1-\sqrt{\lambda'})}
    \end{equation*}
    In light of Lemma 5 in \cite{pu2020distributed}, when $k\geq K_1-m$, the spectral radisu of $A_k$ is less than one, i.e., $\rho(A_k)<1$. By such a choice of $K_1$ and $\gamma$, we also derive the recursion for $\E\left[\Vert \check{z}_k\Vert^2\right]$. 
\end{proof}

\section[Proof of Lemma]{Proof of \cref{lem:hat_W}}\label{sec:prf_hat_W}
\begin{proof}
    From \cref{barzk,lem:checkz}, we have 
    \begin{equation}
        \begin{aligned}
            \Omega(k+1) &\leq 
            \left[(1-\frac{1}{2}\alpha_k\mu)+\omega(k)\frac{\alpha_k^2\lambda'(3+\sqrt{\lambda'})L^2\Vert K_L\Vert^2n}{c^2(1-\lambda')}\right] M(k)\\
            &\quad +
            \left[\frac{2\alpha_k c^2 \Vert K_R\Vert^2 L^2}{n\mu}+\omega(k)\frac{3+\sqrt{\lambda'}}{4}\right] T(k)\\ 
            &\quad + 
            \left(\frac{1}{n} + \frac{\omega(k) \lambda' n\Vert K_L\Vert^2}{c^2}\right)\alpha_k^2\sigma^2\\
            &\quad + \omega(k)\frac{(\sqrt{\lambda'}+3)\Vert K_L\Vert^2\Vert V'^{-1}\Vert^2\Vert\nabla F(\x^*)\Vert^2}{c^2(1-\sqrt{\lambda'})}\Vert\alpha_{k+1}-\alpha_k\Vert^2
        \end{aligned}
    \end{equation}
    We choose $K_1$ so that the following inequalities hold for all $k\geq K_1-m$, 
    \begin{equation}\label{U_ineq}
        (1-\frac{3}{2}\alpha_k\mu)+\omega(k)\frac{\alpha_k^2\lambda'(3+\sqrt{\lambda'})\delta^2n}{1-\lambda'} \leq 1-\frac{4}{3}\alpha_k\mu
    \end{equation}
    \begin{equation}\label{V_ineq}
        \frac{3\alpha_k\Vert K_L\Vert^2\Vert K_R\Vert^2\delta^2}{n\mu}+\omega(k)\frac{3+\sqrt{\lambda'}}{4}\leq (1-\frac{4}{3}\alpha_k\mu)\omega(k)
    \end{equation}
    Let 
    \begin{equation}\label{K1_3}
        \frac{4}{3}\alpha_k\mu \frac{1-\sqrt{\lambda'}}{24},\quad i.e. K_1\geq \frac{32\theta}{1-\sqrt{\lambda'}}
    \end{equation}
    then, 
    \begin{equation*}
        1-\frac{4}{3}\alpha_k\mu-\frac{3+\sqrt{\lambda'}}{4}\geq \frac{1-\sqrt{\lambda'}}{4}-\frac{1-\sqrt{\lambda'}}{24}=\frac{5(1-\sqrt{\lambda'})}{24}
    \end{equation*}
    Thus for inequality \ref{V_ineq} to hold, it is sufficient that, 
    \begin{equation}\label{omegak_1}
        \omega(k)\geq \frac{72\alpha_k\Vert K_L\Vert^2\Vert K_R\Vert^2\delta^2}{5n\mu(1-\sqrt{\lambda'})}
    \end{equation} 
    For inequality \ref{U_ineq} to hold, 
    \begin{equation}\label{omegak_2}
        \omega(k)\leq \frac{(1-\lambda')\mu}{6\lambda'(3+\sqrt{\lambda'})\delta^2 n}\frac{1}{\alpha_k}
    \end{equation}
    We choose $\omega(k)$ to be, 
    \begin{equation}
        \omega(k):= \frac{72\alpha_k\Vert K_L\Vert^2\Vert K_R\Vert^2\delta^2}{5n\mu(1-\sqrt{\lambda'})}
    \end{equation}
    Next we verify that for $K_1$ defined in \ref{K1}, such an $\omega(k)$ satisfy the inequalities \ref{omegak_1} nad \ref{omegak_2}. 
    \begin{equation}
        \begin{aligned}
            \frac{(1-\lambda')\mu}{6\lambda'(3+\sqrt{\lambda'})\delta^2 n}\frac{1}{\alpha_k} &\geq 
            \frac{(1-\lambda')\mu}{6\lambda'(3+\sqrt{\lambda'})\delta^2 n}\frac{1}{\alpha_{K_1-m}}\\
            &\geq \frac{(1-\lambda')\mu}{6\lambda'(3+\sqrt{\lambda'})\delta^2 n} \frac{12\sqrt{15}\delta^2\Vert K_L\Vert\Vert K_R\Vert}{5\mu}\sqrt{\frac{\lambda'(3+\sqrt{\lambda'})}{1+\sqrt{\lambda'}}}\frac{1}{1-\sqrt{\lambda'}}\\
            &=\sqrt{\frac{1+\sqrt{\lambda'}}{\lambda'(3+\sqrt{\lambda'})}}\frac{2\sqrt{15}\Vert K_L\Vert \Vert K_R\Vert}{5n}:=a
        \end{aligned}
    \end{equation}
    In addition, 
    \begin{equation}
        \begin{aligned}
            \frac{72\alpha_k\Vert K_L\Vert^2\Vert K_R\Vert^2\delta^2}{5n\mu(1-\sqrt{\lambda'})} &\leq 
            \frac{72\alpha_{K_1-m}\Vert K_L\Vert^2\Vert K_R\Vert^2\delta^2}{5n\mu(1-\sqrt{\lambda'})}\\
            &= \frac{72\Vert K_L\Vert^2\Vert K_R\Vert^2\delta^2}{5n\mu(1-\sqrt{\lambda'})}\frac{5\mu}{12\sqrt{15}\Vert K_L\Vert \Vert K_R\Vert\delta^2}\sqrt{\frac{(1-\lambda')(1-\sqrt{\lambda'})}{\lambda'(3+\sqrt{\lambda'})}}\\
            &= \sqrt{\frac{1-\lambda'}{\lambda'(3+\sqrt{\lambda'})}}\frac{2\sqrt{15}\Vert K_L\Vert\Vert K_R\Vert}{5n}:=b
        \end{aligned}
    \end{equation}
    Since $a/b=\sqrt{\frac{1+\sqrt{\lambda'}}{1-\lambda'}}>1$, $\omega(k)$ defined in \ref{omegak} satisfies the inequalities \ref{omegak_1} and \ref{omegak_2}. 
    Hence, we have, 
    \begin{equation}
        \begin{aligned}
            \Omega(k+1)&\leq (1-\frac{4}{3}\alpha_k\mu)\Omega(k) + (\omega(k)\lambda'+1)\frac{\alpha_k^2\sigma^2}{n}\\
            &\quad + \omega(k)\frac{3+\sqrt{\lambda'}}{1-\sqrt{\lambda'}}\frac{\Vert K_L\Vert^2\Vert \overline{W}\Vert^2\Vert\nabla f(x^*)\Vert^2}{\Vert V\Vert^2}\Vert \alpha_{k+1}-\alpha_k\Vert^2\\
            &\leq (1-\frac{4\theta}{3(k+m)})\Omega(k) + \left[\frac{72\theta\Vert K_L\Vert^2\Vert K_R\Vert^2\delta^2}{5n\mu^2(1-\sqrt{\lambda'})(k+m)}\lambda' + 1\right]\frac{\theta^2\sigma^2}{n\mu^2(k+m)^2}\\
            &\quad +\frac{72\theta\Vert K_L\Vert^2\Vert K_R\Vert^2\delta^2}{5n\mu^2(1-\sqrt{\lambda'})(k+m)}\frac{3+\sqrt{\lambda'}}{1-\sqrt{\lambda'}}\frac{\Vert K_L\Vert^2\Vert \overline{W}\Vert^2\Vert\nabla f(x^*)\Vert^2}{\Vert V\Vert^2} \frac{\theta^2}{\mu^2(k+m)^4}\\
            &:= (1-\frac{4\theta}{3(k+m)})\Omega(k)+\frac{p_2}{(k+m)^2}+\frac{p_3}{(k+m)^3}+\frac{p_5}{(k+m)^5}
        \end{aligned}
    \end{equation}
    Where 
    \begin{equation}\label{p_form}
        \begin{aligned}
            p_2 &= \frac{\theta^2\sigma^2}{n\mu^2},\\
            p_3 &= \frac{72\theta^3\sigma^2\Vert K_L\Vert^2\Vert K_R\Vert^2\delta^2}{5n^2\mu^4(1-\sqrt{\lambda'})},\\
            p_5 &= \frac{72\theta^3\Vert K_L\Vert^2\Vert K_R\Vert^2\delta^2}{5n\mu^4}\frac{3+\sqrt{\lambda'}}{(1-\sqrt{\lambda'})^2}\frac{\Vert \overline{W}\Vert^2\Vert\nabla f(x^*)\Vert^2}{\Vert V\Vert^2}
        \end{aligned}
    \end{equation}
    Thus, 
    \begin{equation}
        \begin{aligned}
            \Omega(k)&\leq\prod_{t=K_1-m}^{k-1}(1-\frac{4\theta}{3(t+m)})\Omega(K_1-m)\\
            &\quad +\sum_{t=K_1-m}^{k-1}\left[\prod_{j=t+1}^{k-1}(1-\frac{4\theta}{3(j+m)})\right]\left(\frac{p_2}{(t+m)^2}+\frac{p_3}{(t+m)^3}+\frac{p_5}{(t+m)^5}\right)\\
            &\leq \frac{K_1^{\frac{4\theta}{3}}}{(k+m)^{\frac{4\theta}{3}}}\Omega(K_1-m)\\
            &\quad + \frac{1}{(k+m)^{\frac{4\theta}{3}}}\left[\sum_{t=K_1-m}^{k-1}\frac{(m+t+1)^{\frac{4\theta}{3}}p_2}{(t+m)^2}+\sum_{t=K_1-m}^{k-1}\frac{(m+t+1)^{\frac{4\theta}{3}}p_3}{(t+m)^3}+\sum_{t=K_1-m}^{k-1}\frac{(m+t+1)^{\frac{4\theta}{3}}p_5}{(t+m)^5}\right]\\
            &\leq \frac{K_1^{\frac{4\theta}{3}}}{(k+m)^{\frac{4\theta}{3}}}\Omega(K_1-m)\\
            &\quad + \frac{2}{(k+m)^{\frac{4\theta}{3}}}\left[\sum_{t=K_1-m}^{k-1}p_2(m+t)^{\frac{4\theta}{3}-2}+\sum_{t=K_1-m}^{k-1}p_3(m+t)^{\frac{4\theta}{3}-3}+\sum_{t=K_1-m}^{k-1}p_5(m+t)^{\frac{4\theta}{3}-5}\right]\\
        \end{aligned}
    \end{equation}
    Suppose $\theta>3$, then 
    \begin{equation}
        \sum_{t=K_1-m}^{k-1}(m+t)^{\frac{4\theta}{3}-5}\leq \int_{-1}^k(m+t)^{\frac{4\theta}{3}-5}\mathrm{dt}\leq \frac{3}{4\theta-12}(m+k)^{\frac{4\theta}{3}-4}
    \end{equation}
    Similarly, 
    \begin{equation}
        \begin{aligned}
        \sum_{t=K_1-m}^{k-1}(m+t)^{\frac{4\theta}{3}-3}&\leq \frac{3}{4\theta-6}(m+k)^{\frac{4\theta}{3}-2}\\
        \sum_{t=K_1-m}^{k-1}(m+t)^{\frac{4\theta}{3}-2}&\leq \frac{3}{4\theta-3}(m+k)^{\frac{4\theta}{3}-1}
        \end{aligned}
    \end{equation}
    Thus, 
    \begin{equation}
        \begin{aligned}
            \Omega(k)\leq \frac{K_1^{\frac{4\theta}{3}}}{(k+m)^{\frac{4\theta}{3}}}\Omega(K_1-m) + \frac{6p_2}{4\theta-3}\frac{1}{k+m} 
            +\frac{3p_3}{2\theta-3}\frac{1}{(m+k)^2}+\frac{3p_5}{2\theta-6}\frac{1}{(m+k)^4}
        \end{aligned}
    \end{equation}
    In addition, for all $k\geq K_1-m$, $\frac{K_1^{\frac{4\theta}{3}}}{(k+m)^{\frac{4\theta}{3}}}\leq \frac{K_1}{k+m}$ since $\frac{K_1}{k+m}\leq 1\quad \forall k\geq K_1-m$, we have, 
    \begin{equation}
        \begin{aligned}
            \Omega(k)\leq \frac{K_1}{k+m} \Omega(K_1-m) + \frac{6p_2}{4\theta-3}\frac{1}{k+m} 
            +\frac{3p_3}{2\theta-3}\frac{1}{(m+k)^2}+\frac{3p_5}{2\theta-6}\frac{1}{(m+k)^4}
        \end{aligned}
    \end{equation}
    \textcolor{red}{Suppose $\Omega(K_1-m)$ is finite}. For all $k\geq K_1-m$,
    \begin{equation}
        \begin{aligned}
            \Omega(k) &\leq \frac{1}{k+m}\left[K_1 \Omega(K_1-m)+\frac{3}{2\theta-6}\left(\frac{2p_2(2\theta-6)}{4\theta-3}+\frac{p_3(2\theta-6)}{(2\theta-3)(k+m)}+\frac{p_5}{(k+m)^3}\right)\right]\\
            &\leq \frac{1}{k+m}\left[K_1 \Omega(K_1-m) +\frac{3}{2\theta-6}\left(2p_2+\frac{p_3}{K_1}+\frac{p_5}{K_1^3}\right)\right]\\
            &:=\frac{\hat W}{k+m}
        \end{aligned}
    \end{equation}
    Since $\Omega(k)=M(k)+\omega(k)T(k)$, we have $M(k)\leq\frac{\hat W}{k+m}$.

    From lemma \ref{Echeckz}, we have, 
    \begin{equation}
        \begin{aligned}
            T(k+1)&\leq \frac{3+\sqrt{\lambda'}}{4}T(k)+\frac{\lambda'(3+\sqrt{\lambda'}\delta^2n\theta^2\hat{W})}{\mu(1-\lambda')}\frac{1}{(k+m)^3}+\frac{\theta^2\lambda'\sigma^2}{n\mu^2}\frac{1}{(k+m)^2}\\
            &\quad + \frac{(3+\sqrt{\lambda'})\theta^2\Vert K_L\Vert^2\Vert\overline{W}\Vert^2\Vert\nabla f(x^*)\Vert^2}{(1-\sqrt{\lambda'}\Vert V\Vert^2\mu^2)}\frac{1}{(k+m)^4}\\
            &:= q_0V(k)+\frac{q_2}{(k+m)^2}+\frac{q_3}{(k+m)^3}+\frac{q_4}{(k+m)^4}\\
            &\leq q_0^{k-(K_1-m)+1}T(K_1-m)+\sum_{t=K_1-m}^k q_0^{k-(t+1)+1}\left[\frac{q_2}{(m+t)^2}+\frac{q_3}{(m+t)^3}+\frac{q_4}{(m+t)^4}\right]\\
            &:= q_0^{k-(K_1-m)+1}T(K_1-m) + T_2(k+1) + T_3(K+1)+T_4(k+1)
        \end{aligned}
    \end{equation}
    Where 
    \begin{equation}\label{q_form}
        \begin{aligned}
            q_0&= \frac{3+\sqrt{\lambda'}}{4}\\
            q_2&= \frac{\theta^2\lambda'\sigma^2}{n\mu^2}\\
            q_3 &= \frac{\lambda'(3+\sqrt{\lambda'})\delta^2n\theta^2\hat{W}}{\mu(1-\lambda')}\\
            q_4 &= \frac{(3+\sqrt{\lambda'})\theta^2\Vert\overline{W}\Vert^2\Vert\nabla f(x^*)\Vert^2}{(1-\sqrt{\lambda'})\Vert V\Vert^2\mu^2}\\
            T_2(k) &= \sum_{t=K_1-m}^{k-1}\frac{q_0^{k-t-1}q_2}{(m+t)^2}\\
            T_3(k) &= \sum_{t=K_1-m}^{k-1}\frac{q_0^{k-t-1}q_3}{(m+t)^3}\\
            T_4(k) &= \sum_{t=K_1-m}^{k-1}\frac{q_0^{k-t-1}q_4}{(m+t)^4}
        \end{aligned}
    \end{equation}
    We also have, 
    \begin{equation*}
        T_i(k+1) = q_0\left[\sum_{t=K_1-m}^{k-1}q_0^{k-t-1}\frac{q_i}{(m+t)^i}+\frac{q_i}{q_0(m+k)^i}\right]= q_0T_i(k)+\frac{q_i}{(m+k)^i}, i=2,3,4
    \end{equation*}
    Note that $K_1\geq \frac{32\theta}{1-\sqrt{\lambda'}}$, for all $k\geq K_1-m$, we can verify,
    \begin{equation*}
        \left(1-\frac{1}{k+m+1}\right)^4-q_0\geq\left(1-\frac{1}{K_1+1}\right)^4-q_0\geq \frac{1-q_0}{2}=\frac{1-\sqrt{\lambda'}}{8}
    \end{equation*}
    By induction we obtain, 
    \begin{equation}
        T_i(k)\leq \frac{1}{(k+m)^i}\frac{q_i}{(1-\frac{1}{K_1+1})^i-q_0},i=2,3,4
    \end{equation}
    Hence, 
    \begin{equation}
        T(k)\leq q_0^{k-(K_1-m)}T(K_1-m) + \frac{8q_2}{(1-\sqrt{\lambda'})(k+m)^2}+\frac{8q_3}{(1-\sqrt{\lambda'})(k+m)^3}+\frac{8q_4}{(1-\sqrt{\lambda'})(k+m)^4}
    \end{equation}
    which leads to the result.
\end{proof}

\section[Proof of Thm]{Proof of \cref{thm:bigthm}}
\label{sec:proof}
                                                                                             
\lipsum[106-112]

\section{Additional experimental results}
\Cref{tab:foo} shows additional
supporting evidence. 

\begin{table}[htbp]
{\footnotesize
  \caption{Example table}  \label{tab:foo}
\begin{center}
  \begin{tabular}{|c|c|c|} \hline
   Species & \bf Mean & \bf Std.~Dev. \\ \hline
    1 & 3.4 & 1.2 \\
    2 & 5.4 & 0.6 \\ \hline
  \end{tabular}
\end{center}
}
\end{table}

\bibliographystyle{siamplain}
\bibliography{references}